\documentclass{amsart}

\usepackage{amssymb}
\usepackage[dvips]{epsfig}
\usepackage{graphicx}
\usepackage{color}

\newcommand{\R}{{\mathbb R}}

\newtheorem{theorem}{Theorem}[section]

\theoremstyle{definition}

\newtheorem{conjecture}[theorem]{Conjecture}

\newtheorem*{merci}{Acknowledgements}

\theoremstyle{remark}
\newtheorem{remark}{Remark}[section]
\def\R{{\mathbb R}}

\numberwithin{equation}{section}

\begin{document}
\title[Amick-Schonbek system in 2D]{Numerical study of the 
Amick-Schonbek system in 2D}

\author[C. Klein]{Christian Klein}
\address{Institut de Math\'ematiques de Bourgogne,  UMR 5584;\\
Institut Universitaire de France \\
Universit\'e de Bourgogne Europe, 9 avenue Alain Savary, 21078 Dijon
                Cedex, France} 
\email{Christian.Klein@u-bourgogne.fr}

\author[J.-C. Saut]{Jean-Claude Saut}
\address{Laboratoire de Math\' ematiques, UMR 8628,\\
Universit\' e Paris-Saclay et CNRS\\ 91405 Orsay, France}
\email{jean-claude.saut@universite-paris-saclay.fr}

\date{\today}

\maketitle
\begin{center}
\emph{This paper is dedicated to the memory of Vladimir E. Zakharov}

\end{center}

\begin{abstract}
A numerical study of the 2D Amick-Schonbek Boussinesq system is 
presented. Numerical evidence is given for the transverse stability 
of the 1D solitary waves that are line solitary waves of the 2D 
equations. It is shown that initial data not satisfying the 
non-cavitation condition can lead to the formation of a gradient 
catastrophe in finite time. The numerical propagation of localised 
smooth initial data does not lead to the formation of stable 
structures localised in both spatial directions. For initial data 
satisfying the non-cavitation condition, smooth solutions appear to 
exist for all times.

\end{abstract}

\section{Introduction}
This paper presents numerical simulations of a special case of 
two-dimensional Boussinesq systems for the propagation of long, weakly nonlinear waves in the so-called  Boussinesq regime. Those systems degenerate to the Korteweg-de Vries equation when restricted to one-dimensional, one way propagating waves. 
They  model  water waves on a flat bottom propagating in 
both directions in the aforementioned  regime (see \cite{BCL, BCS1, 
BCS2}). The rigorous derivation is obtained by expanding the 
Dirichlet to Neumann operator occurring in the Zakharov-Craig-Sulem 
formulation of the water waves system  \cite{Za2,CSS2} with respect to the small 
parameter $\epsilon$ (see below).

\subsection{Background}

More precisely we  will consider a particular case of the so-called 
abcd Boussinesq systems for surface water waves, see \cite{BCL, 
BCS1,BCS2}. Note that Boussinesq \cite{Bou2} was the first to derive  
particular Boussinesq systems, but not in the class of those studied here. 
We refer to \cite{Da,Da2,Da3} for details and the  history of 
hydrodynamics in the nineteenth century. The systems can be put into 
the form
\begin{equation}
    \label{abcd2}
    \left\lbrace
    \begin{array}{l}
    \eta_t+\nabla \cdot {\bf v}+\epsilon \nabla\cdot(\eta {\bf v})+\mu\lbrack a \nabla\cdot \Delta{\bf v}-b\Delta \eta_t\rbrack=0 \\
    {\bf v}_t+\nabla \eta+\epsilon \frac{1}{2}\nabla |{\bf v}|^2+\mu\lbrack c\nabla \Delta \eta-d\Delta {\bf v}_t\rbrack=0.
\end{array}\right.
    \end{equation}
Here $\eta= \eta(x,t), x\in \R^d, d=1,2, t\in \R$ is the elevation of the wave, ${\bf v}={\bf v}(x,t)$ is a measure of the horizontal velocity, $\mu$ and $\epsilon$ are the small parameters (shallowness and nonlinearity parameters respectively) defined as
$$\mu=\frac{h^2}{\lambda^2}, \quad \epsilon= \frac{\alpha}{h}$$
where  $\alpha$ is a typical amplitude of the wave, $h$ a typical depth and $\lambda$ a typical horizontal wavelength.

In the Boussinesq regime, $\epsilon$ and $\mu$ are supposed to be 
small and of same order, $\epsilon\sim\mu\ll1$.  For simplicity we 
will put $\epsilon=\mu,$  writing \eqref{abcd2} as 
\begin{equation}
    \label{abcd}
    \left\lbrace
    \begin{array}{l}
    \eta_t+\nabla \cdot {\bf v}+\epsilon \lbrack\nabla\cdot(\eta {\bf v})+a \nabla\cdot \Delta{\bf v}-b\Delta \eta_t\rbrack=0 \\
    {\bf v}_t+\nabla \eta+\epsilon\lbrack \frac{1}{2}\nabla |{\bf 
	v}|^2+c\nabla \Delta \eta-d\Delta {\bf v }_t\rbrack=0.
\end{array}\right.
    \end{equation}
Particular cases are formally derived in \cite{Broer, Ding, Per}, see also \cite{Ding}.
The coefficients (a, b, c, d) are restricted by the condition 
$$a+b+c+d=\frac{1}{3}-\tau,$$
where $\tau\geq 0$ is the surface tension coefficient.

%

All (well-posed) members of the abcd class provide the same 
approximation of the water wave system in the Boussinesq regime, with 
an error of order $O(\epsilon^2t)$  (see \cite{BCL}), and  their 
dispersion relations are similar in the long wave regime. 
Nevertheless they possess rather  different mathematical properties 
as nonlinear dispersive systems due to the different behavior of the 
dispersion relation at high frequencies. In fact the order of their dispersive part  can vary from $+3$ to $-1$, \cite{BCS1, BCS2}.  We refer to \cite{Mel} for a careful analysis of linear dispersive estimates for the abcd systems.

Two particular one-dimensional cases of the abcd 
systems have remarkable properties and were studied numerically in 
\cite{KSAS,KSKBK}. We focus here on the two-dimensional version of one of them 
This system 
is a particular case of a system derived by Peregrine in \cite{Per}, 
and it is often referred to as the Peregrine or classical Boussinesq 
system, but Schonbek and Amick were the first to recognize its 
remarkable mathematical properties. A variant with slowly varying 
bottom is derived in \cite{TW,W}. This particular case  corresponds to $a=b=c=0$, 
$d=1$ and thus
\begin{equation}
    \label{AS}
    \left\lbrace
    \begin{array}{l}
    \eta_t+v_x+\epsilon(\eta v)_x=0 \\
    v_t+\eta_x+\epsilon(vv_x-v_{xxt})=0.
\end{array}\right.
    \end{equation}
which can be seen as the BBM regularization of a (linearly ill-posed) 
system  obtained from the Zakharov-Craig-Sulem formulation of 
	the water wave system after expanding the Dirichlet-to-Neumann operator at first order in $\epsilon.$ 

	
	It turns out that the hyperbolic structure of the underlying 
	Saint-Venant (Shallow Water) system can be used for solving the 
	Cauchy problem of the one-dimensional Amick-Schonbek system. In 
	particular the entropy of the Saint-Venant system helps to derive 
	an a priori bound for the solutions of the Amick-Schonbek 
	system, leading to the global well-posedness of the Cauchy 
	problem for arbitrary large initial data satisfying the physical 
	non-cavitation condition, see \cite{Am, Sch} and \cite{MTZ}, this 
	last paper presenting the sharpest results, in particular the 
	existence of global weak solutions. Note that the large time 
	behavior of these solutions is not known. We refer to \cite{KSAS} 
	for numerical simulations leading to relevant conjectures.


 There have been many papers dealing with the numerical study of the one-dimensional 
 Amick-Schonbek system, see for instance \cite{An-Dou, An-Dou2, KSAS}. The 
 case of a non-flat bottom is considered in \cite{GAD} and the 
 periodic case in \cite{ADM3}. Solitary waves were constructed in 
 \cite{chen}, where their interaction was also studied. 
 
 The present paper is a continuation of a previous work \cite{KSAS} that considered the 1D case.
 We focus here on the two-dimensional case for which the existence (and asymptotic behavior) of global solutions is still an open issue both for small and large initial data and thus for which  numerical simulations are crucial to present relevant conjectures.

The two-dimensional Amick-Schonbek system writes (we have kept the small parameter $\epsilon$):
\begin{equation}
    \label{2D}
    \left\lbrace
    \begin{array}{l}
    \eta_t+\nabla \cdot {\bf v}+\epsilon \lbrack\nabla\cdot(\eta {\bf v})=0 \\
    {\bf v}_t+\nabla \eta+\epsilon\lbrack \frac{1}{2}\nabla |{\bf v}|^2-d\Delta {\bf v }_t\rbrack=0.
\end{array}\right.
    \end{equation}

We are not aware of a connection between \eqref{2D} and the 
hyperbolic structure the underlying  two-dimensional Saint-Venant 
system similar to the aforementioned results of the 1D case leading 
for instance to the global existence of weak solutions. The local 
well-posedness of the Cauchy problem is relatively standard, (see 
\cite{BCS2}  in the 1D case but the extension to the 2D case is 
straightforward). The  existence  on long times of order  $O(1/\epsilon)$ of solutions to the Cauchy problem for \eqref{2D} that is needed for the full rigorous justification of the system, see \cite{BCL}, was established in \cite{SWX}, Theorem 4.4, under the non-cavitation condition $1+\eta_0\geq H_0>0$ and in \cite{Bu} without this condition.
\begin{remark}
The long time existence for a two-dimensional Amick-Schonbek system in presence of a non-trivial bathymetry is proven in \cite{Meso}.
\end{remark}
On the other hand, as aforementioned, the {\it global} behavior of solutions to the 2D Amick-Schonbek system is unknown.
In particular the possibility of  finite-time blow-up for solutions with initial data satisfying the non-cavitation condition  is an open question. The aim of the present paper is to present numerical simulations leading to relevant conjectures,
on the qualitative properties of solutions of the Cauchy problem.

\subsection{Main results}
In this paper we first study for the 1D solitary waves the stability 
which was also considered in the 1D setting in \cite{KSAS}. Obviously 
these 1D solutions are infinitely extended (in $y$) solutions to the 
2D system (\ref{2D}), so-called line solitary waves. We consider 
numerically 2D perturbations of these line solitary waves and present 
evidence for the 
following 
\begin{conjecture}
	The line solitary waves of the Amick-Schonbek system are 
	asymptotically stable.
\end{conjecture}

Initial data not satisfying the non-cavitation condition were studied 
in 1D numerically in \cite{KSAS}. We expand on these results and 
consider initial data violating this condition in one point in the 2D 
case. This leads to the following
\begin{conjecture}
	Initial data violating the non-cavitation condition in one point 
	in 1D lead to a gradient catastrophe in finite time.
\end{conjecture}
Similar results appear to hold in 2D, but we do not have access to 
sufficient numerical resolution in this case to present as convincing 
evidence as in the 1D case. 

We study the time evolution of 2D localised initial data and show 
that it forms an annular structure which then leads to a similar 
dynamical behavior as initial data that come close to violating the 
non-cavitation condition. In particular we do not find lump solutions to this system as known 
from the Kadomtsev-Petviashvili (KP) I equation. It appears that solutions 
to the 2D system show  a similar behavior as solutions to the 
KP II equation. The Boussinesq equation appears to have a 
\emph{defocusing} effect. We get the following 
\begin{conjecture}
	Localised smooth initial data satisfying the non-cavitation 
	condition do not lead to stable structures localised 
	in both spatial dimension. The solutions exist for all times and are 
	simply dispersed.
\end{conjecture}

Furthermore we study the 
appearance of zones of rapid modulated oscillations in the vicinity 
of shocks of the solutions to the same initial data for the 
corresponding Saint-Venant system, so called \emph{dispersive shock 
waves} (DSWs).

    This paper is organized as follows. In section 2  we discuss the 
	transverse stability of line solitary waves.     Section 3 is 
	devoted to the behavior of solutions corresponding to initial 
	data violating the non-cavitation condition. In section 4 the 
	time evolution of 2D localised smooth initial data is studied. 
	Section 5  is focussed on dispersive shock waves arising in the 
	zero dispersion limit. We add some concluding remarks in section 
	6.

\section{Transverse stability of line solitary waves}
The solitary waves constructed numerically for the 1D AS system, see 
for instance \cite{KSAS} and references therein, are $y$-independent 
traveling wave solutions to the 2D AS system (\ref{2D}). These 
solitary waves infinitely extended in $y$-direction form so-called 
line solitary waves. In this section we will study the transverse 
stability of these waves, i.e., their stability under 
small perturbations with a non-trivial $y$-dependence. We will 
present numerical evidence for the transverse stability of the line 
solitary waves. 

\subsection{Numerical approach}
The numerical approach to be applied in this paper is essentially a 
2D version of the 1D code in \cite{KSAS}. For the spatial coordinates 
we use a discrete Fourier transform which is conveniently computed 
with a 2D fast Fourier transform (FFT). This means we approximate a 
setting in $\mathbb{R}^{2}$ by a setting in $\mathbb{T}^{2}$ where 
the torus is chosen to be sufficiently large in both directions that 
the effects of periodicity are not dominant. In other words the 
period $2\pi L_{x}$, $L_{x}>0$ is large enough that the reentering 
radiation is small compared to the effects we want to study. On the 
other hand the period $2\pi L_{y}$ in $y$-direction is small meaning 
that we essentially consider line solitary waves on $\mathbb{T}_{y}$, 
but we expect that this will not affect the result (for the cases we 
studied, no effect of the choice of $L_{y}$ was noted). 

In other words we use the standard discretisation $x_{n}=-\pi 
L_{x}+n2\pi L_{x}/N_{x}$, $n=1,\ldots,N_{x}$ and $y_{m}=-\pi 
L_{y}+n2\pi L_{y}/N_{y}$, $m=1,\ldots,N_{y}$. For the time 
integration of the resulting ODE system for (\ref{2D}) we use as in 
\cite{KSAS} the classical explicit 4th order Runge-Kutta method. The 
code is tested for the example of the line solitary waves that are 
propagated to essentially the  precision with which they were 
constructed in \cite{KSAS} if the initial data are unperturbed. 

\subsection{Perturbed line solitary wave with $c=2$}
The 1D solitary wave solutions $Q_{c}(x)$ exist for $|c|>1$, the amplitude of 
the waves tending to the trivial solution for $|c|\to1$. We consider here only 
positive values of $c$, i.e., waves traveling to the right. First we 
consider the case $c=2$ with $N_{x}=2^{12}$ Fourier modes for $x\in 
10[-\pi,\pi]$ and $N_{y}=2^{7}$ Fourier modes for $y\in 3[-\pi,\pi]$ 
with $N_{t}=10^{4}$ time steps for $t\leq20$. The first perturbation we 
will look at is of the form 
\begin{equation}
	\eta(x,y,0) = Q_{2}(x) \pm 0.3\exp(-x^{2}-y^{2}),\quad 
	v_{x}(x,y,0) = V_{2}(x),\quad v_{y}(x,y,0)=0.
	\label{c2gauss}.
\end{equation}
This corresponds to perturbations of the line solitary wave of the 
order of $5\%$. Note that we need some finite perturbations to 
numerically see an effect of the perturbation in finite time. 

\begin{figure}[htb!]
 \includegraphics[width=0.49\textwidth]{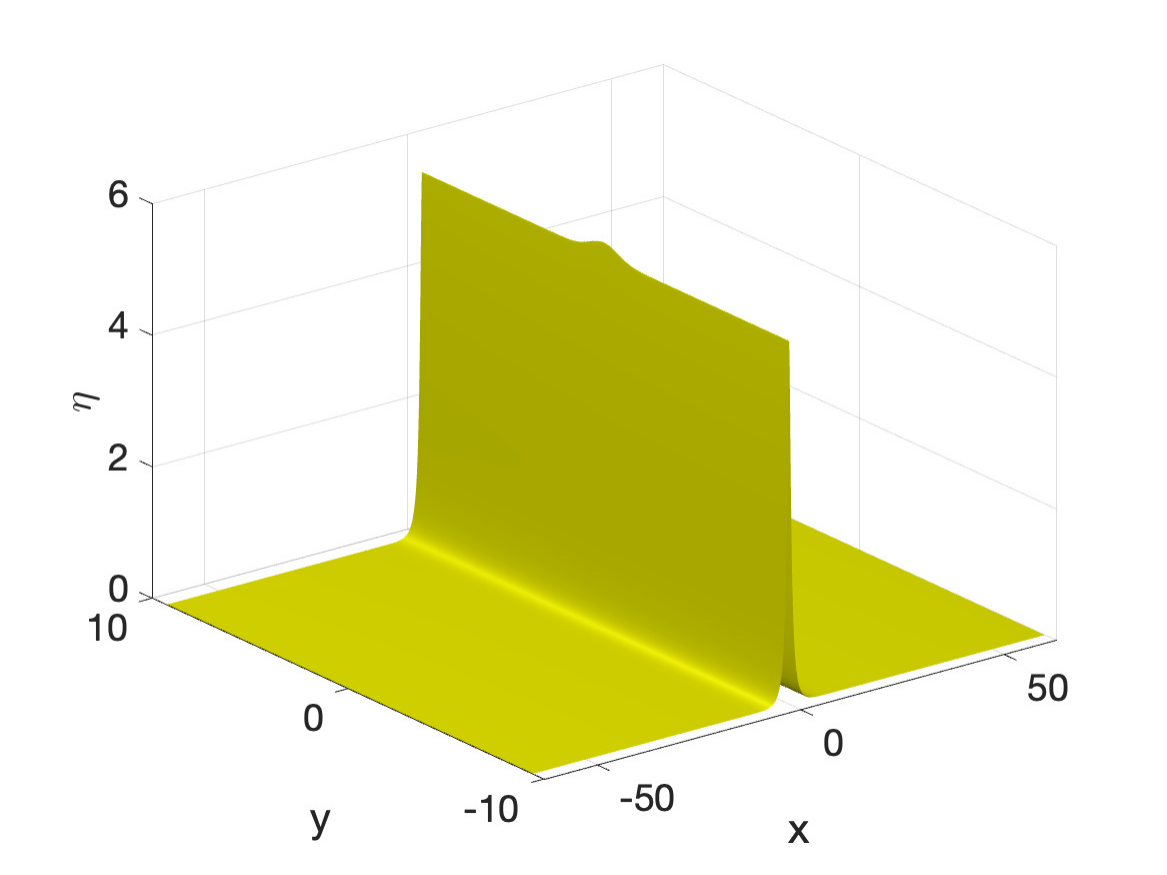}
 \includegraphics[width=0.49\textwidth]{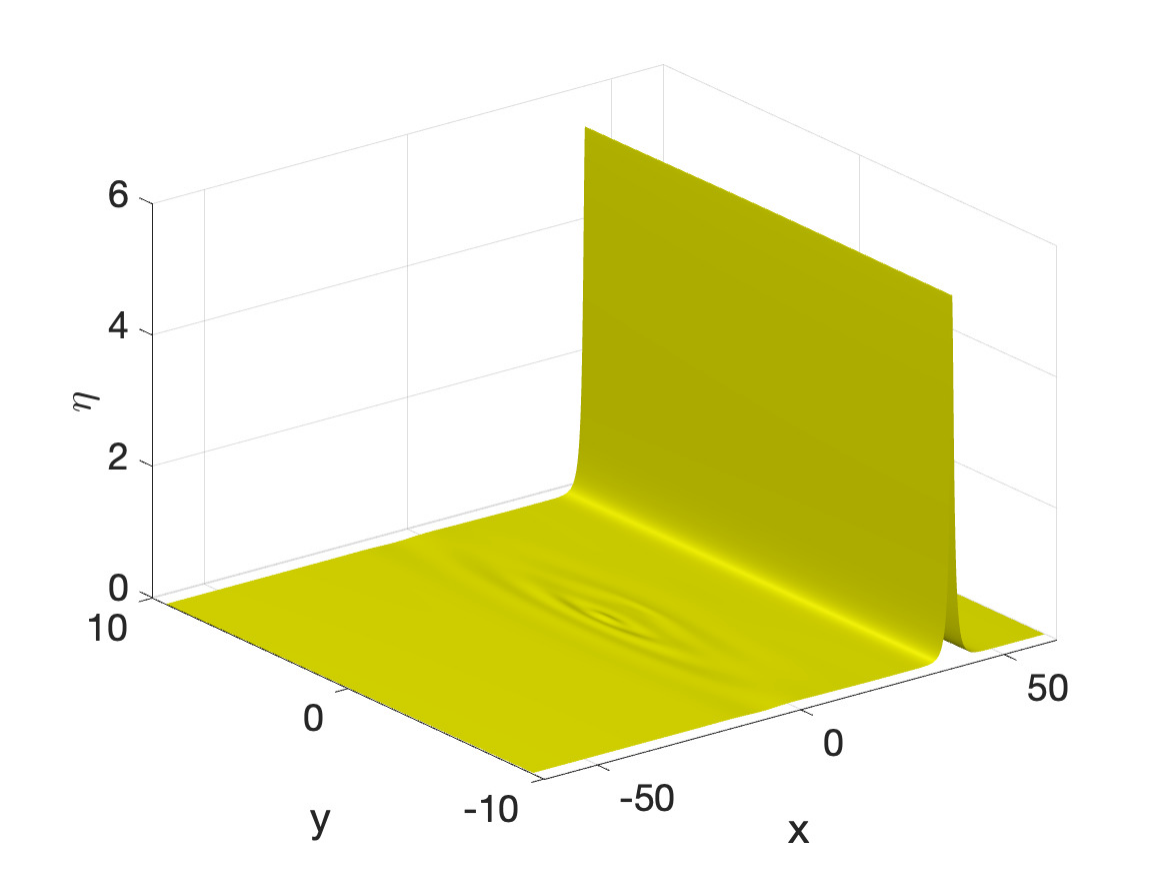}
 \caption{Solution to the AS system (\ref{2D}) for the initial data 
 (\ref{c2gauss}) for the $+$ sign, on the left $\eta$ for $t=0$, on the right for 
 $t=20$.}
 \label{figc2gauss}
\end{figure}

The solution $\eta$ for the initial data with the $+$ sign can be 
seen in Fig.~\ref{figc2gauss} for $t=0$ and $t=20$. It appears 
that the final state  corresponds to a line solitary wave 
with slightly larger velocity $c>2$ plus radiation. 

\begin{figure}[htb!]
 \includegraphics[width=0.49\textwidth]{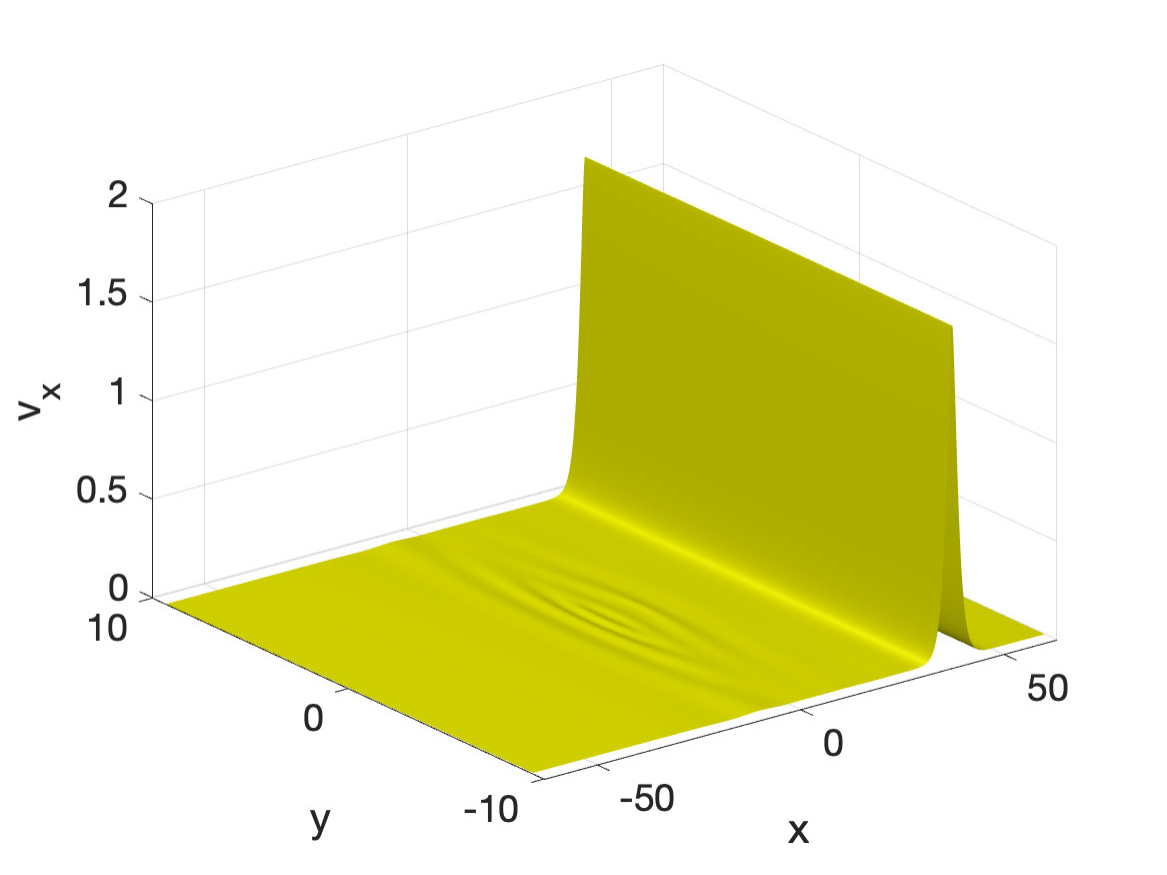}
 \includegraphics[width=0.49\textwidth]{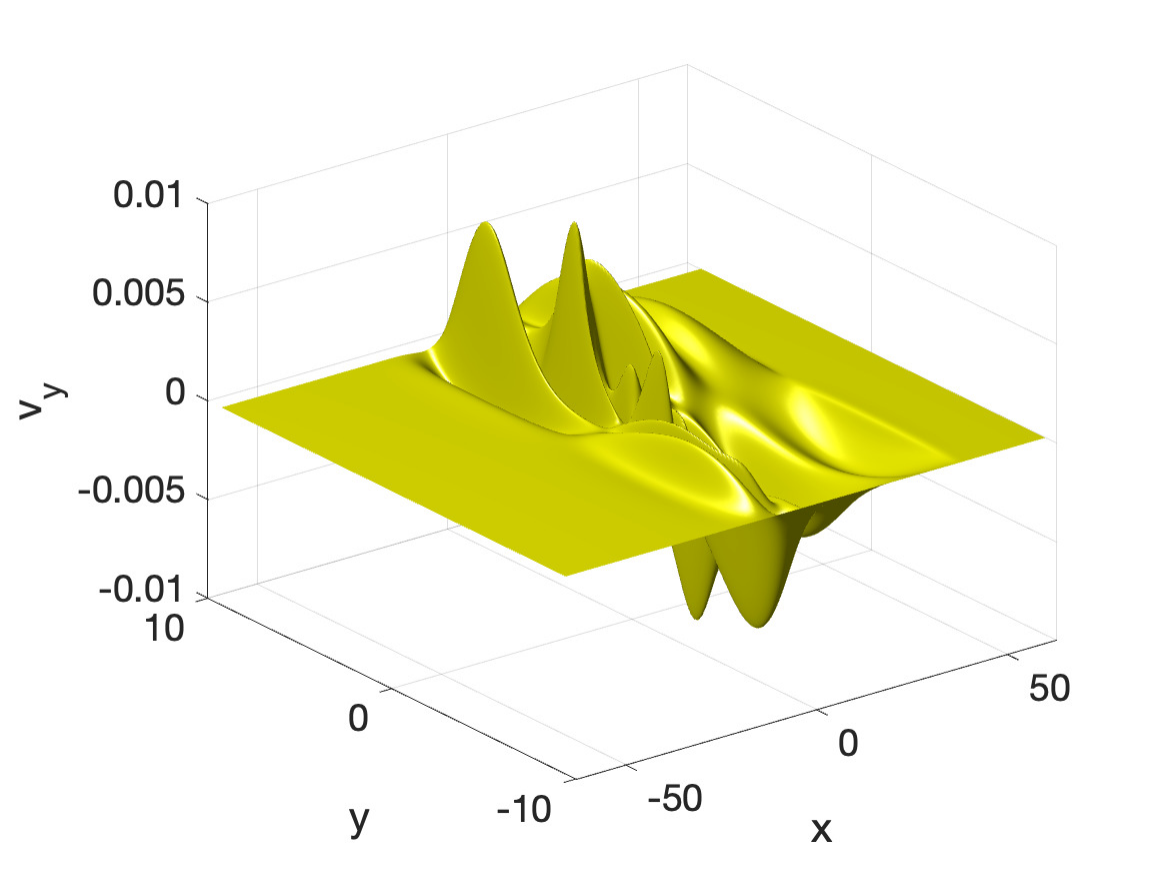}
 \caption{Solution to the AS system (\ref{2D}) for the initial data 
 (\ref{c2gauss}) for $t=20$, on the left $v_{x}$, on the right 
 $v_{y}$. }
 \label{figc2gaussv}
\end{figure}

This interpretation is confirmed by the velocities $v_{x}$ and $v_{y}$ 
which are shown for the same initial 
data for $t=20$ in Fig.~\ref{figc2gaussv}.

\begin{figure}[htb!]
 \includegraphics[width=0.49\textwidth]{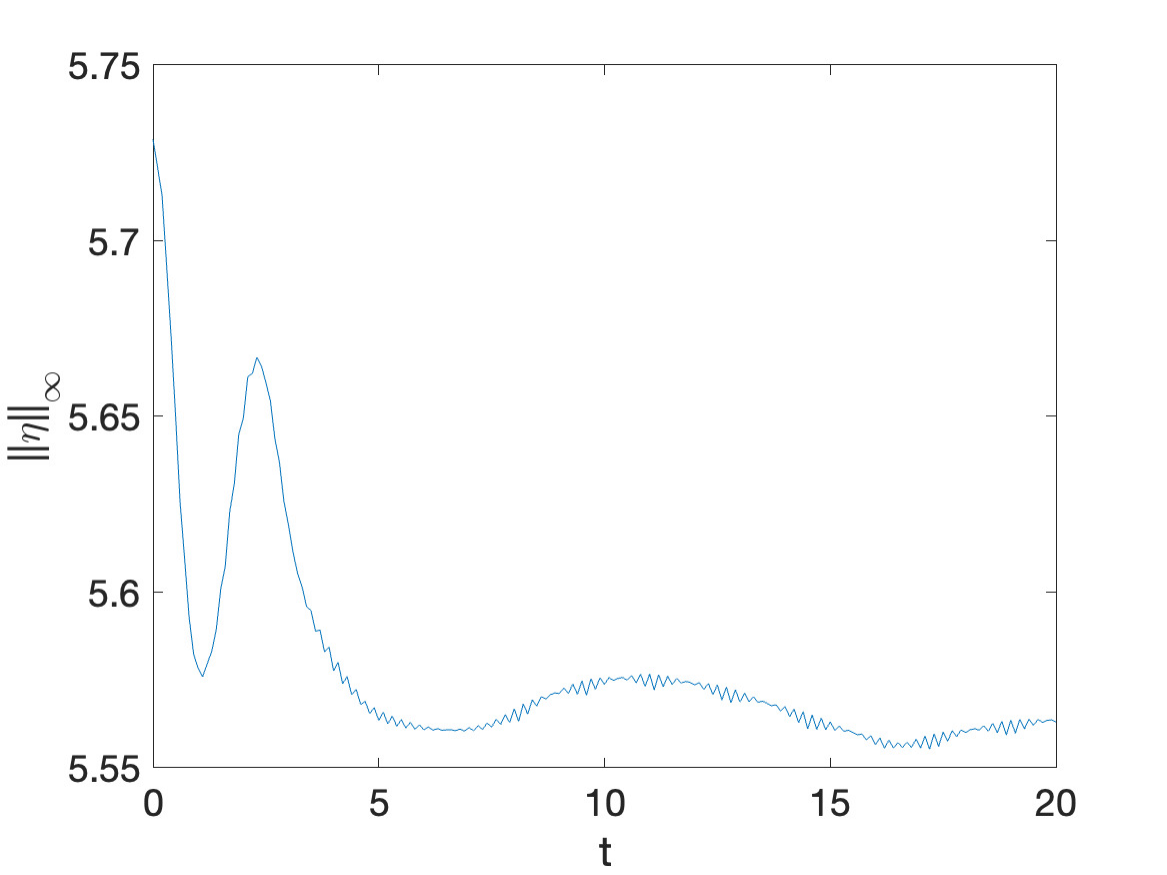}
 \includegraphics[width=0.49\textwidth]{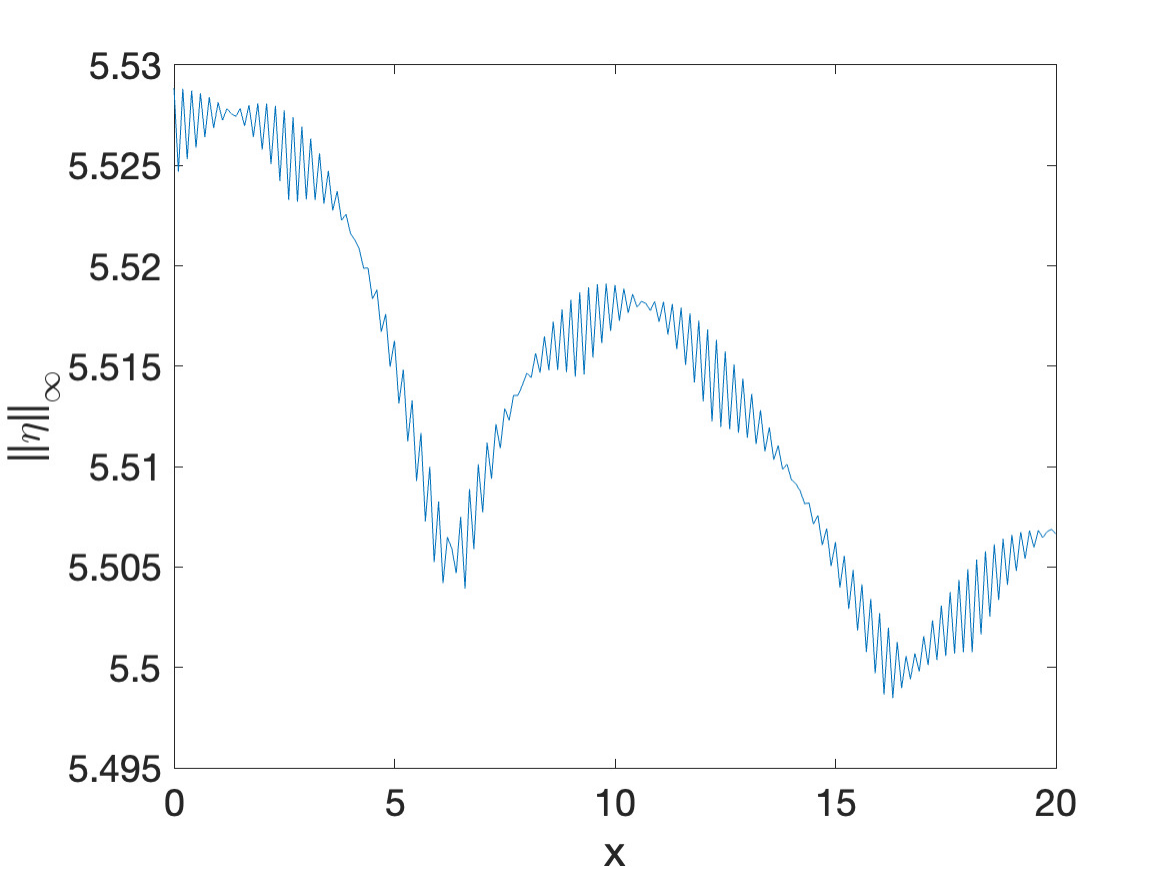}
 \caption{$L^{\infty}$ norm of the solution $\eta$ to the AS system (\ref{2D}) for the initial data 
 (\ref{c2gauss}), on the left for the $+$ sign, on the right for the 
 $-$ sign.}
 \label{figc2gaussmax}
\end{figure}

The $L^{\infty}$ norm of the solution $\eta$ in this case on the left 
of Fig.~\ref{figc2gaussmax} shows oscillations around what appears to 
be a constant final state, a line solitary wave. The above figures are very similar for the $-$sign 
in (\ref{c2gauss}). Therefore we only show the $L^{\infty}$ norm of 
$\eta$ in this case on the right of Fig.~\ref{figc2gaussmax}. The 
final state is reached more slowly, but it appears to be again a line 
solitary wave, this time with slightly smaller velocity than the 
unperturbed case with $c=2$. 

\begin{figure}[htb!]
 \includegraphics[width=0.49\textwidth]{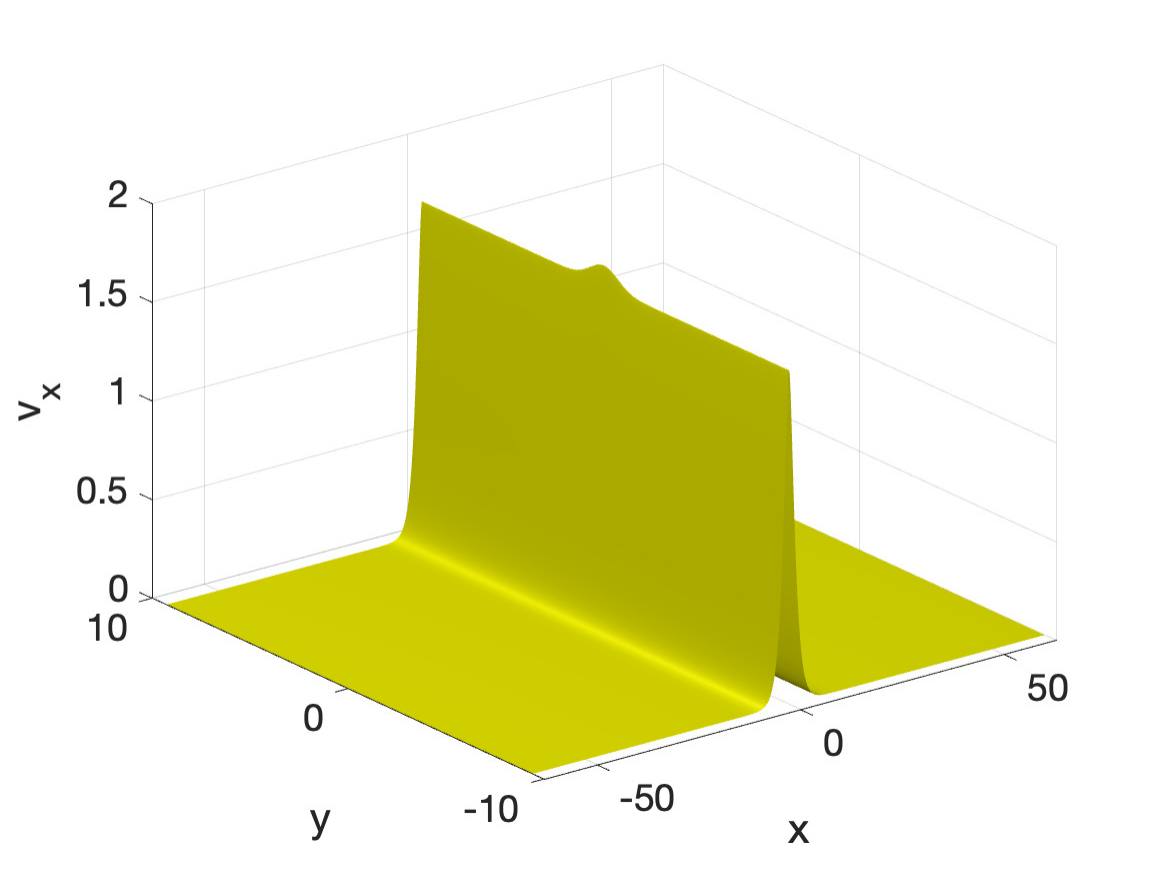}
 \includegraphics[width=0.49\textwidth]{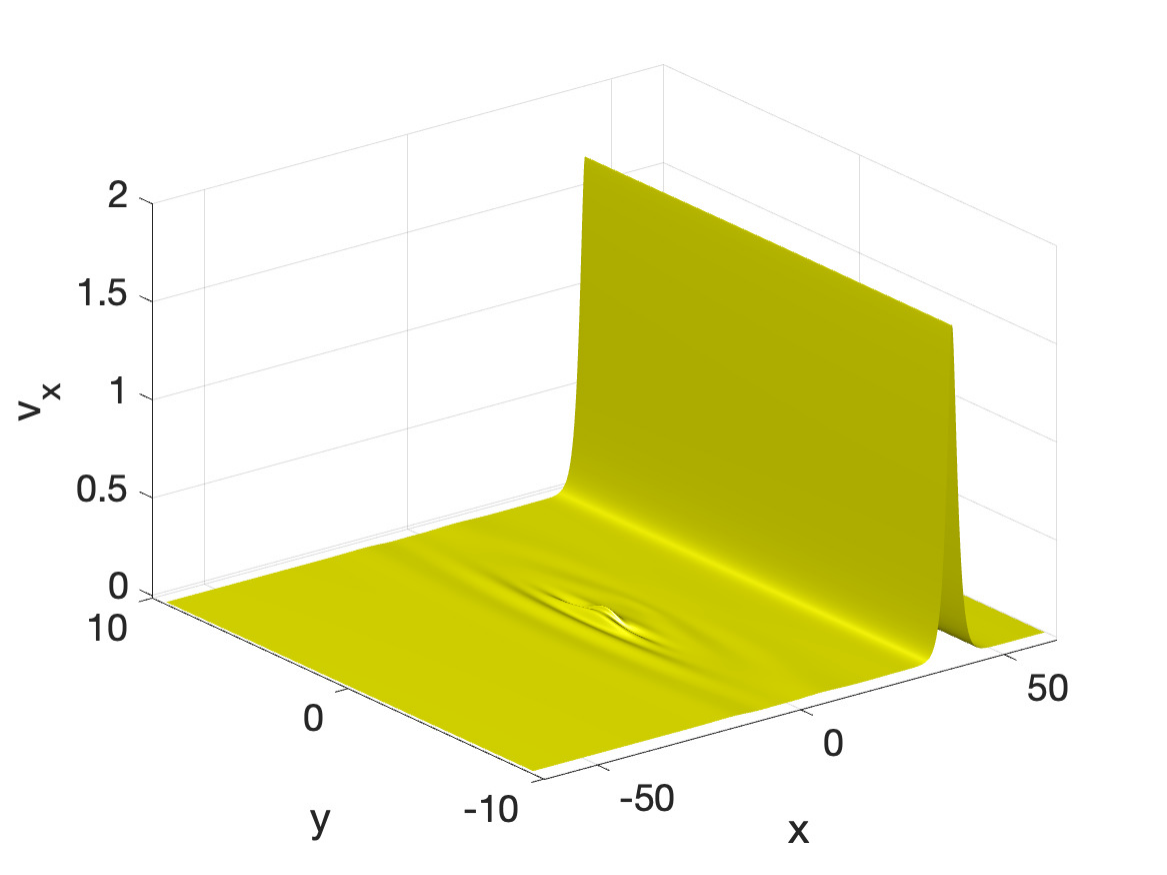}\\
 \includegraphics[width=0.49\textwidth]{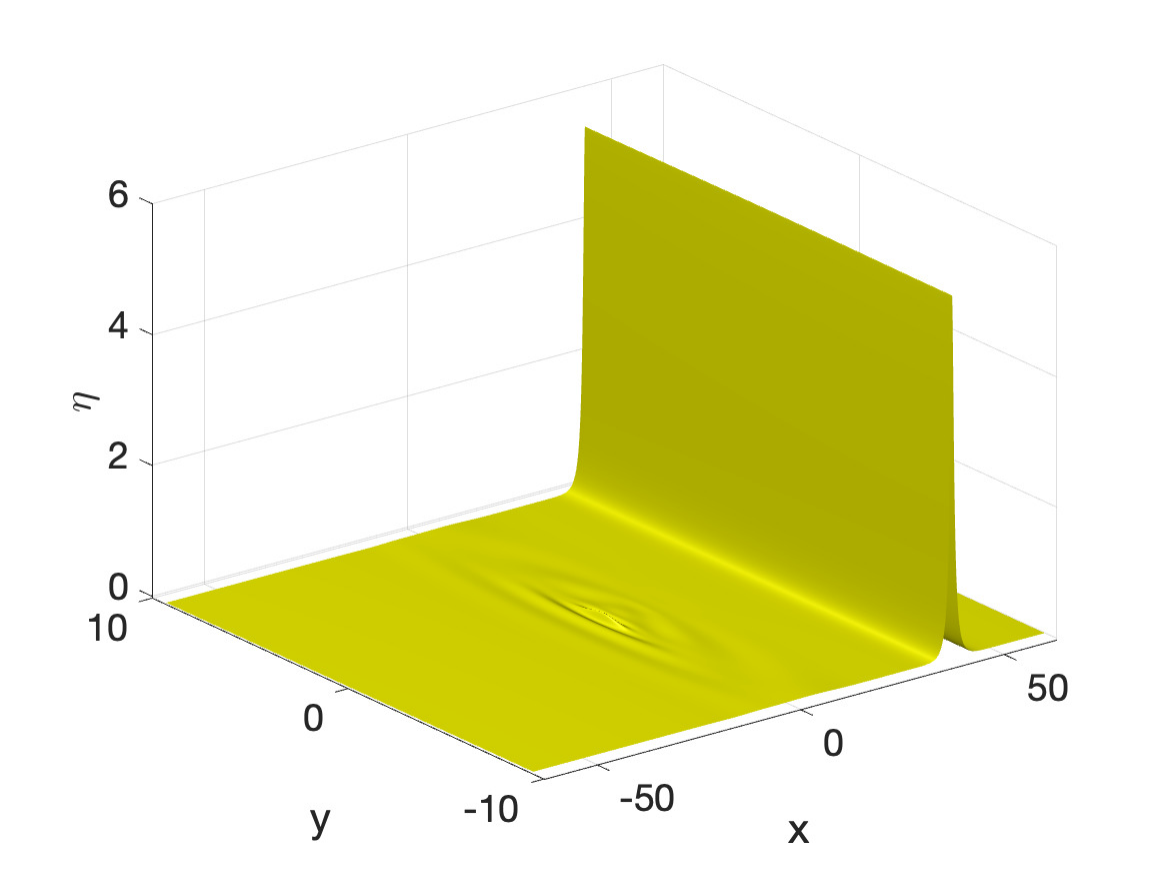}
 \includegraphics[width=0.49\textwidth]{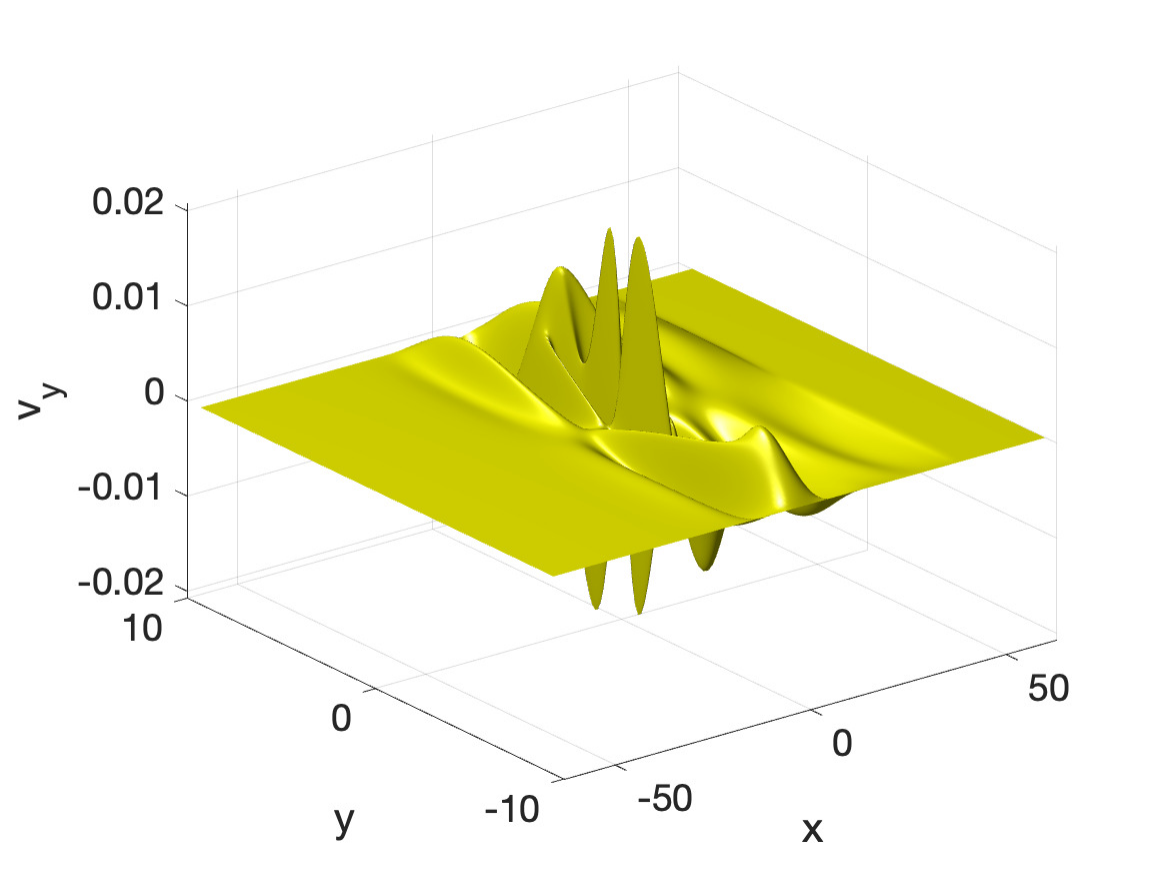}
 \caption{Solution to the AS system (\ref{2D}) for the initial data 
 (\ref{c2gaussvx}), in the upper row on the left $v_{x}$ for $t=0$, on the right for 
 $t=20$, in the lower row on the left $\eta$ and on the right 
 $v_{y}$, both for t=20.}
 \label{figc2gaussvx}
\end{figure}

If we perturb the velocity $v_{x}$ by a small Gaussian, i.e., if we 
consider initial data of the form 
\begin{equation}
	\eta(x,y,0) = Q_{2}(x) ,\quad 
	v_{x}(x,y,0) = V_{2}(x)+ 0.1\exp(-x^{2}-y^{2}),\quad 
	v_{y}(x,y,0)=0,
	\label{c2gaussvx}
\end{equation}
we get the solution shown in Fig.~\ref{figc2gaussvx}. Again the final 
state appears to be a line solitary wave plus radiation. 

\begin{figure}[htb!]
 \includegraphics[width=0.32\textwidth]{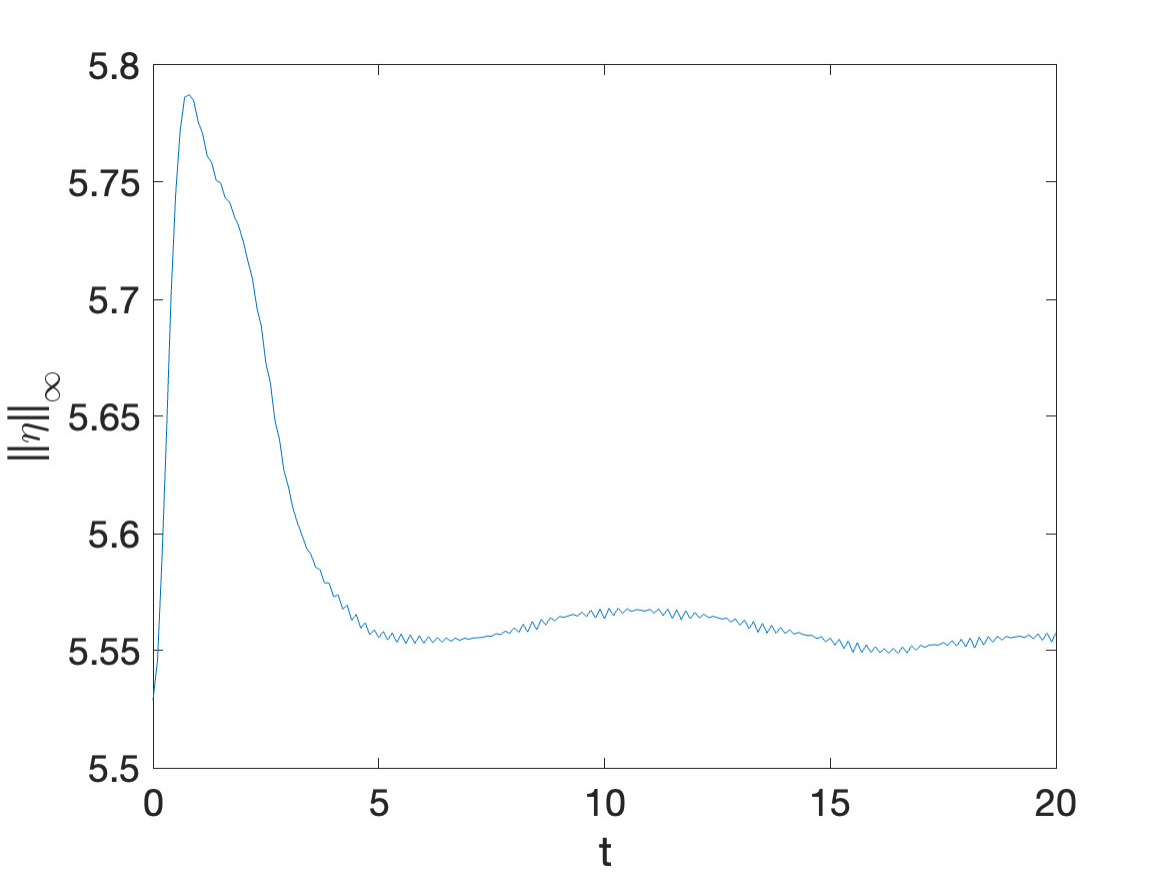}
 \includegraphics[width=0.32\textwidth]{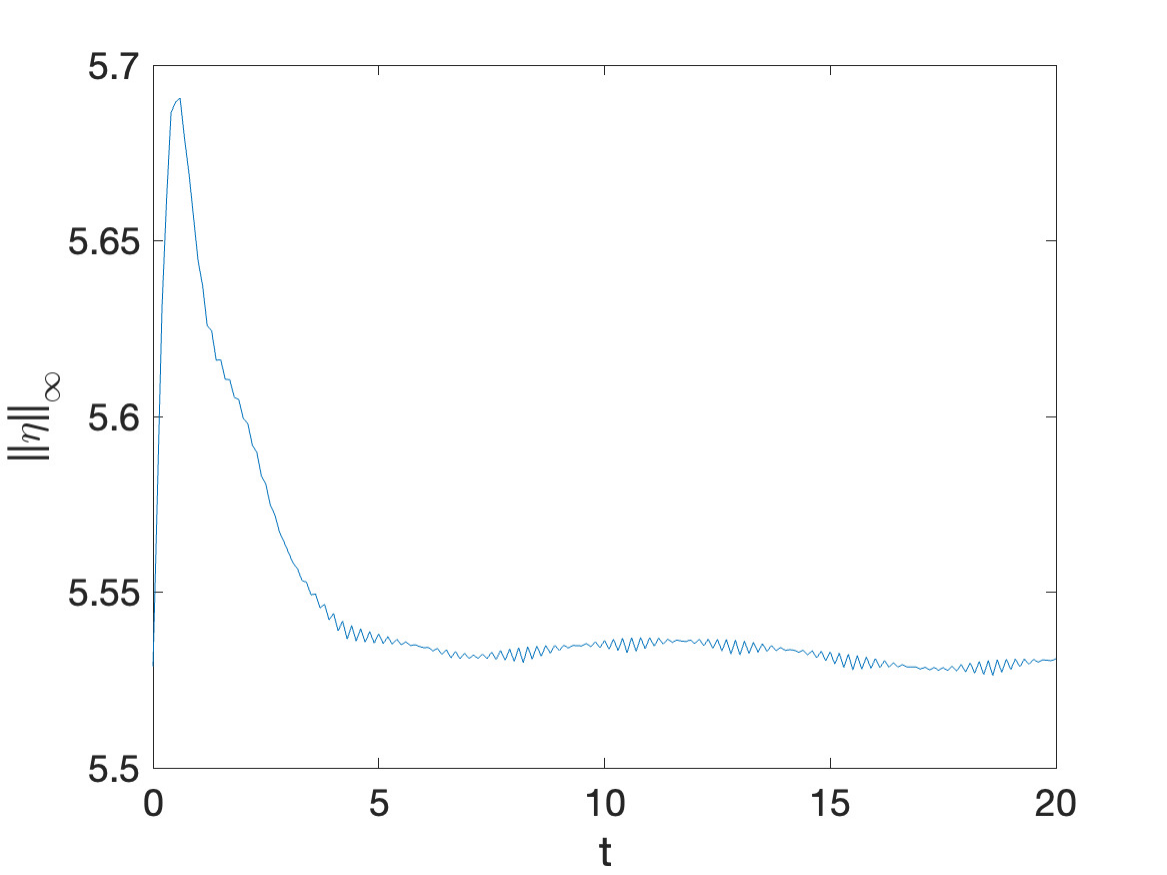}
  \includegraphics[width=0.32\textwidth]{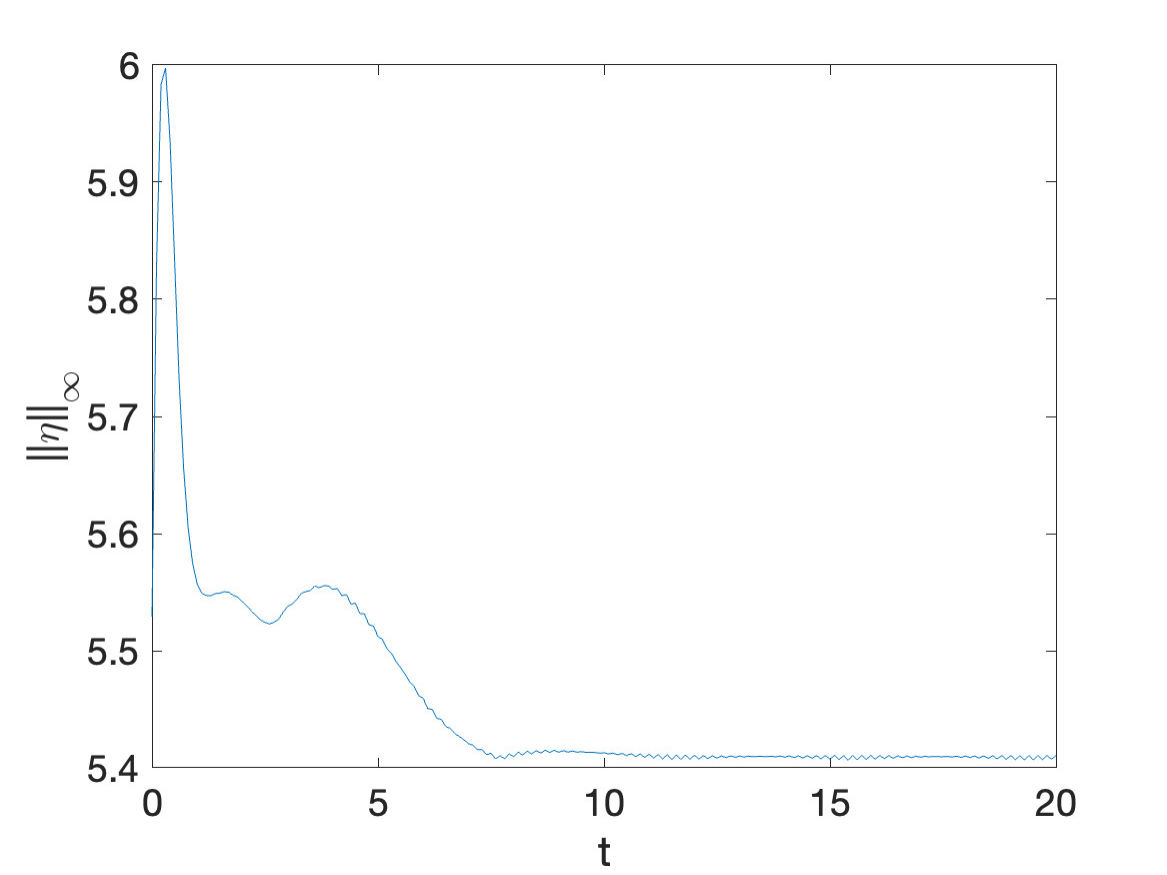}
 \caption{$L^{\infty}$ norms of solutions $\eta$ to the AS system 
 (\ref{2D}) for various initial data, on the left for the solution 
 shown in Fig.~\ref{figc2gaussvx}, in the middle for the solution 
 shown in Fig.~\ref{figc2gaussvy}, on the right for the solution 
 shown in Fig.~\ref{figc2cosv}. }
 \label{figc2max}
\end{figure}

The $L^{\infty}$ norm of the solution $\eta$ in 
Fig.~\ref{figc2gaussvx} can be seen on the left of 
Fig.~\ref{figc2max}. It appears to oscillate around a final state 
corresponding to a line solitary wave. 

\begin{figure}[htb!]
 \includegraphics[width=0.49\textwidth]{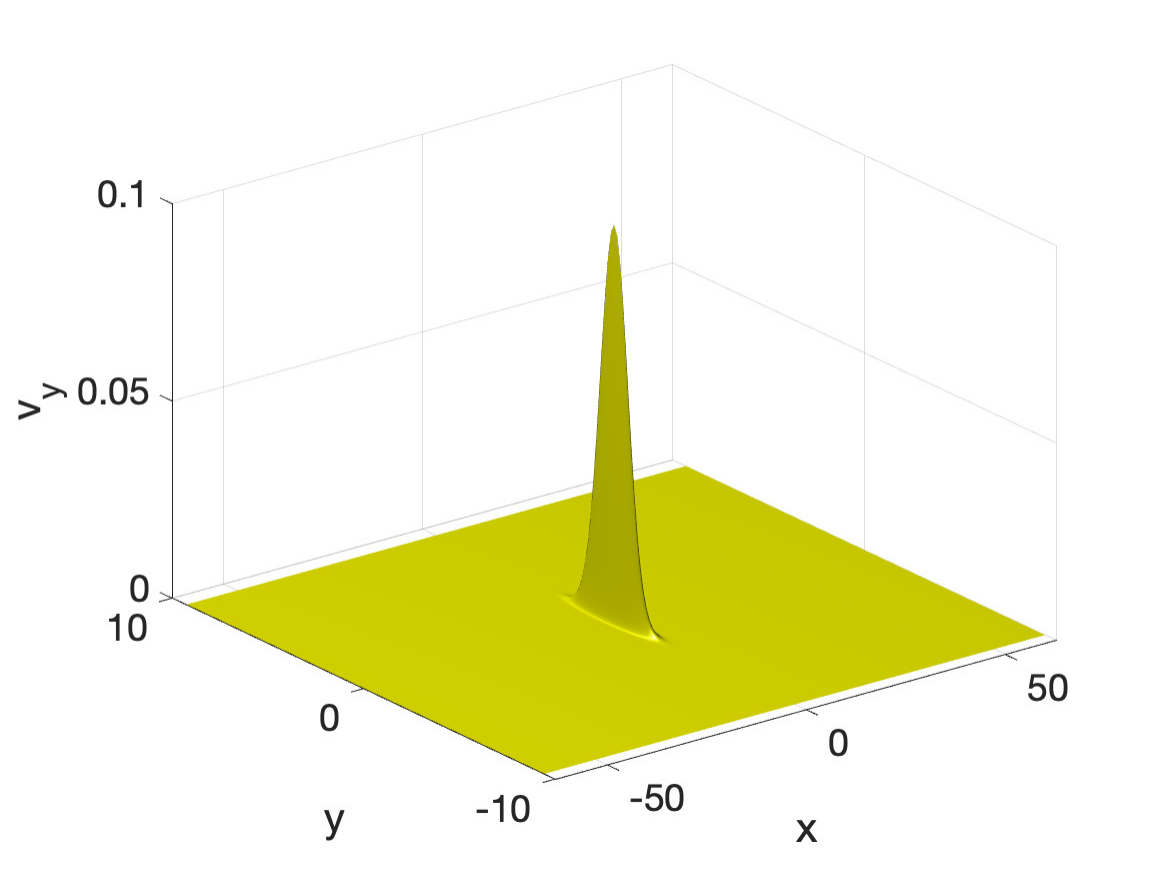}
 \includegraphics[width=0.49\textwidth]{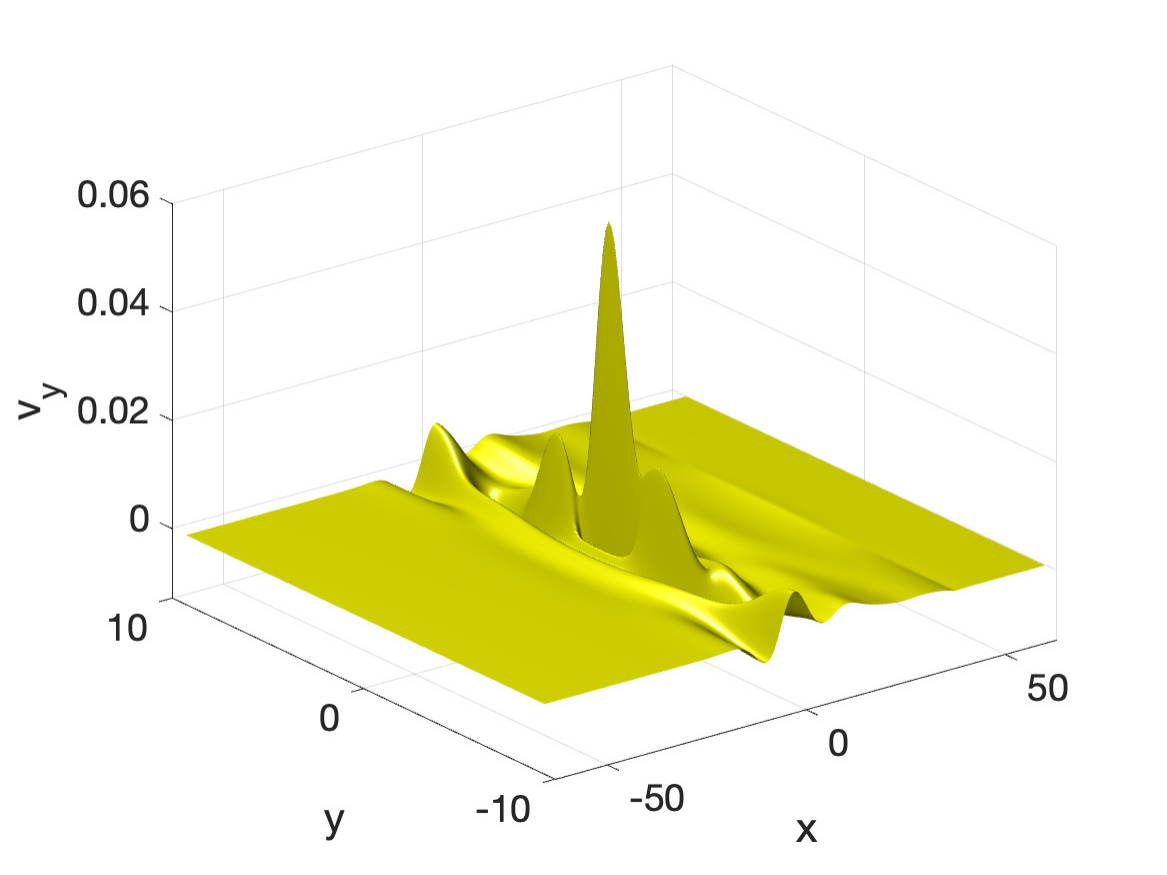}\\
 \includegraphics[width=0.49\textwidth]{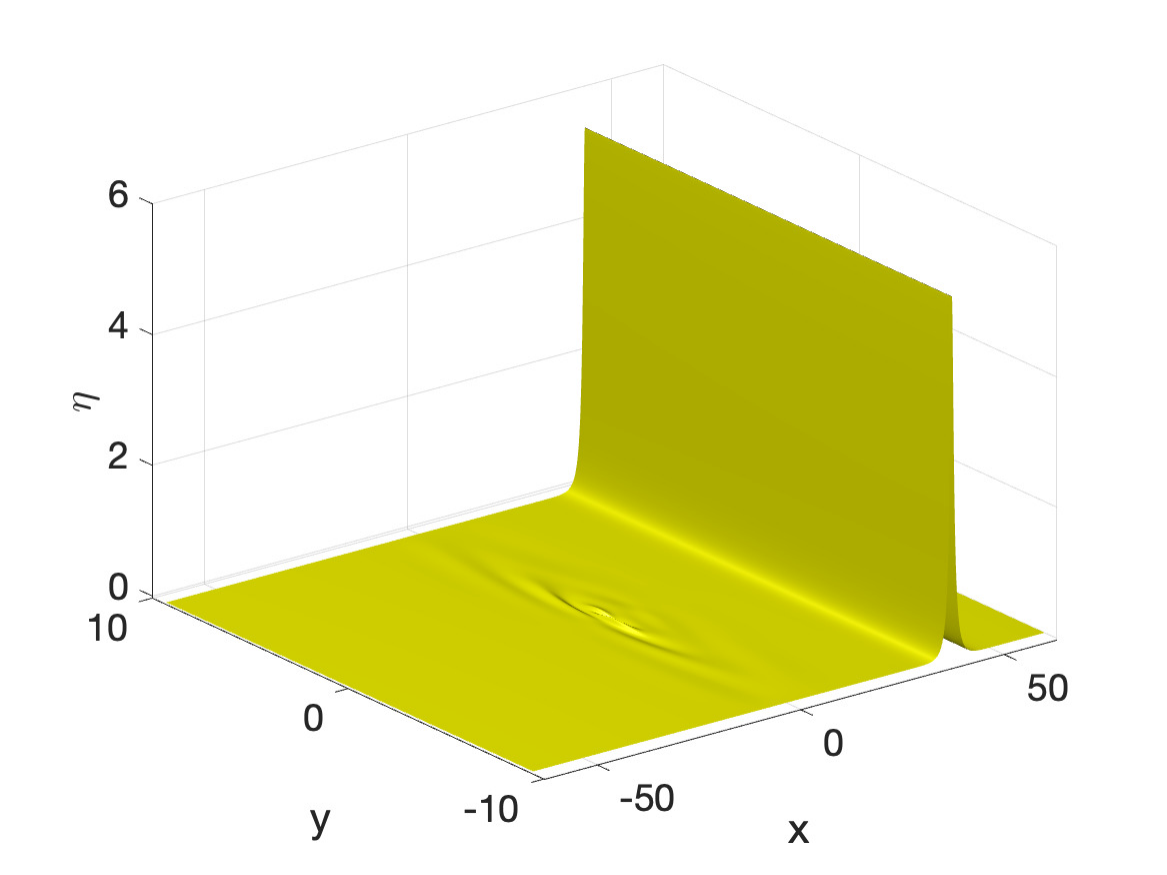}
 \includegraphics[width=0.49\textwidth]{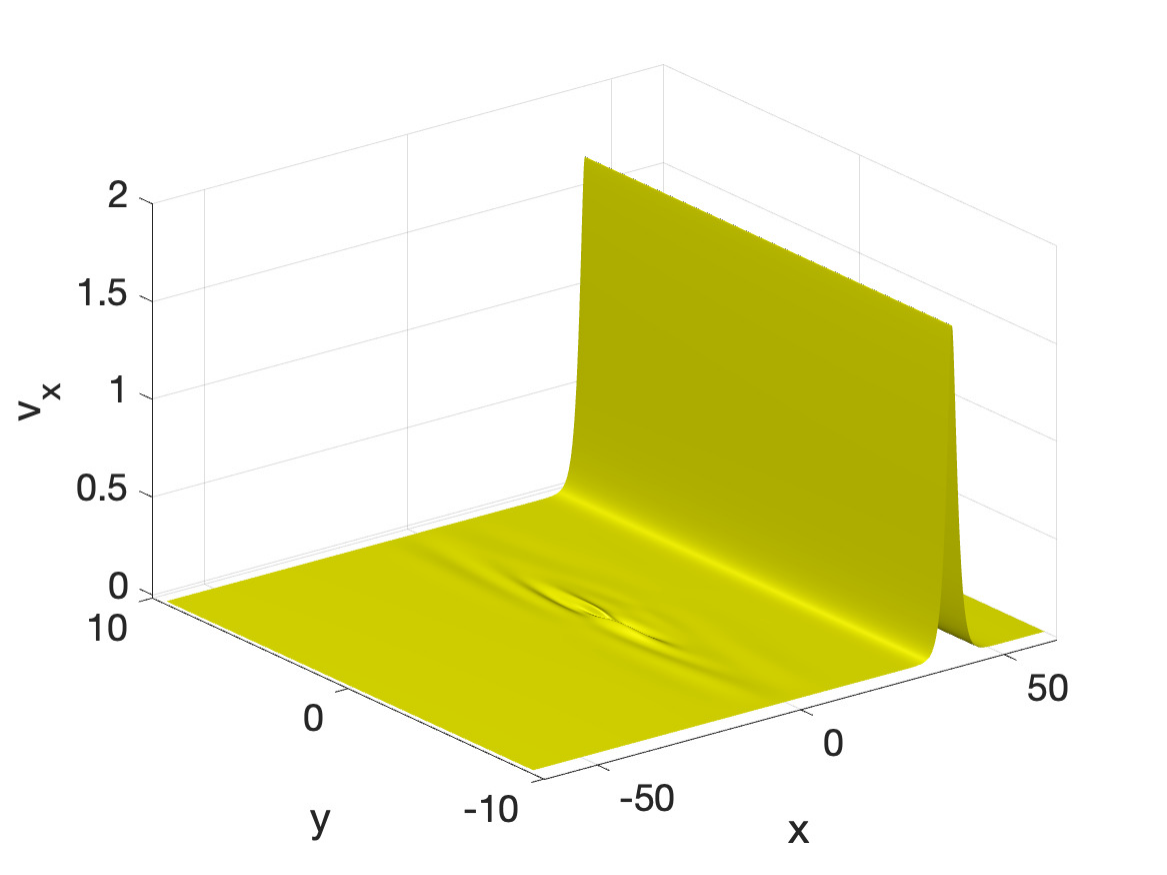}
 \caption{Solution to the AS system (\ref{2D}) for the initial data 
 (\ref{c2cos}), in the upper row on the left $v_{y}$ for $t=0$, on the right for 
 $t=20$, in the lower row on the left $\eta$ and on the right 
 $v_{x}$, both for t=20.}
 \label{figc2gaussvy}
\end{figure}

If we perturb $v_{y}$ with a small Gaussian, i.e., if we 
consider initial data of the form 
\begin{equation}
	\eta(x,y,0) = Q_{2}(x) ,\quad 
	v_{x}(x,y,0) = V_{2}(x),\quad v_{y}(x,y,0)=0.1\exp(-x^{2}-y^{2}),
	\label{c2gaussvy}
\end{equation}
we get the solution shown in Fig.~\ref{figc2gaussvy}. Once more the 
final state of the solution seems to be a line solitary wave plus 
radiation. The $L^{\infty}$ norm of the solution $\eta$ in this case 
can be seen in the middle of Fig.~\ref{figc2max}. It confirms the 
interpretation of the final state being a line solitary wave plus 
radiation. 

\begin{figure}[htb!]
 \includegraphics[width=0.49\textwidth]{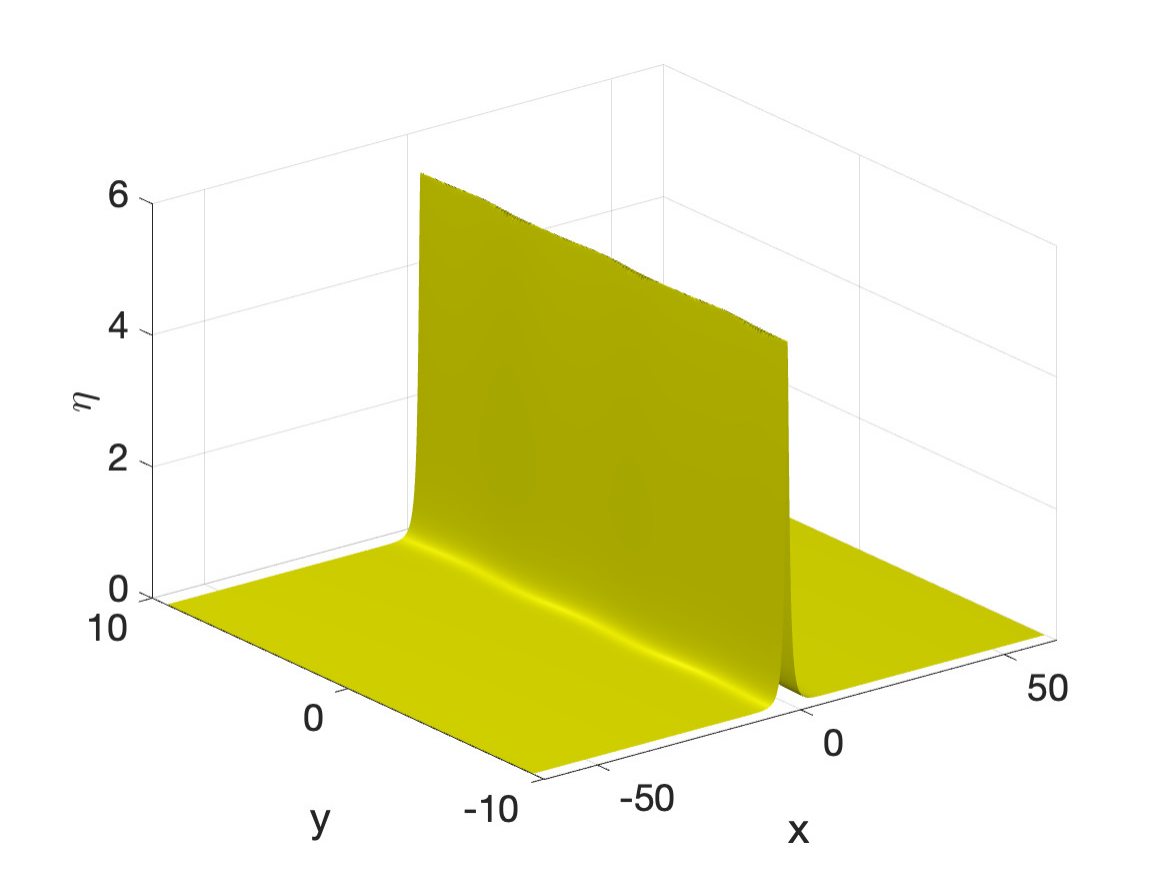}
 \includegraphics[width=0.49\textwidth]{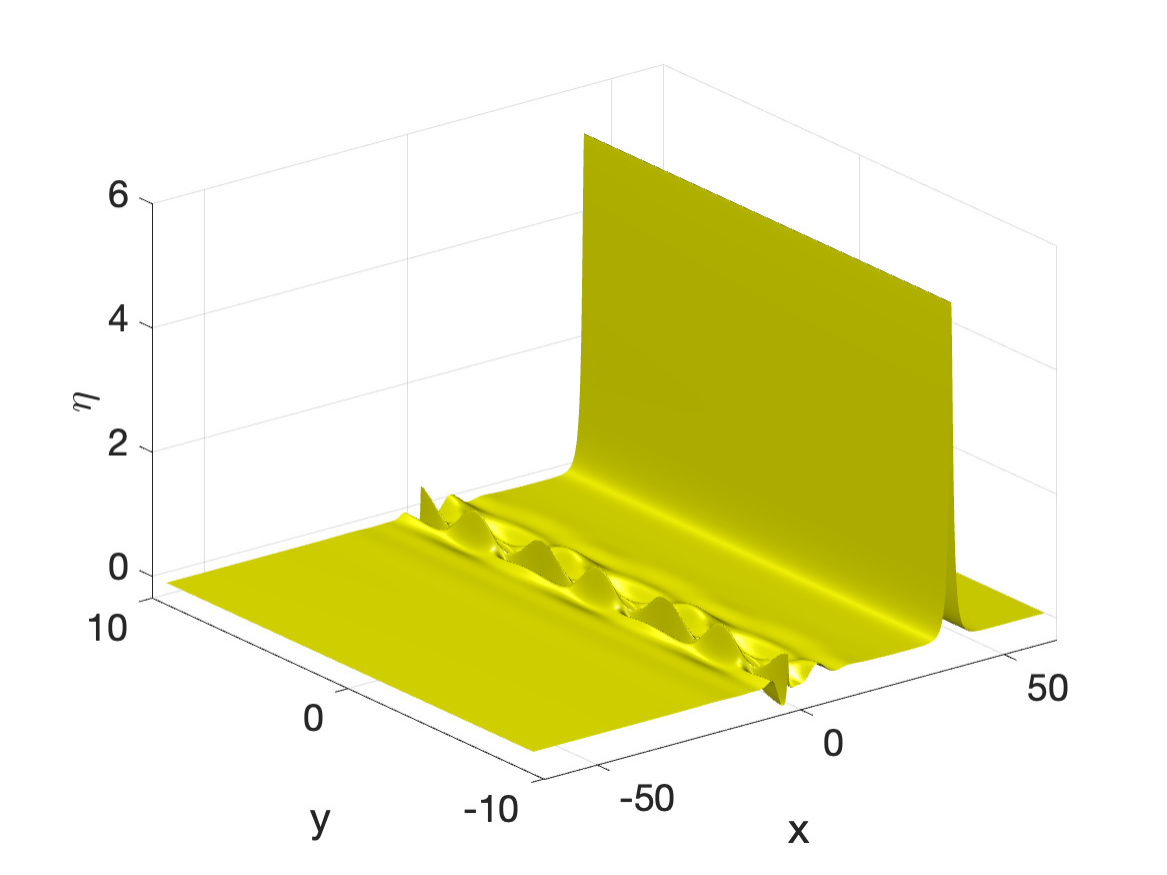}
 \caption{Solution to the AS system (\ref{2D}) for the initial data 
 (\ref{c2cos}), on the left $\eta$ for $t=0$, on the right for 
 $t=20$.}
 \label{figc2cos}
\end{figure}

A different kind of perturbation is to deform the line solitary wave 
periodically. Concretely we consider initial data of the form 
\begin{equation}
	\eta(x,y,0) = Q_{2}(x+0.4\cos(y)),\quad 
	v_{x}(x,y,0) = V_{2}(x),\quad v_{y}(x,y,0)=0.
	\label{c2cos}
\end{equation}
The initial condition for $\eta$ can be seen on the left of 
Fig.~\ref{figc2cos}. The deformed solitary wave emits some radiation 
and appears to reach another solitary wave as a final state as 
suggested by the figure on the right of Fig.~\ref{figc2cos} for 
$t=20$. 

\begin{figure}[htb!]
 \includegraphics[width=0.49\textwidth]{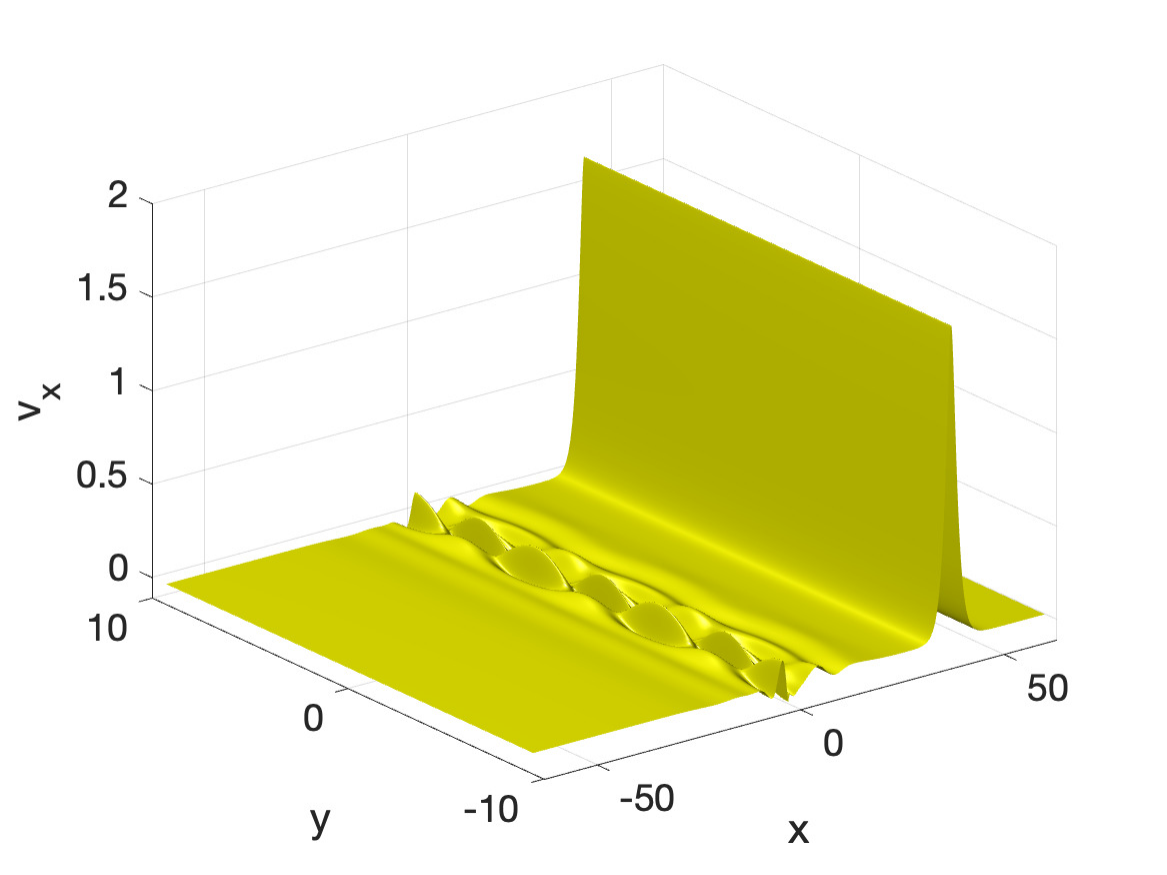}
 \includegraphics[width=0.49\textwidth]{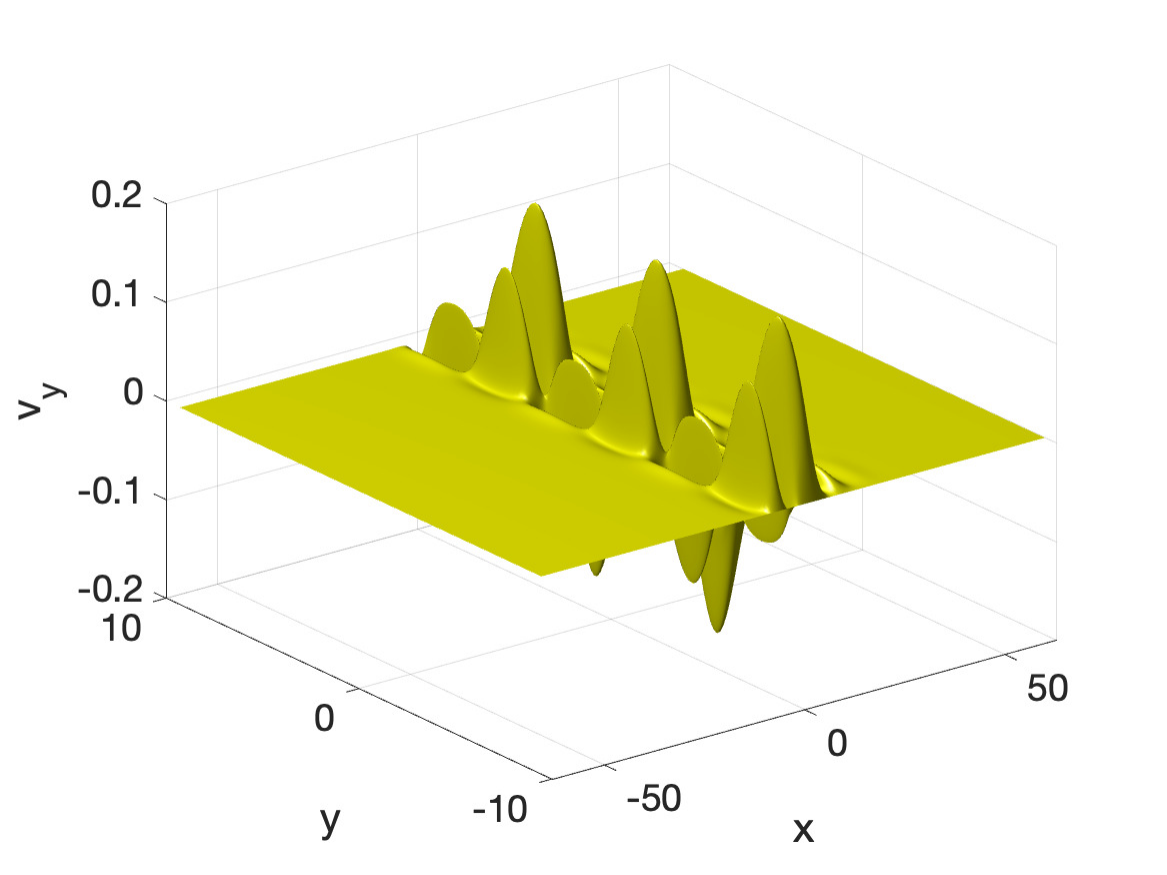}
 \caption{Solution to the AS system (\ref{2D}) for the initial data 
 (\ref{c2cos}) for $t=100$, on the left $v_{x}$, on the right 
 $v_{y}$. }
 \label{figc2cosv}
\end{figure}

This is confirmed by the velocities at the final time $t=20$ shown in 
Fig.~\ref{figc2cosv}. The $L^{\infty}$ norm of $\eta$ in the same 
figure appears to approach rapidly the final state as can be seen on 
the right of Fig.~\ref{figc2max}. 

\subsection{Perturbed line solitary wave with $c=1.1$}

\begin{figure}[htb!]
 \includegraphics[width=0.49\textwidth]{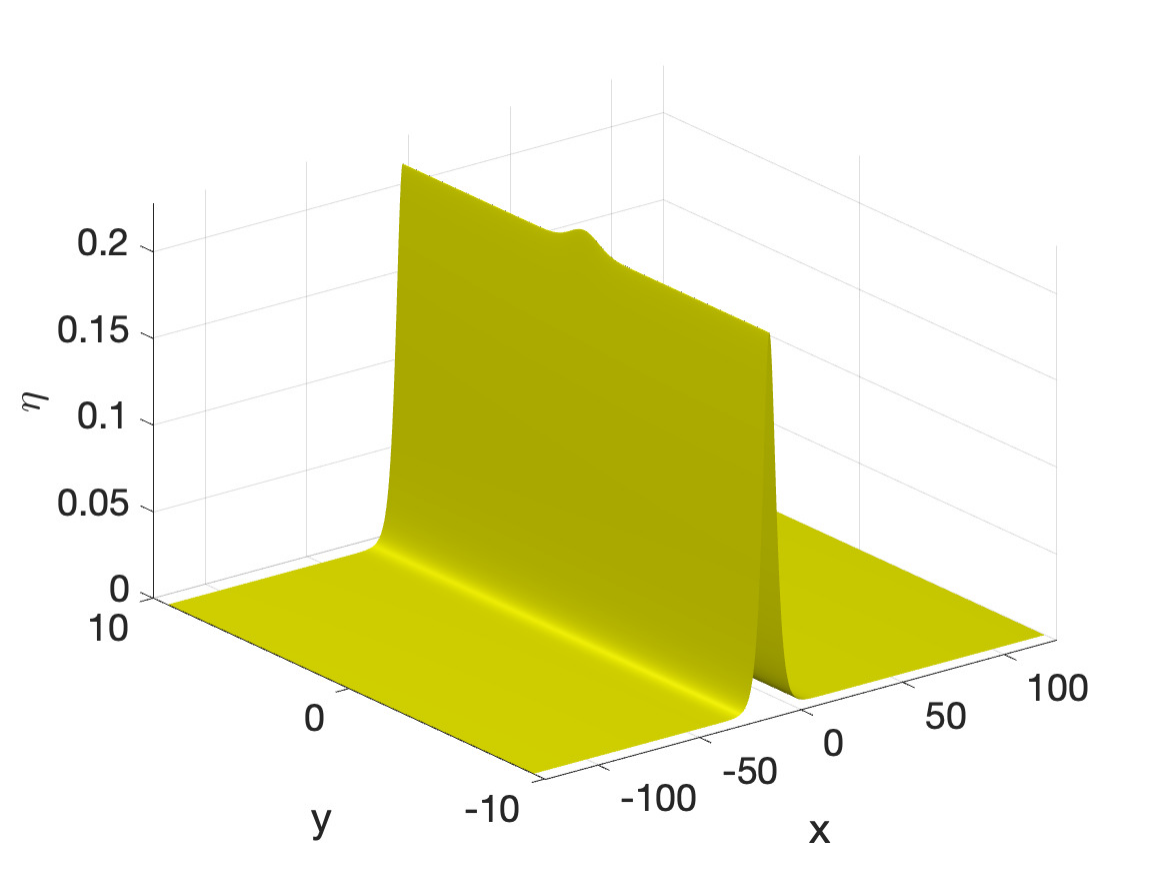}
 \includegraphics[width=0.49\textwidth]{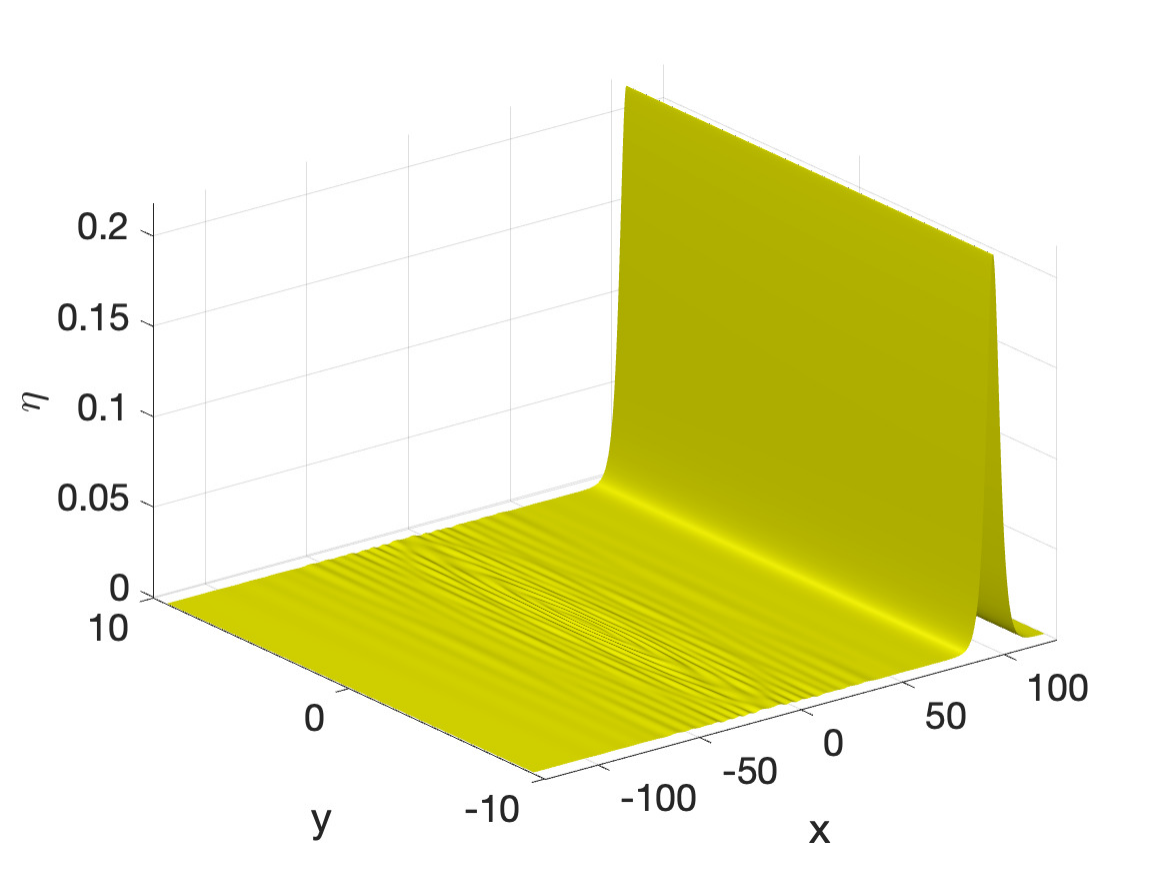}
 \caption{Solution to the AS system (\ref{2D}) for the initial data 
 (\ref{c11gauss}), on the left $\eta$ for $t=0$, on the right for 
 $t=100$.}
 \label{figc11gauss}
\end{figure}

As an example of solitary waves with velocities close to the 
threshold $c=1$, i.e., small slow solitary waves, we consider the 
case $c=1.1$. For the numerical experiments we use $N_{x}=2^{14}$ Fourier modes for $x\in 
40[-\pi,\pi]$ and $N_{y}=2^{7}$ Fourier modes for $y\in 3[-\pi,\pi]$ 
with $N_{t}=2*10^{4}$ time steps for $t<100$. We consider again 
localised perturbations with initial data of the form ($\kappa>0$, 
$\alpha>0$)
\begin{equation}
	\eta(x,y,0) = Q_{1.1}(x) + 0.01\exp(-x^{2}-\alpha y^{2}),\quad 
	v_{x}(x,y,0) = V_{1.1}(x),\quad v_{y}(x,y,0)=0.
	\label{c11gauss}
\end{equation}
This means the perturbation has a maximum of the order of $5\%$ of the 
wave crest. The initial data for $\eta$ and the solution for $t=20$ 
are shown in Fig.~\ref{figc11gauss}. It can be seen that the final 
state appears to be a line solitary wave of slightly different 
velocity as in the case for $c=2$. 

\begin{figure}[htb!]
 \includegraphics[width=0.49\textwidth]{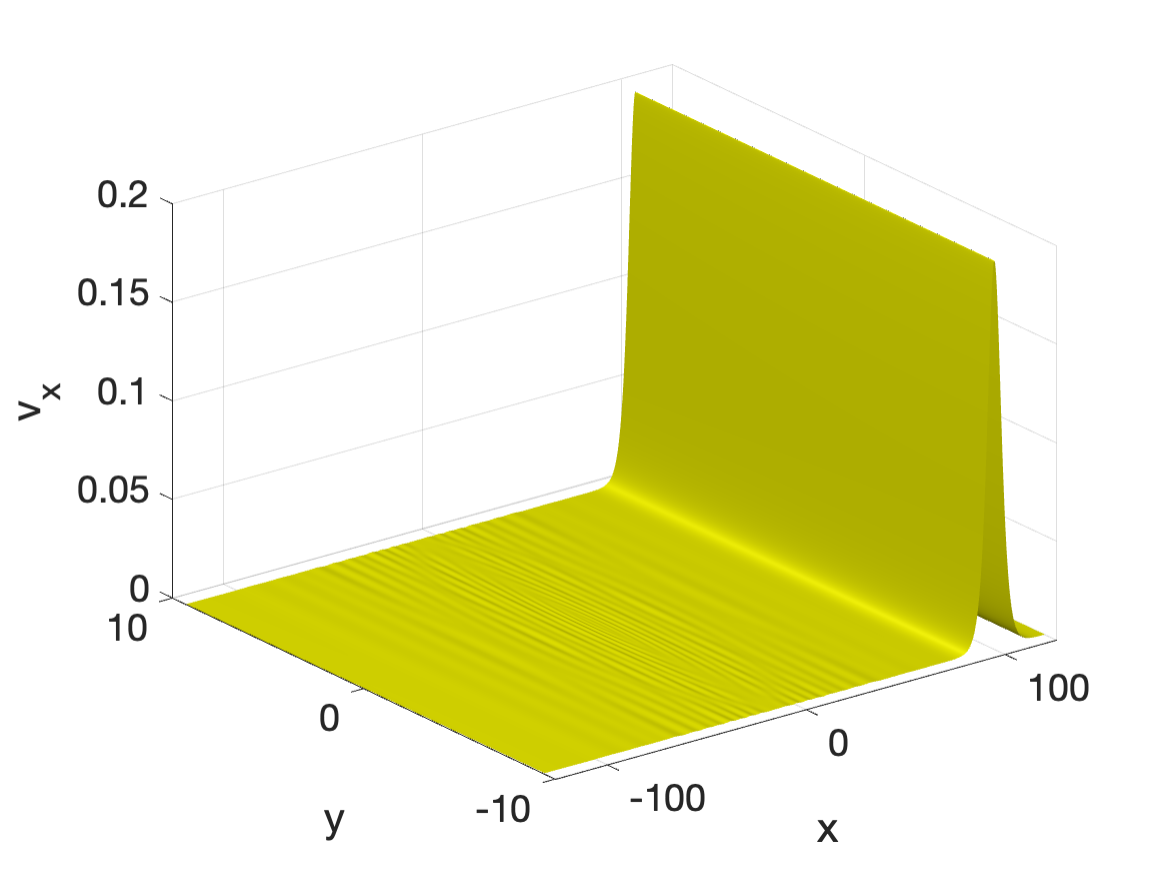}
 \includegraphics[width=0.49\textwidth]{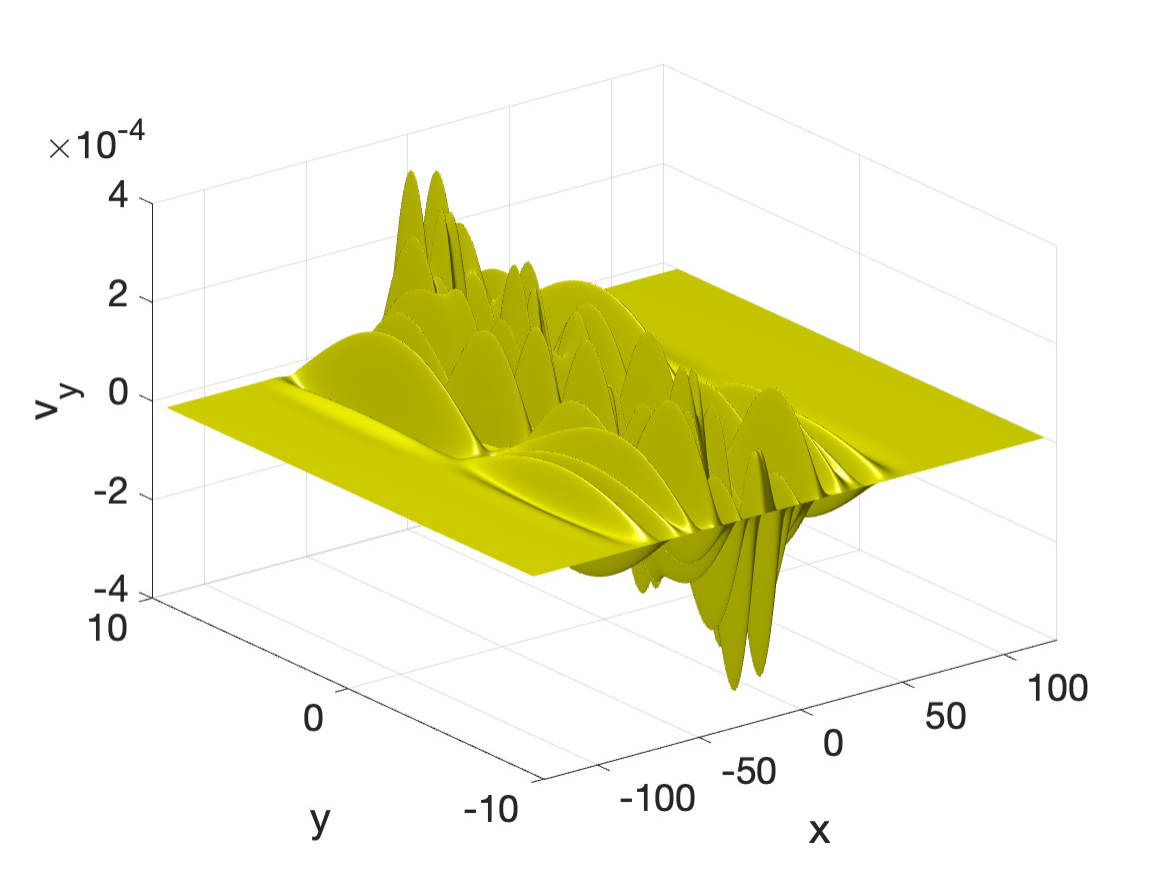}
 \caption{Solution to the AS system (\ref{2D}) for the initial data 
 (\ref{c11gauss}) for $t=100$, on the left $v_{x}$, on the right 
 $v_{y}$. }
 \label{figc11gaussv}
\end{figure}

The plots for the velocities  $v_{x}$ and $v_{y}$ are presented in 
Fig.~\ref{figc11gaussv}. Once more the final state appears to be 
a line solitary wave plus radiation. 

\begin{figure}[htb!]
 \includegraphics[width=0.49\textwidth]{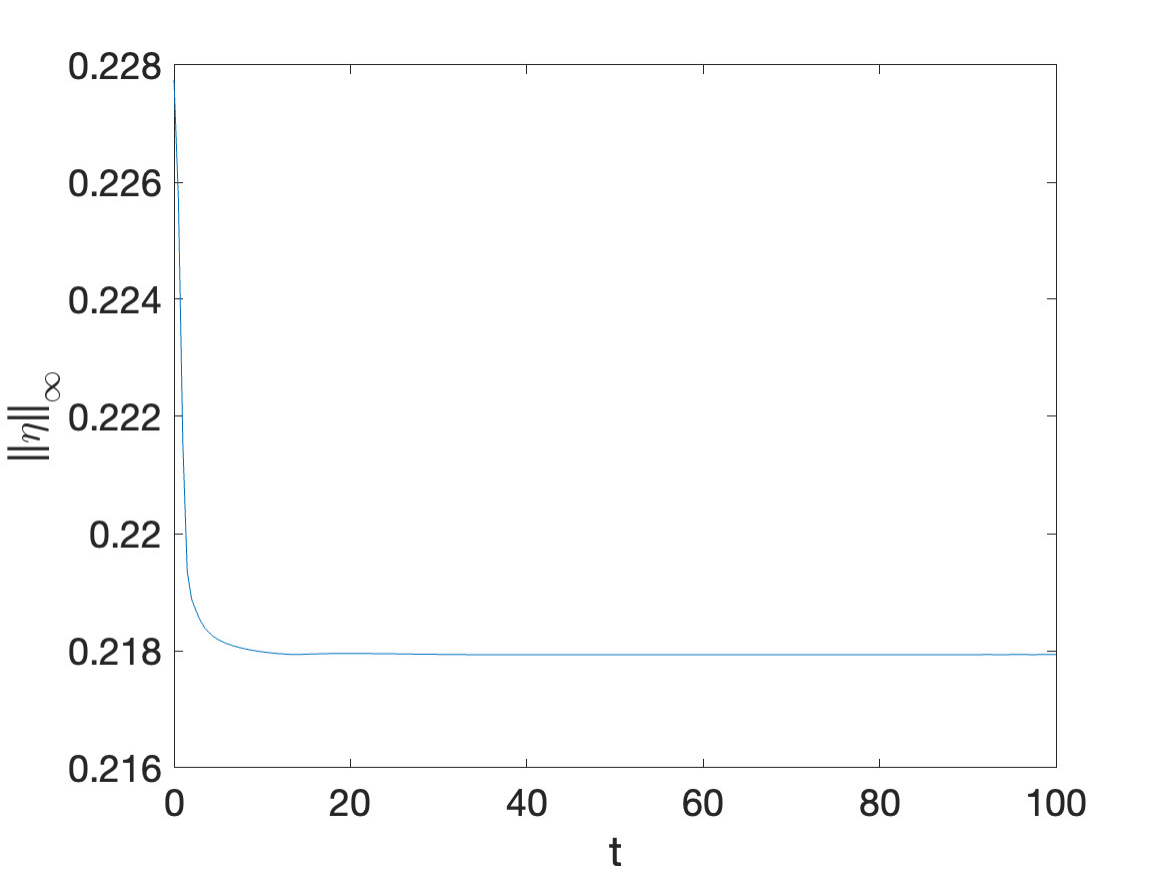}
 \includegraphics[width=0.49\textwidth]{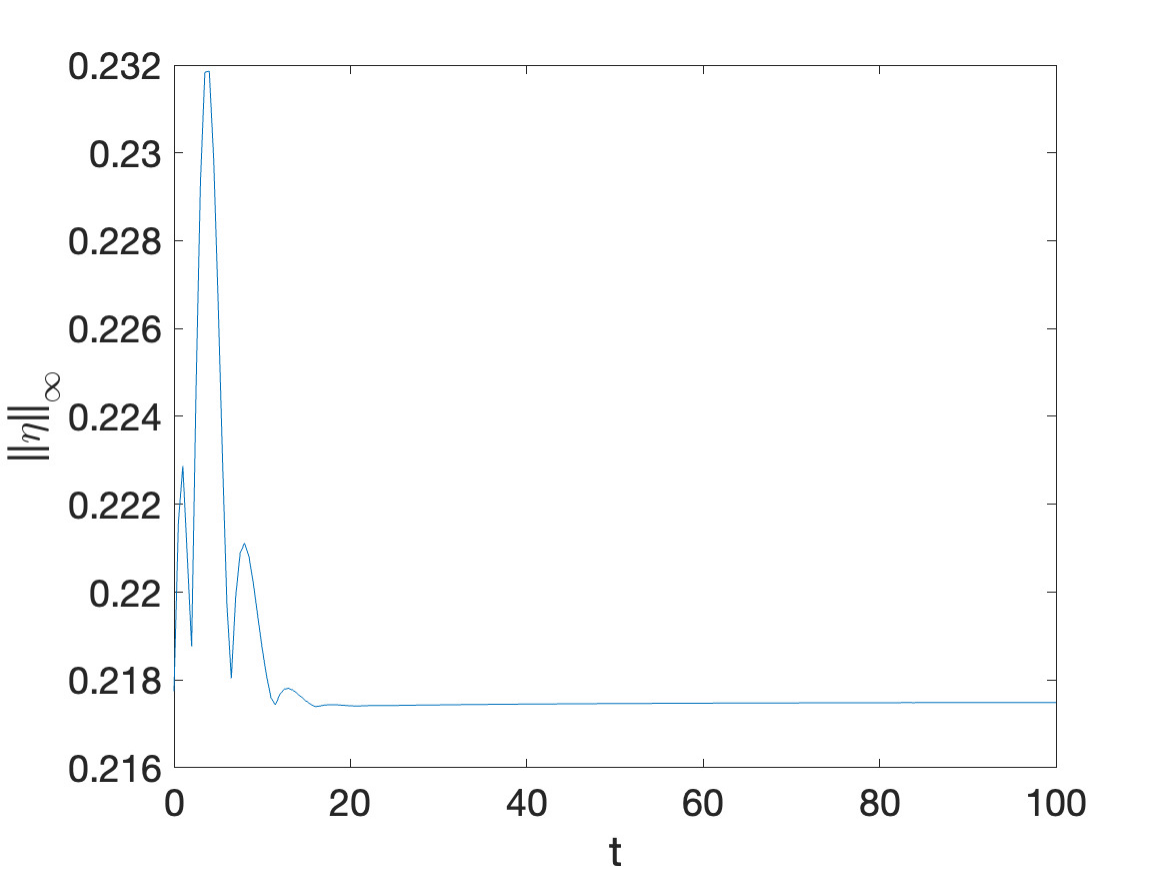}
 \caption{$L^{\infty}$ norm of the solution to the AS system (\ref{2D}) for the initial data 
 (\ref{c11gauss}) on the left, on the right for the initial data 
 (\ref{c11cos}).}
 \label{figc11gaussmax}
\end{figure}

This interpretation is confirmed by the $L^{\infty}$ norm of $\eta$ 
shown on the left of Fig.~\ref{figc11gaussmax}. It appears to settle 
on a constant level for larger $t$ indicating that the final state is 
a line solitary wave of slightly higher velocity than the unperturbed 
wave with $c=1.1$. 

\begin{figure}[htb!]
 \includegraphics[width=0.49\textwidth]{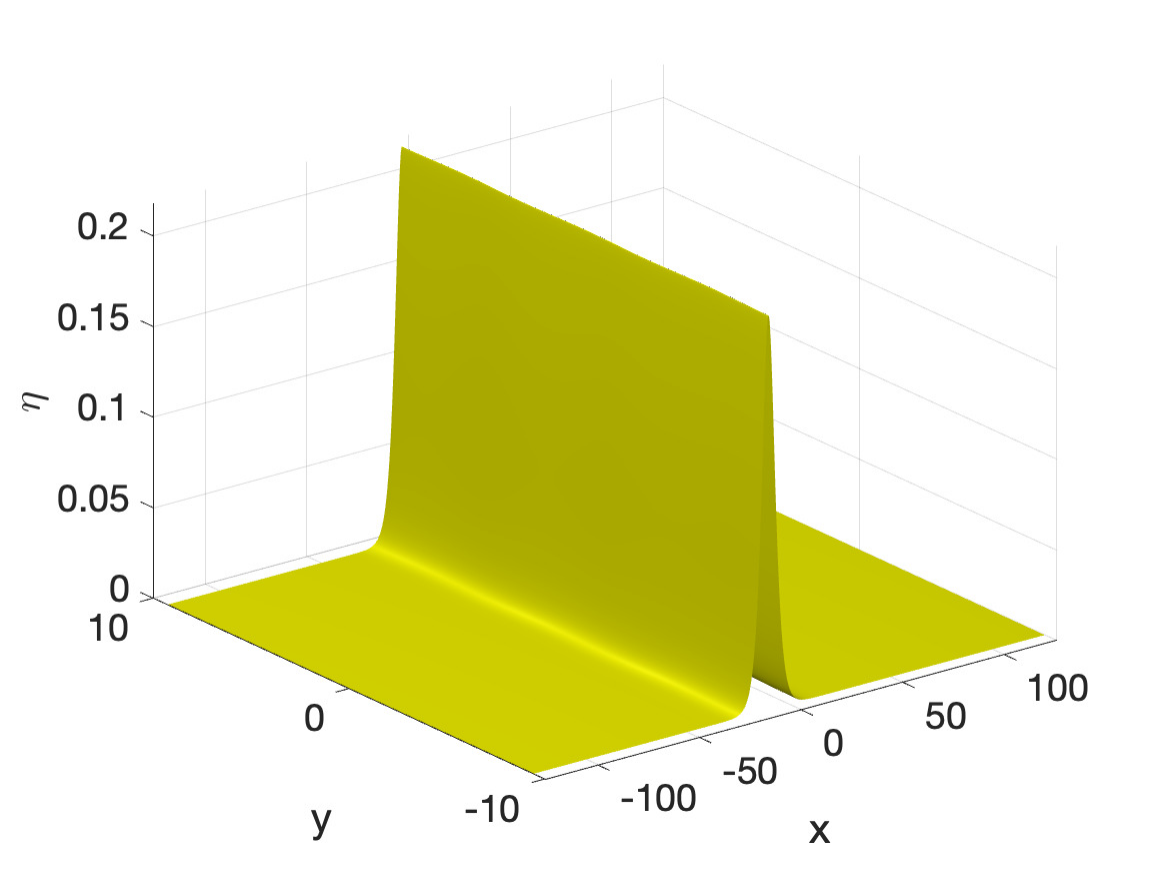}
 \includegraphics[width=0.49\textwidth]{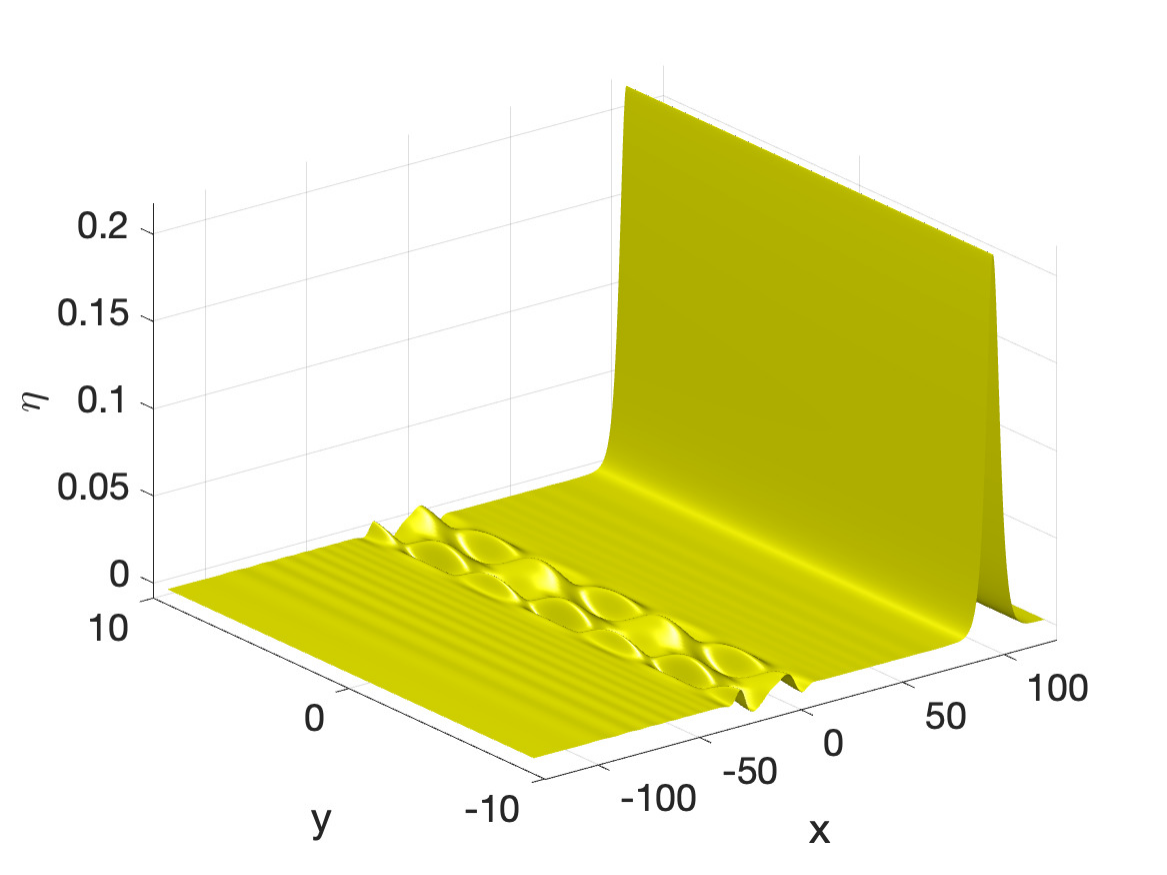}\\
 \includegraphics[width=0.49\textwidth]{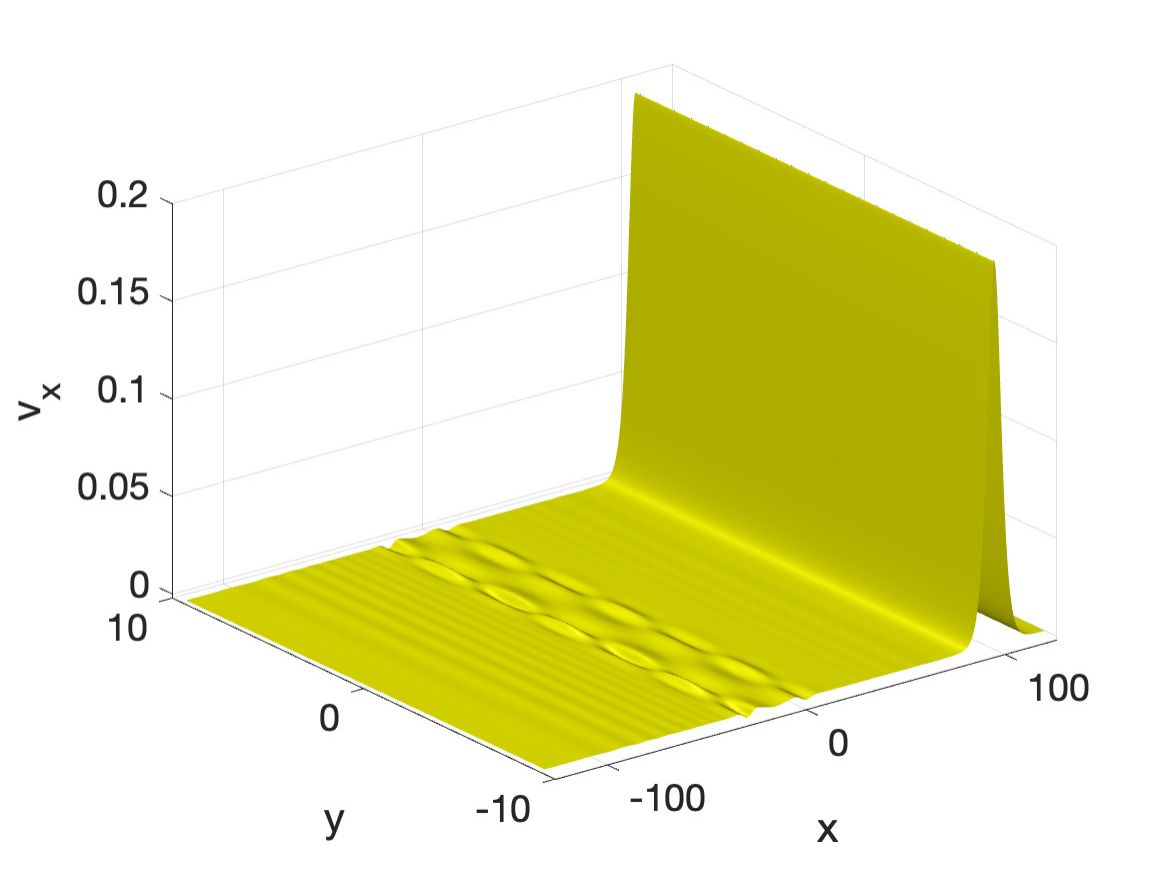}
 \includegraphics[width=0.49\textwidth]{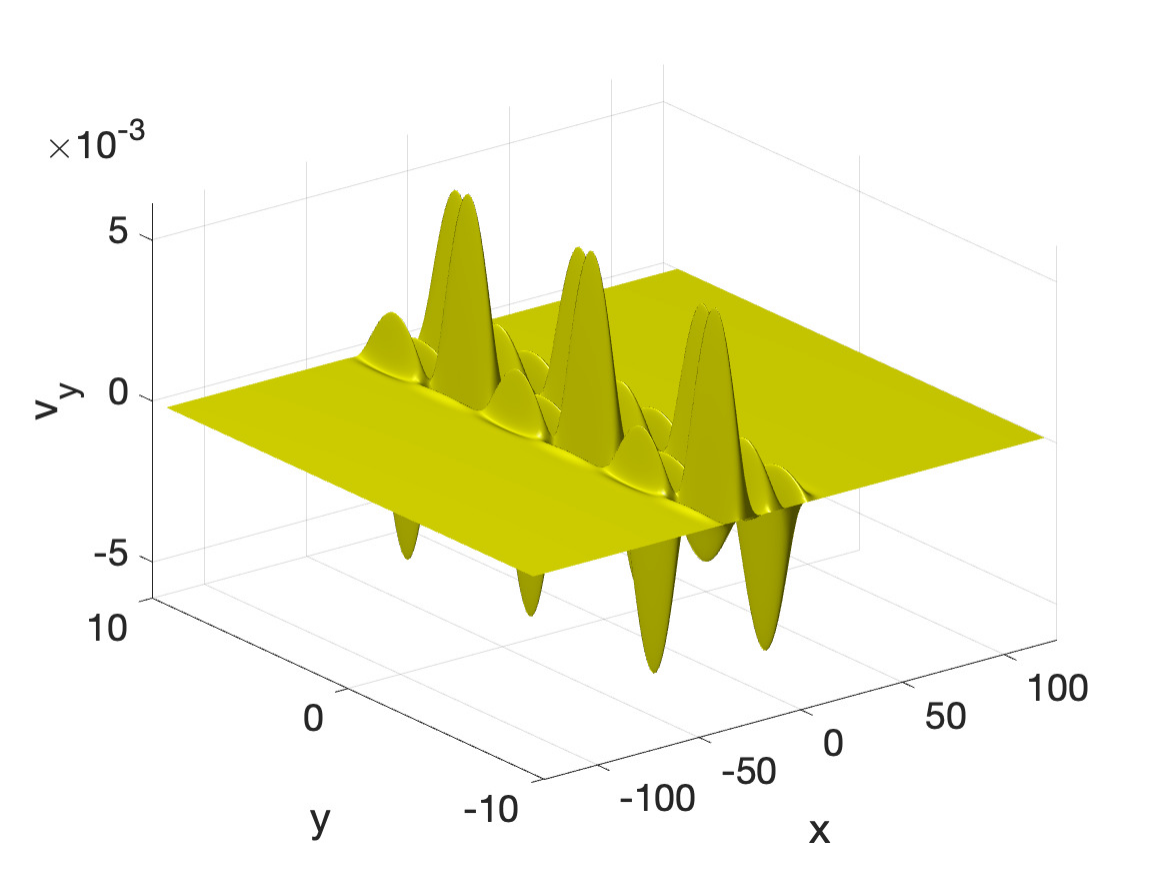}
 \caption{Solution to the AS system (\ref{2D}) for the initial data 
 (\ref{c11cos}), in the upper row on the left $\eta$ for $t=0$, on the right for 
 $t=100$, in the lower row on the left $v_{x}$ and on the right 
 $v_{x}$, both for $t=100$.}
 \label{figc11cos}
\end{figure}

A periodic deformation of the line solitary wave similar to 
(\ref{c2cos}), 
\begin{equation}
	\eta(x,y,0) = Q_{1.1}(x+0.4\cos(y)),\quad 
	v_{x}(x,y,0) = V_{1.1}(x),\quad v_{y}(x,y,0)=0,
	\label{c11cos}
\end{equation}
gives similar results as in the case $c=2$ as can be seen in 
Fig.~\ref{figc11cos}. The final state is once more a line solitary 
wave plus radiation. This is confirmed by the $L^{\infty}$ norm of 
$\eta$ in Fig.~\ref{figc11gaussmax} on the right. 

\section{Non-cavitation condition}
In this section we consider initial data which violate the 
non-cavitation condition or come close to doing so. We first recall 
and expand on the case in 1D discussed in \cite{KSAS}. Then we 
discuss similar examples in 2D. 

\subsection{Non-cavitation condition in 1D}
In \cite{KSAS}, we have considered initial data of the form 
$$\eta(x,0)=-\exp(-x^{2}),\quad v(x,0)=0.$$
We applied  the approach by Sulem, Sulem and Frisch in 
\cite{SSF} to numerically identify singularities.  It uses the fact
that an essential singularity in the complex plane of the form $u\sim 
(z-z_{j})^{\mu_{j}}$, $\mu_{j}\notin \mathbb{Z}$, 
with $z_{j}=\alpha_{j}-i\delta_{j}$ in the lower half 
plane ($\delta_{j}\geq 0$) implies    for $k\to\infty$  
the following asymptotic behavior of the Fourier 
transform (see e.g.~\cite{asymbook}, here denoted in the same way as the DFT),
\begin{equation}
    \hat{u}\sim 
    \sqrt{2\pi}\mu_{j}^{\mu_{j}+\frac{1}{2}}e^{-\mu_{j}}\frac{(-i)^{\mu_{j}+1}}{k^{\mu_{j}+1}} e^{-ik\alpha_{j}-k\delta_{j}}
    \label{fourierasym}.
\end{equation}
For a single such singularity with positive $\delta_{j}$, the modulus of the Fourier 
transform decreases exponentially for large $|k|$ until
$\delta_{j}=0$, when 
this modulus  has an algebraic dependence 
on $|k|$ for large $|k|$. The same behavior is expected in the discrete version of the 
Fourier transform, the DFT. In a dynamic situation, both $\delta_{j}$ 
and $\mu_{j}$ are expected to depend on $t$, and the time that 
$\delta_{j}$ vanished corresponds to the singularity hitting the real
axis, i.e., to the solution forming a singularity. This will be 
identified by fitting the DFT coefficients to the (\ref{fourierasym}) 
via some linear regression.  
Note that this fitting  is much less reliable for the factor $\mu_{j}$ 
than for the factor $\delta_{j}$ since this is an algebraic correction to 
an exponential term. In the 1D case \cite{KSAS}, we used $2^{18}$ 
DFT modes, but in 2D below we will be limited on the used computers to values 
of $2^{12}$ to $2^{13}$. Note that if there are several singularities 
of the form (\ref{fourierasym}), $j=1,\ldots,N_{s}$, the asymptotic 
expression (\ref{fourierasym}) will involve a summation over $j$, but 
this will not be considered further (therefore we will drop the index 
$j$ in the following). 

\begin{figure}[htb!]
 \includegraphics[width=0.49\textwidth]{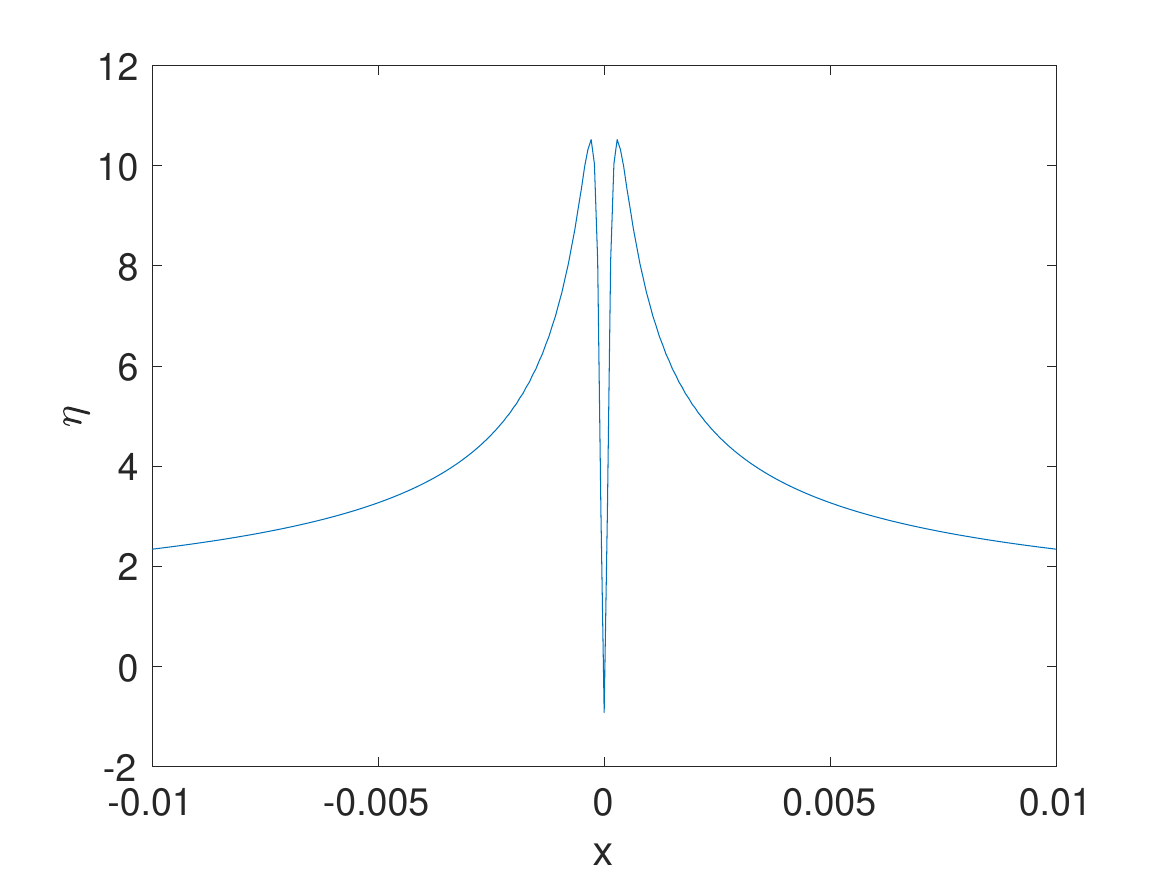}
 \includegraphics[width=0.49\textwidth]{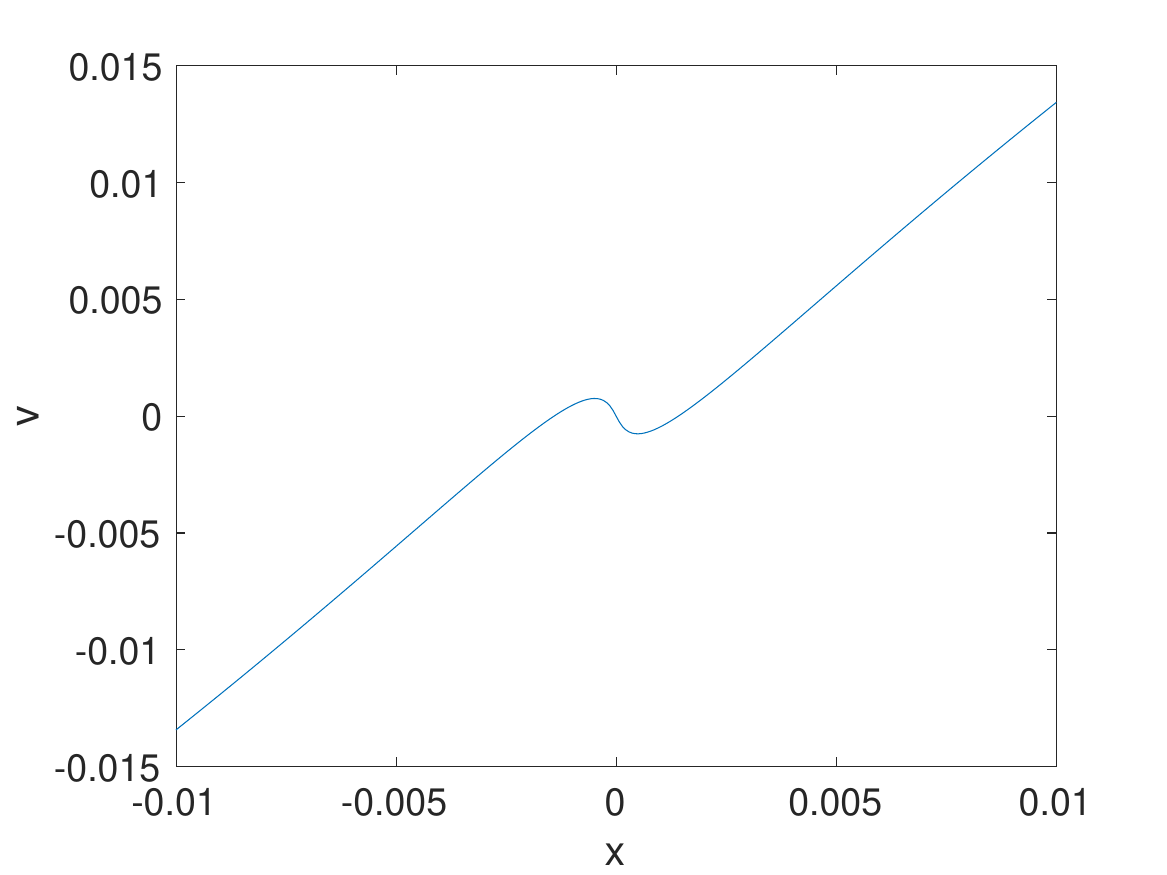}
  \caption{Solution to the AS system in 1D for the initial data 
 $\eta(x,0)=-\exp(-x^{2}),\quad v(x,0)=0$ for $t=4.681$ in a 
 close-up, on the left $\eta$, on the right $v$.}
 \label{AS1dcloseup}
\end{figure}

For the Gaussian example, the code in \cite{KSAS} broke for $t\sim 
4.681$ since the fitted value of $\delta$ vanished. We show the 
solution at the final recorded time in a close-up in 
Fig.~\ref{AS1dcloseup}. It can be seen that though the $L^{\infty}$ 
norm of $\eta$ gets large, a possible cusp near the origin appears to 
be the singular feature. The function $v$ on the right of the same 
figure appears to stay regular though the fitted value of $\delta$ 
for $v$ vanishes at the same time. 

\begin{figure}[htb!]
 \includegraphics[width=0.49\textwidth]{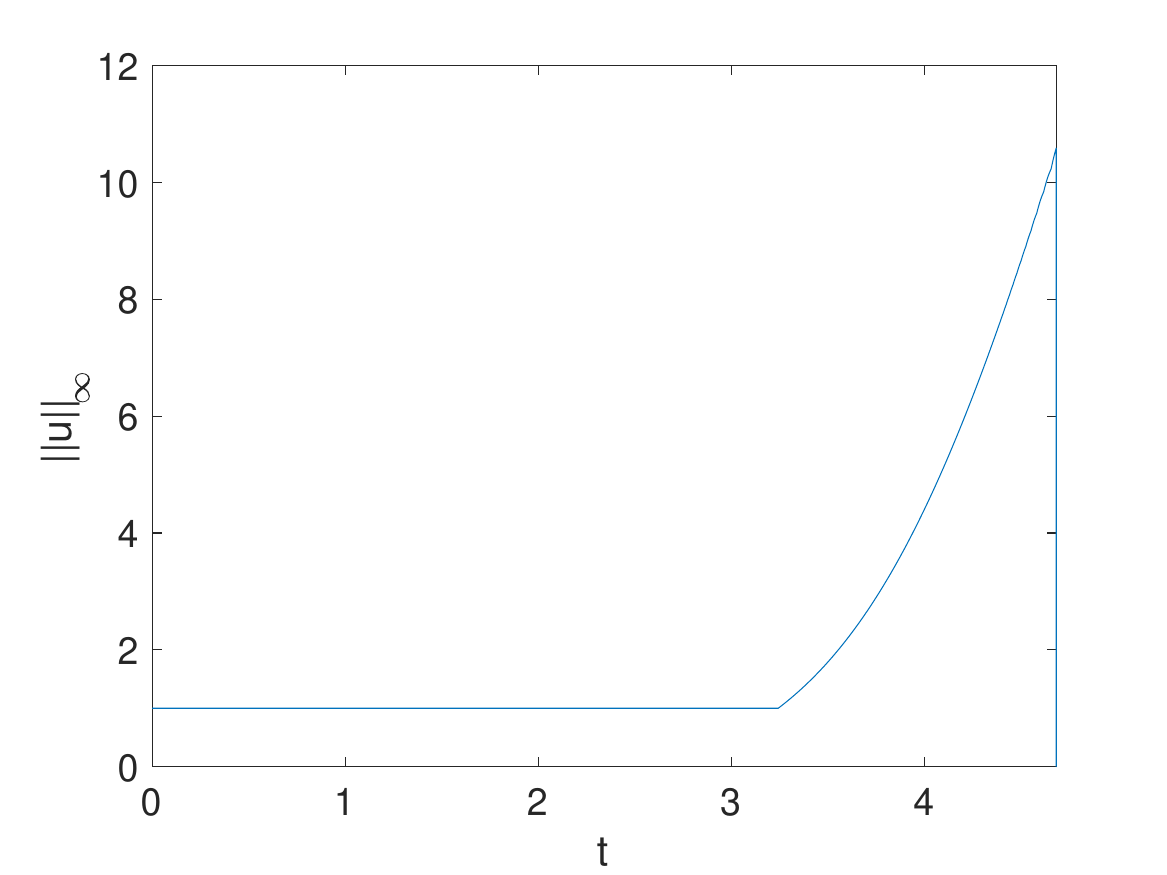}
 \includegraphics[width=0.49\textwidth]{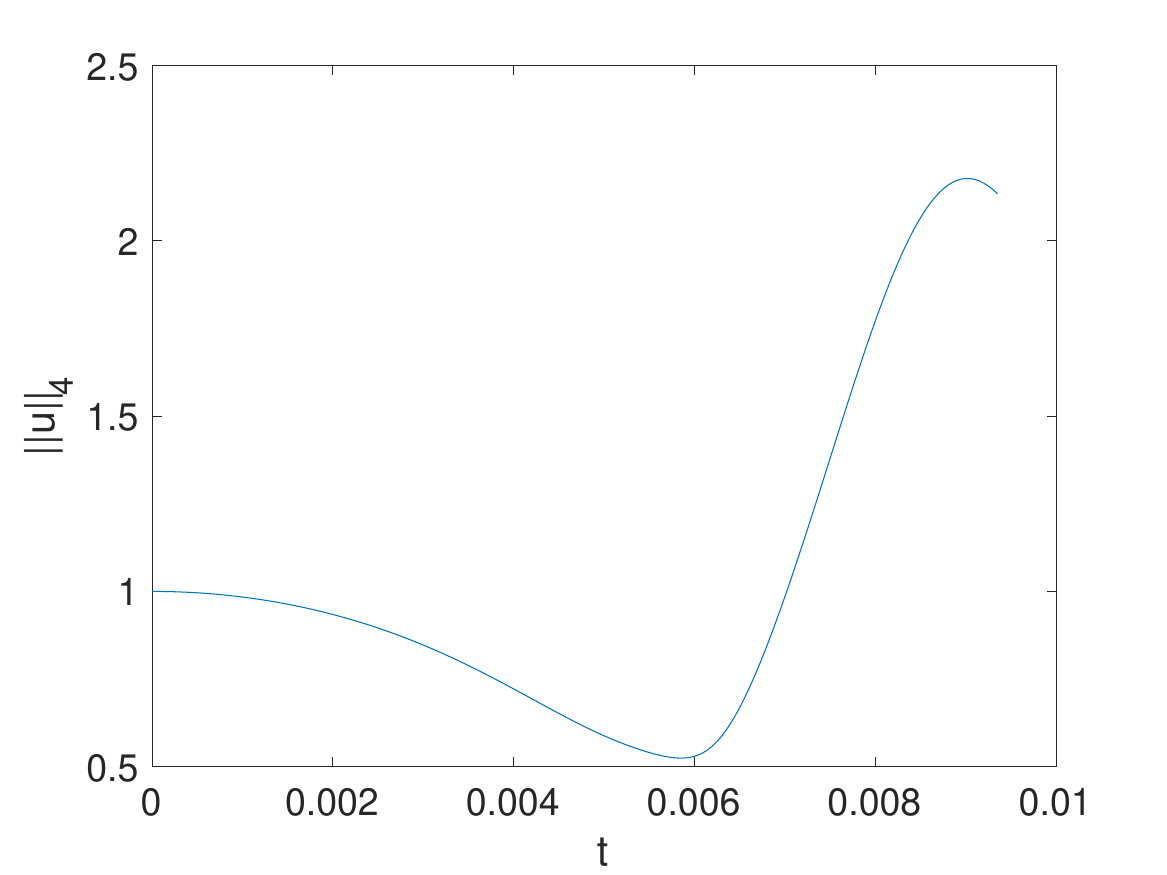}\\
 \includegraphics[width=0.49\textwidth]{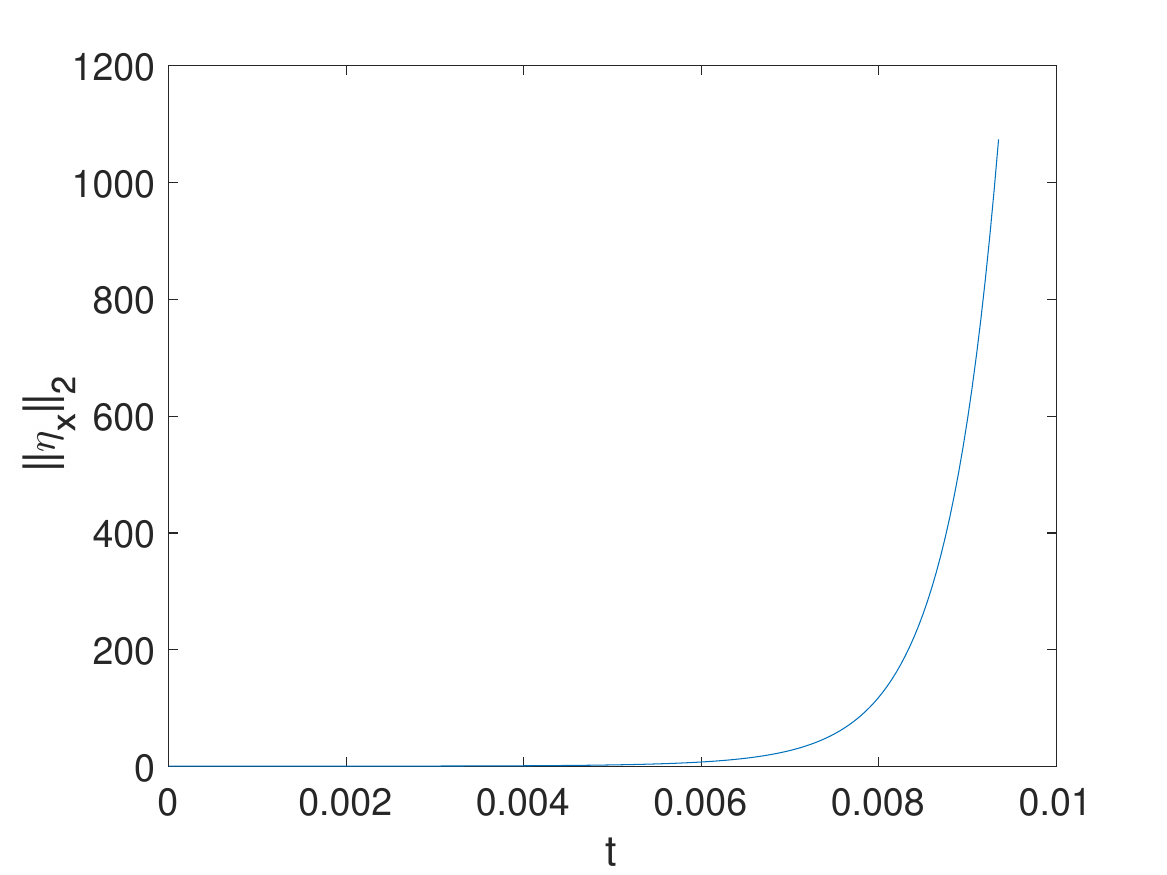}
 \includegraphics[width=0.49\textwidth]{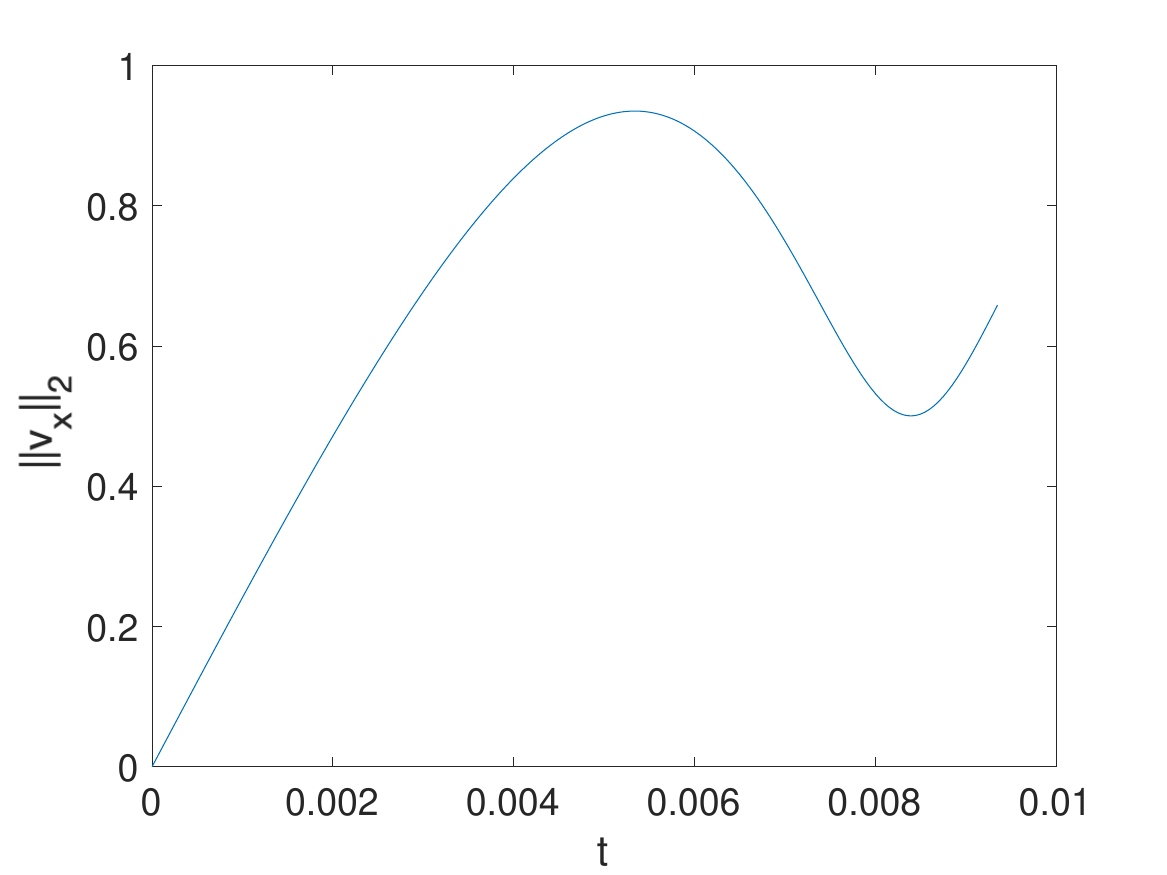}
  \caption{Solution to the AS system in 1D for the initial data 
 $\eta(x,0)=-\exp(-x^{2}),\quad v(x,0)=0$, in the upper row on the 
 left the $L^{\infty}$ norm, on the right the $L^{4}$ norm, in the 
 lower row the $L^{2}$ norm of the gradient of $\eta$ on the left 
 (normalized to 1 for $t=0$ as the $L^{4}$ norm) and 
 of $v$ on the right.}
 \label{AS1dnorms}
\end{figure}

The found value for $\mu$ for the function $\eta$ is negative ($\sim 
-0.85$) which would indicate an $L^{\infty}$ blow-up, whereas the 
corresponding value for $v$ is positive ($\sim1/3$) which would 
indicate the formation of a cusp. However as mentioned, these 
values have to be taken with a grain of salt since they are 
numerically problematic as an algebraic correction to an exponential 
decrease. To get a better insight into the behavior of the solution, 
we therefore consider various norms in Fig.~\ref{AS1dnorms}. Whereas 
the $L^{\infty}$ norm of $\eta$ appears to grow, this is not the case 
for various $L^{p}$ norms, for instance the $L^{4}$ norm which 
appears to be bounded. Though numerically one cannot get close enough 
to a potential singularity, this is clearly an indication that there 
is no $L^{\infty}$ blow-up to be observed in this example. The 
$L^{2}$ norms of the gradients in the lower row of the figure 
indicate, however, the blow-up of the gradient of $\eta$, but not of 
$v$. Thus it appears that there is a cusp formation in $\eta$ in 
finite time whereas $v$ stays regular at this time (of course we 
cannot rule out that higher Sobelev norms explode, it is just that 
$v$ is more regular at this time than $\eta$). 

\subsection{Non-cavitation condition in 2D}

Concretely we 
look at initial data of the form (here $\kappa<0$, $\alpha>0$)
\begin{equation}
	\eta(x,y,0) =  \kappa \exp(-(x^{2}+\alpha y^{2}),\quad v_{x}(x,y,0) = 
	v_{y}(x,y,0) = 0.
	\label{NCini}
\end{equation}
Note that we cannot characterize potential singularities in the 
solutions due to a lack of resolution on the available computers.

\begin{figure}[htb!]
 \includegraphics[width=0.49\textwidth]{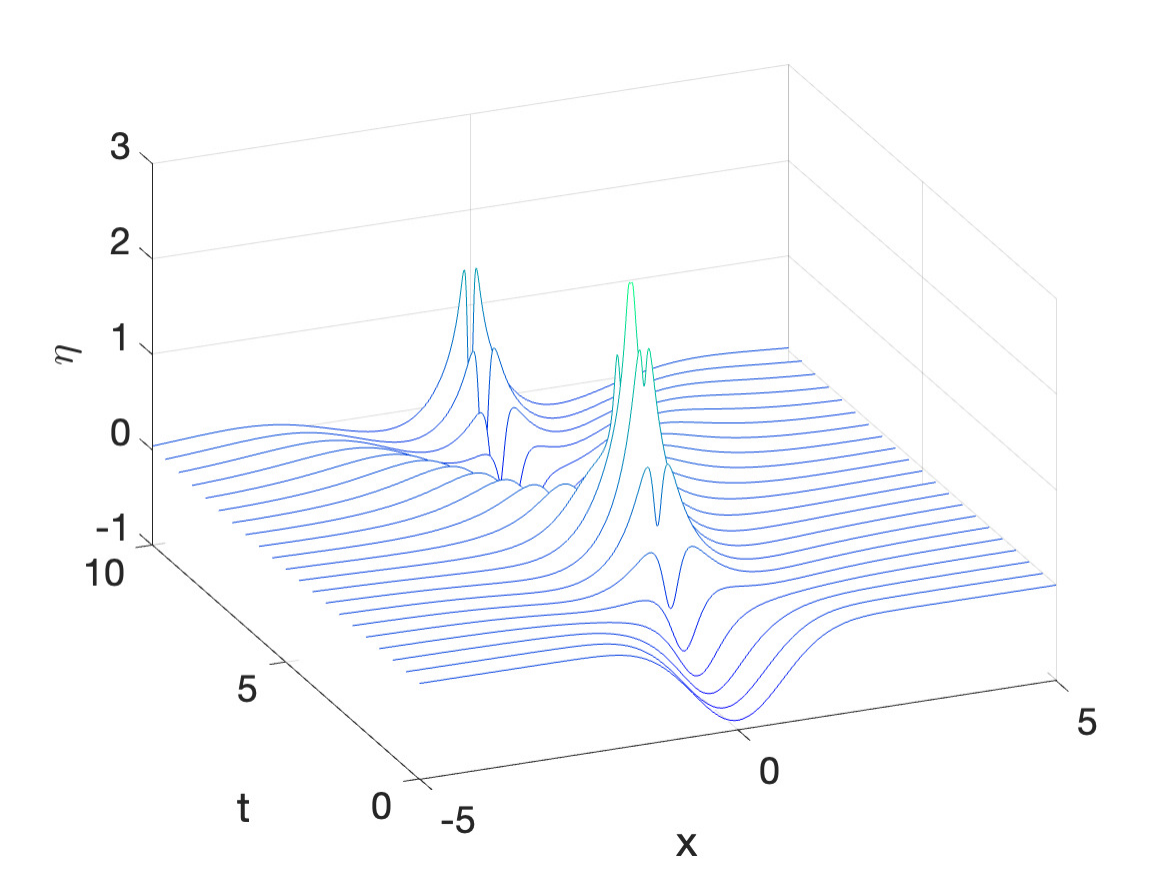}
 \includegraphics[width=0.49\textwidth]{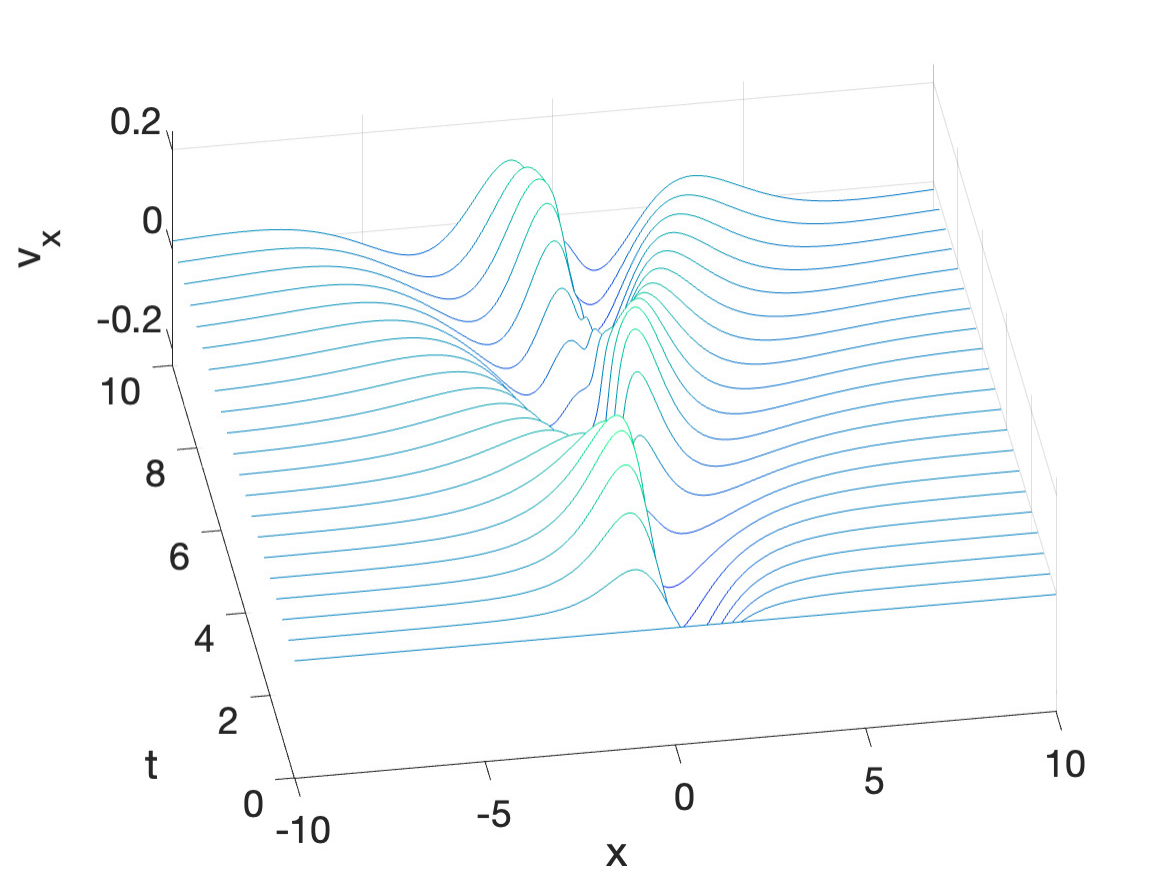}
 \caption{Solution to the AS system (\ref{2D}) for the initial data 
 (\ref{NCini}) with $\kappa=-0.9$ and $\alpha=1$ on the $x$-axis, on 
 the left $\eta$, on the right $v_{x}$.}
 \label{figm09gaussaxis}
\end{figure}

First we study the case $\kappa=-0.9$, i.e., initial data satisfying 
the non-cavitation condition, but close to violating it. 
We use $N_{x}=N_{y}=2^{12}$ DFT modes for $x\in 5[-\pi,\pi]$, $y\in 
5[-\pi,\pi]$ and $N_{t}=10^{4}$ time steps for the considered time 
interval. The solution 
stays smooth for all $t\leq10$, there is as expected no indication of 
the appearance of a singularity. To illustrate the interesting 
dynamics of the solution, we first show it on the $x$-axis in 
Fig.~\ref{figm09gaussaxis} (recall the radial symmetry of the 
solution). For the function $\eta$ shown on the left of the figure, 
there is a closing of the initial gap in the form of two peaks coming 
from the edges and eventually colliding to a central peak of high 
elevation, that then leads to another valley before recombining. The 
velocity $v_{x}$ on the right of the same figure shows a similar 
agitated behavior. 

\begin{figure}[htb!]
 \includegraphics[width=0.49\textwidth]{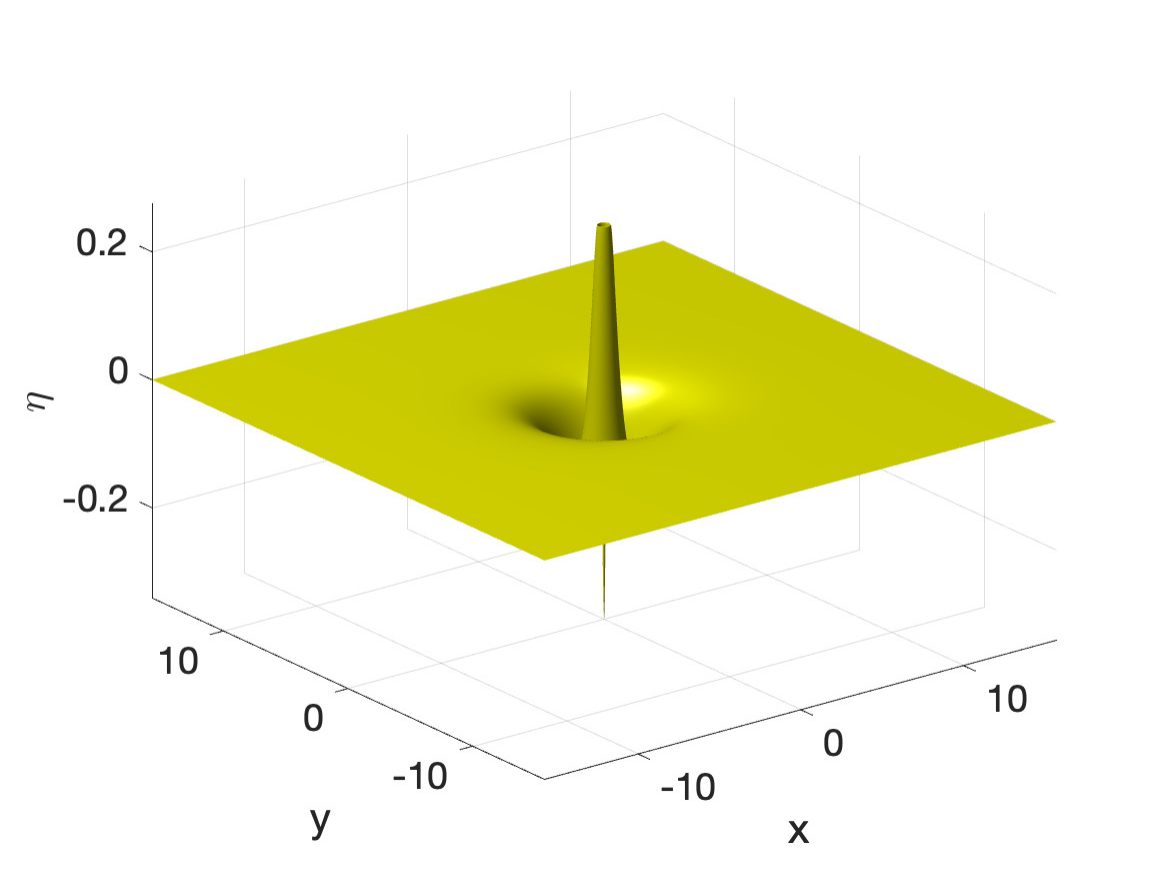}
 \includegraphics[width=0.49\textwidth]{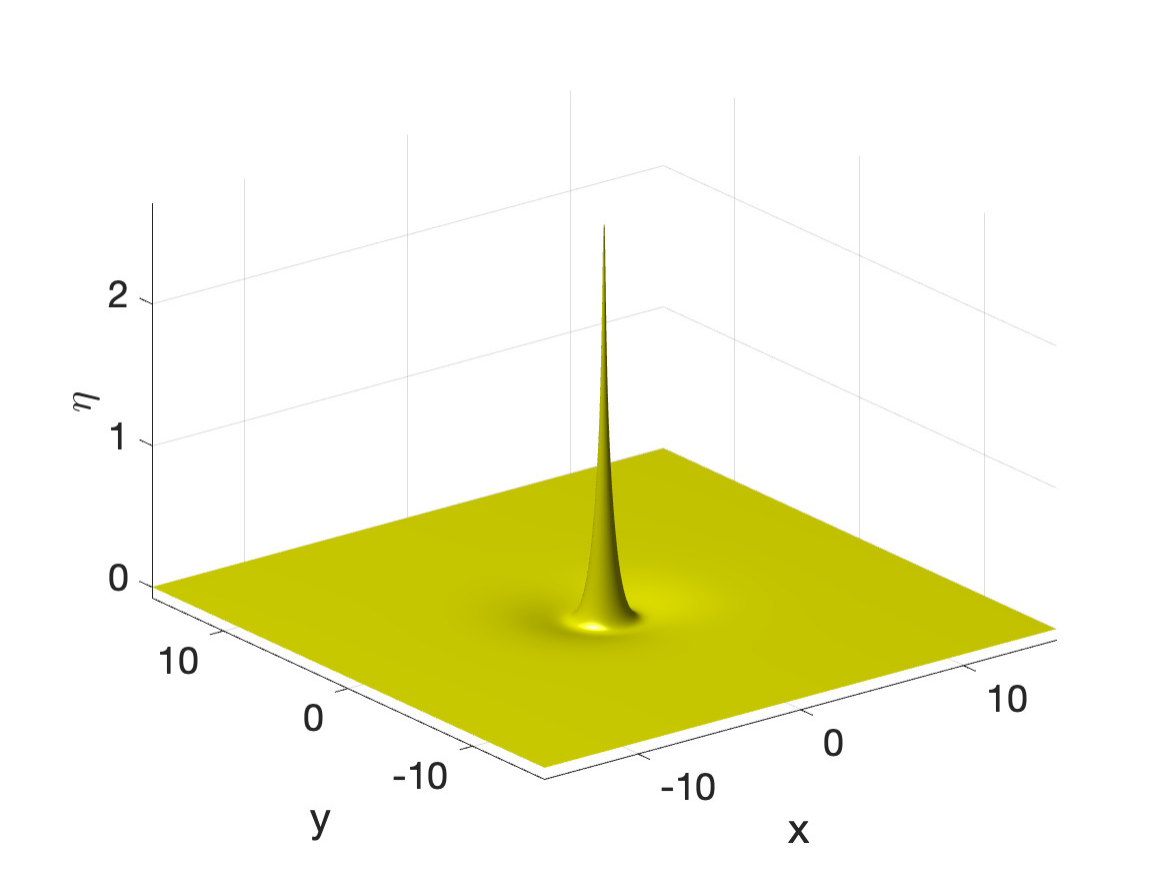}\\
 \includegraphics[width=0.49\textwidth]{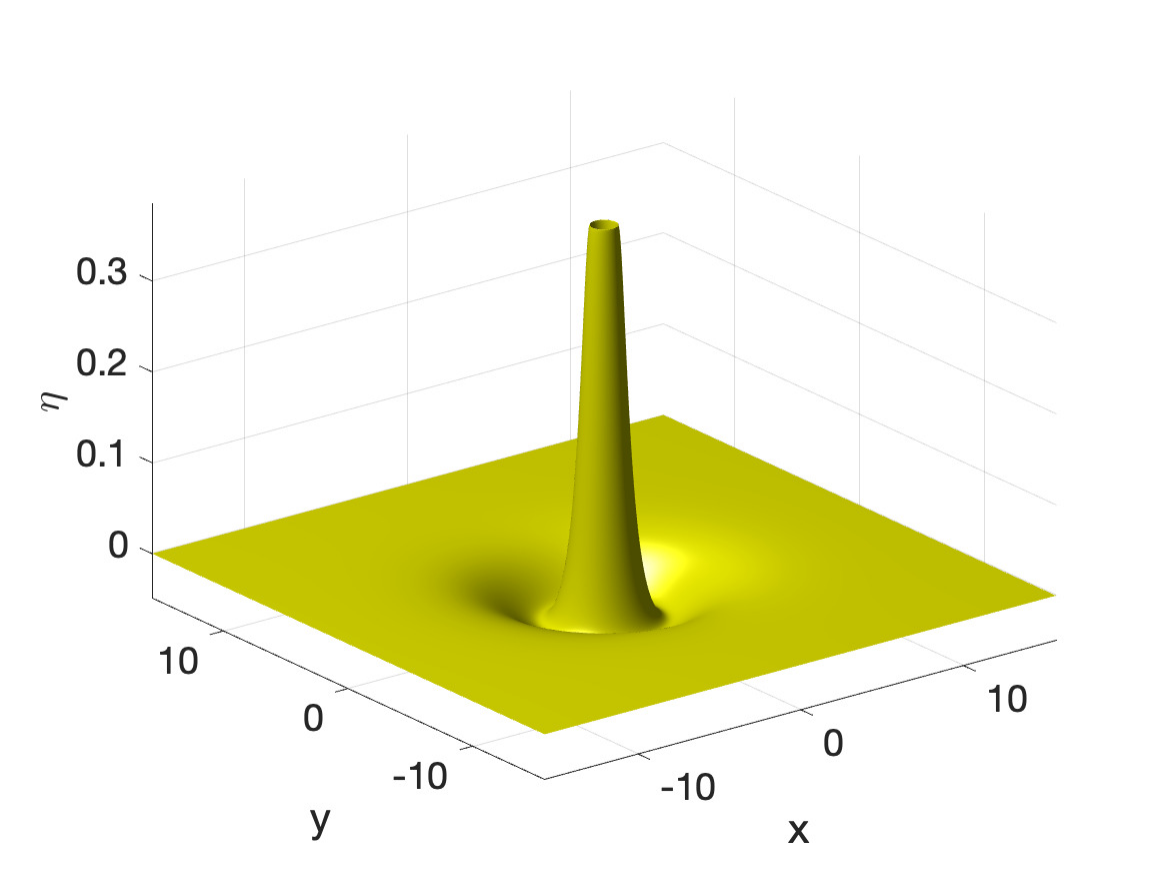}
 \includegraphics[width=0.49\textwidth]{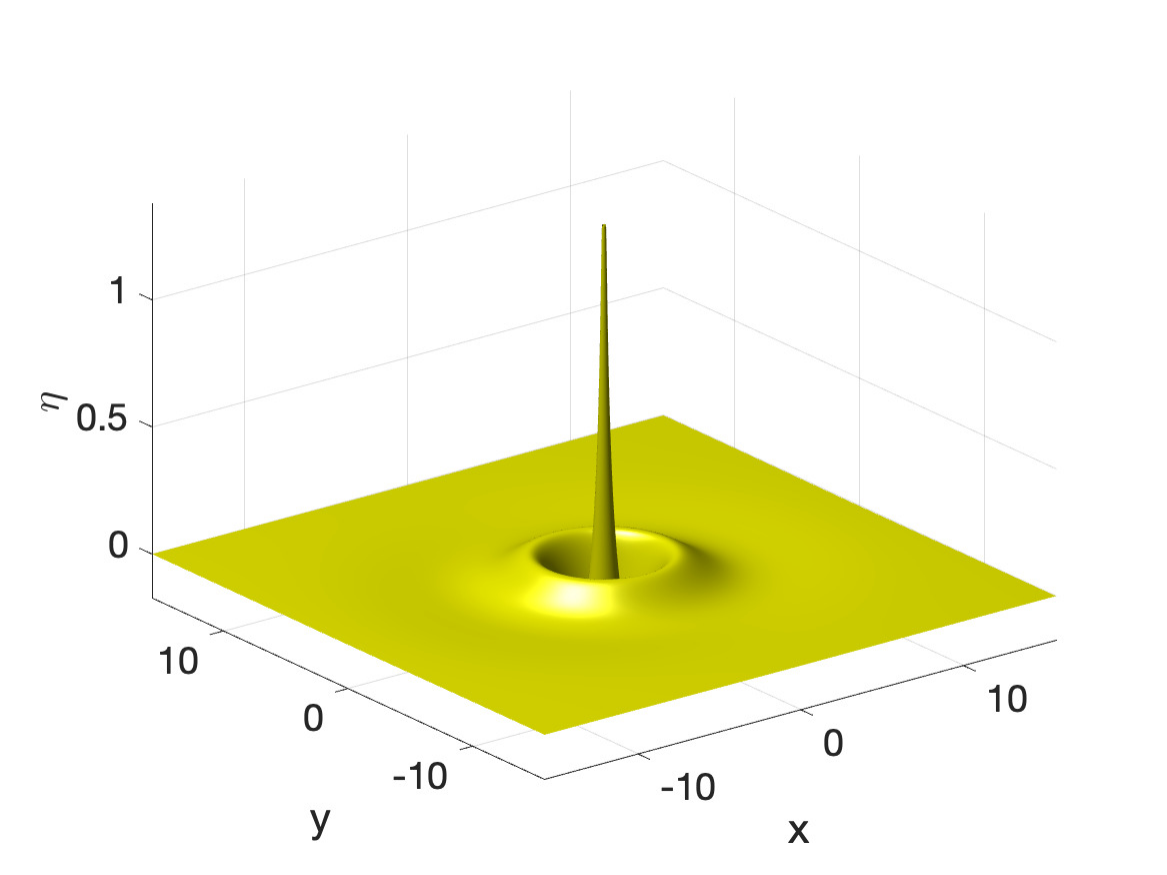}
 \caption{Solution $\eta$ to the AS system (\ref{2D}) for the initial data 
 (\ref{NCini}) with $\kappa=-0.9$ and $\alpha=1$, in the upper row on the left for 
 $t=2.5$,  on the right for 
 $t=4$, in the lower row on the left for $t=5.1$ and on the right 
for $t=10$.}
 \label{figm09gausseta}
\end{figure}

To illustrate the 2D situation, we show the solution $\eta$ for 
several values of time in Fig.~\ref{figm09gausseta}. It is clearly 
visible that the peak for $t=4$ is almost cusp-like, but then the 
solution forms an annular structure and a pronounced minimum appears 
inside the former peak. 

\begin{figure}[htb!]
 \includegraphics[width=0.49\textwidth]{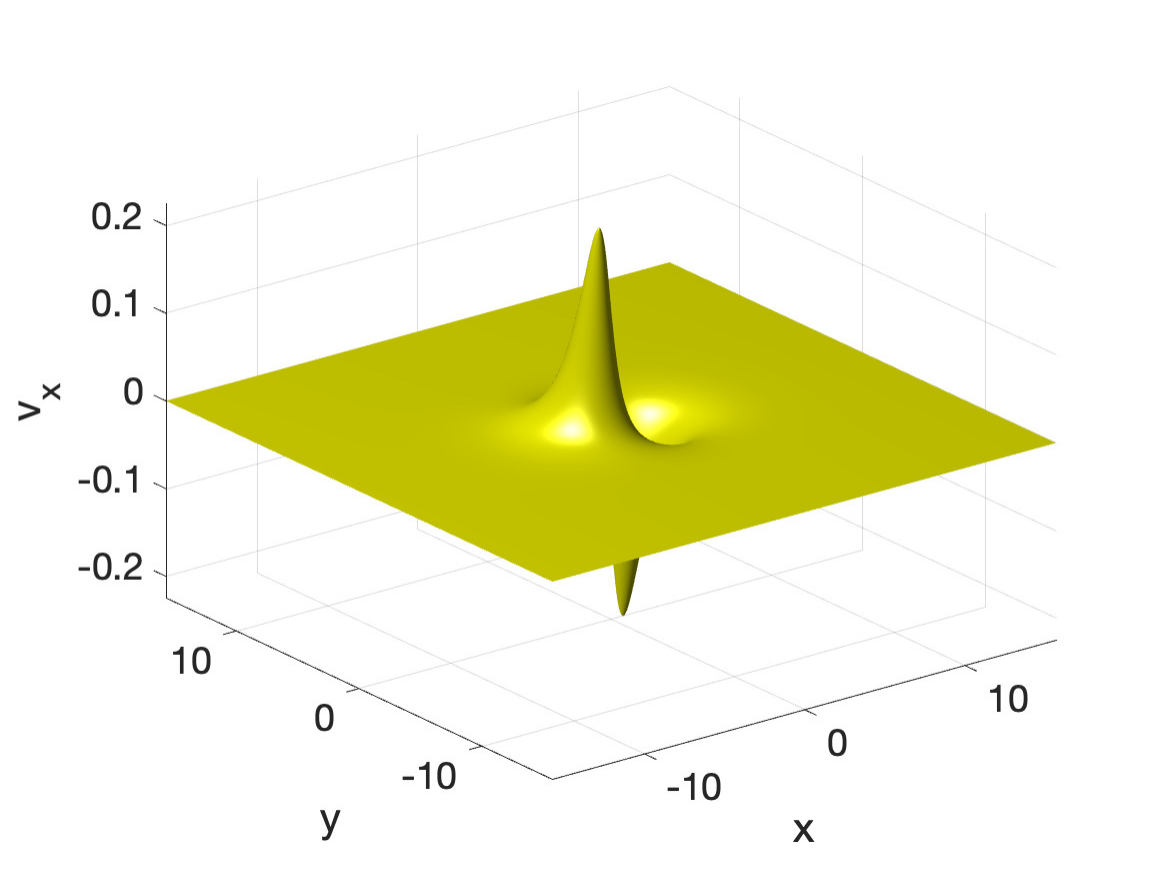}
 \includegraphics[width=0.49\textwidth]{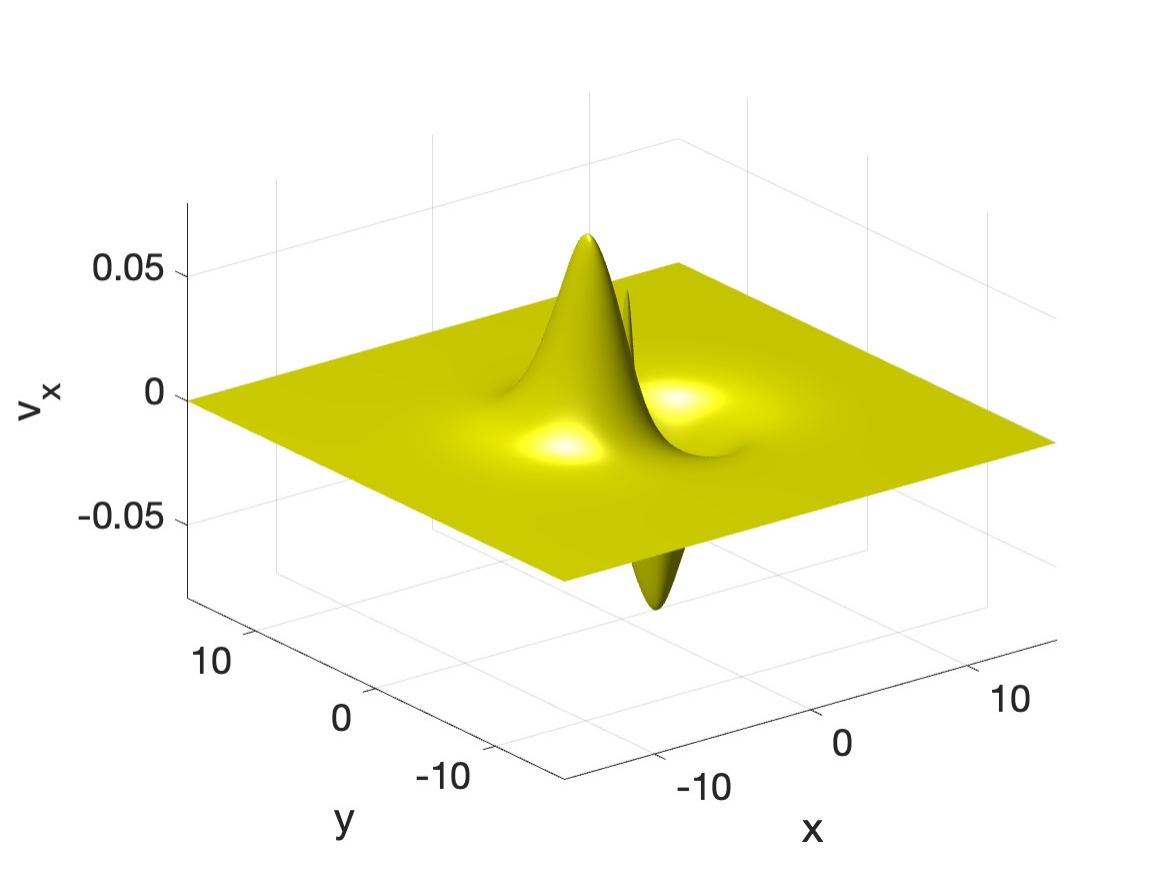}\\
 \includegraphics[width=0.49\textwidth]{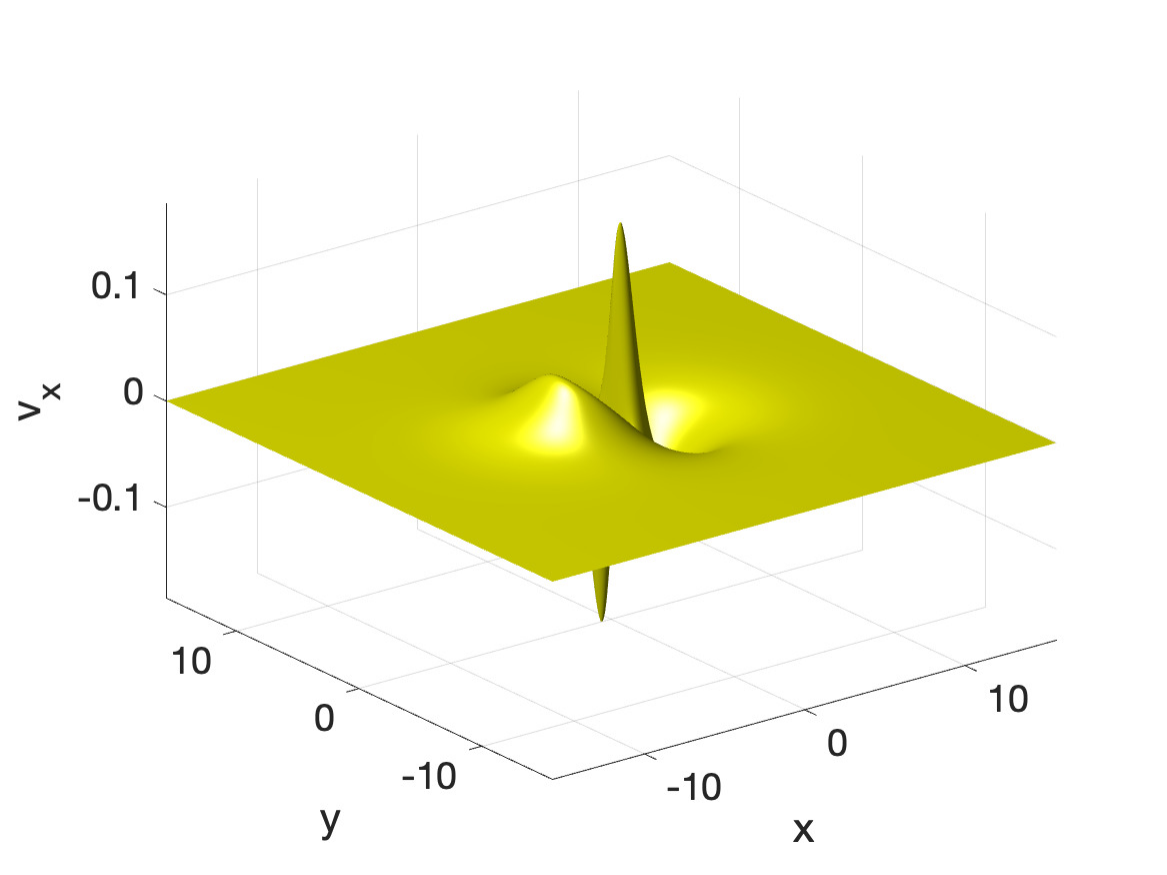}
 \includegraphics[width=0.49\textwidth]{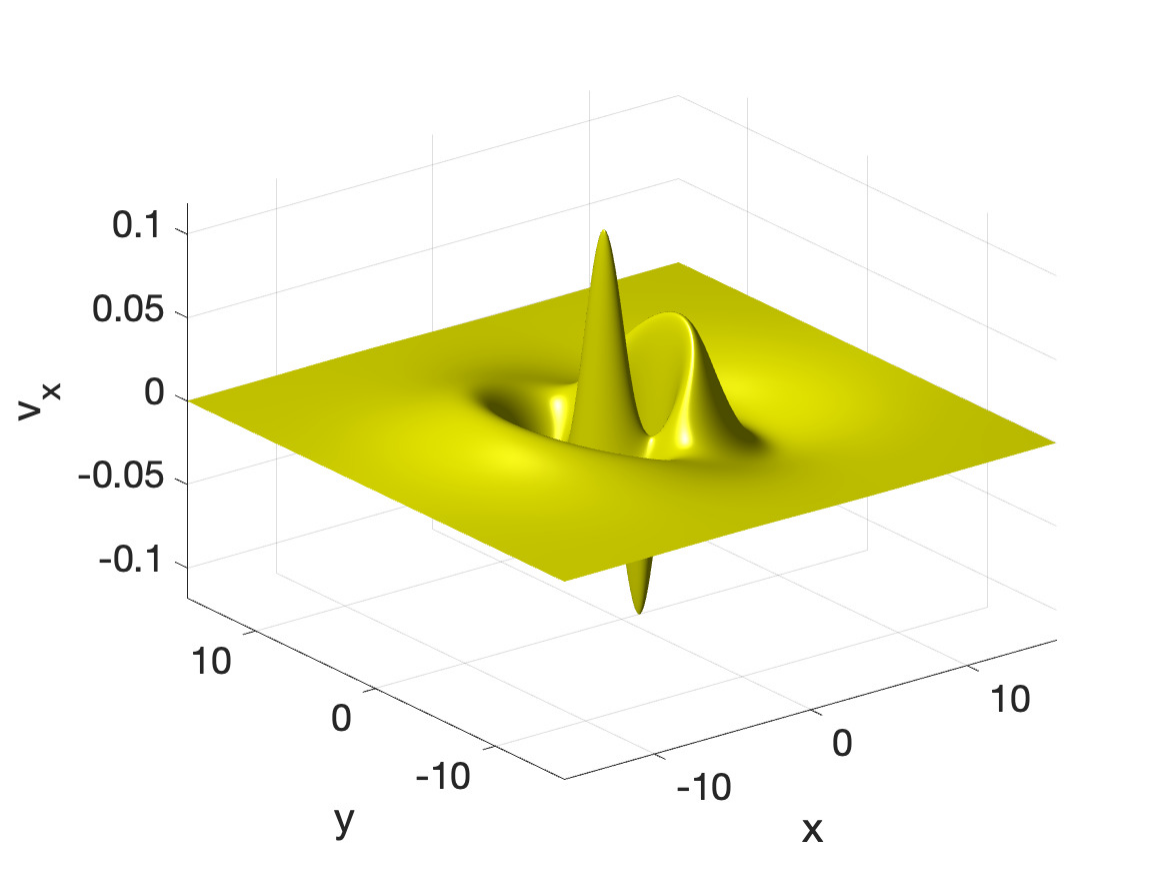}
 \caption{Solution $v_{x}$ to the AS system (\ref{2D}) for the initial data 
 (\ref{NCini}) with $\kappa=-0.9$ and $\alpha=1$, in the upper row on the left for 
 $t=2.5$,  on the right for 
 $t=4$, in the lower row on the left for $t=5.1$ and on the right 
for $t=10$.}
 \label{figm09gaussvx}
\end{figure}

The corresponding plots for $v_{x}$ can be seen in 
Fig.~\ref{figm09gaussvx}. The solution appears to be always smooth 
though strong gradients appear. Note that the solution is always well 
resolved both in space and in time. The plots for $v_{y}$ are very 
similar, just rotated by 90 degrees, and are therefore not shown. 

The situation is different if the non-cavitation condition is not 
satisfied for the initial data, for instance for initial data of the 
form (\ref{NCini}) with $\kappa=-1$, $\alpha=1$ where there is one point with 
$\eta(x,y,0)=0$. We consider $N_{x}=N_{y}=2^{12}$ DFT modes for $x\in 
3[\pi,\pi]$ and $x\in 
3[\pi,\pi]$. 
Nonetheless we find in the case $\kappa=-1$, $\alpha=1$ that the code breaks for 
$t=4.0857$ since the fitting of the DFT coefficients on the $k_{x}$ 
axis leads to a vanishing $\delta$ in (\ref{fourierasym}) 
indicating that a singularity 
appears on the $x$-axis. The fitted value for $\mu$ is $0.604$, i.e., 
compatible with a cusp, but as we mentioned we do not have the needed 
numerical resolution to make a precise statement. We can only report 
that there is an indication of a singularity which needs to be 
checked on more powerful computers. Note that the same 
values of the fitted parameters are obtained on the $y$-axis. This 
would indicate the formation of a cusp for $\eta$ of the form
\begin{equation}
	\eta\sim (x^{2}+y^{2})^{\mu/2}
	\label{cusp}
\end{equation}
at the critical time. 

\begin{figure}[htb!]
 \includegraphics[width=0.49\textwidth]{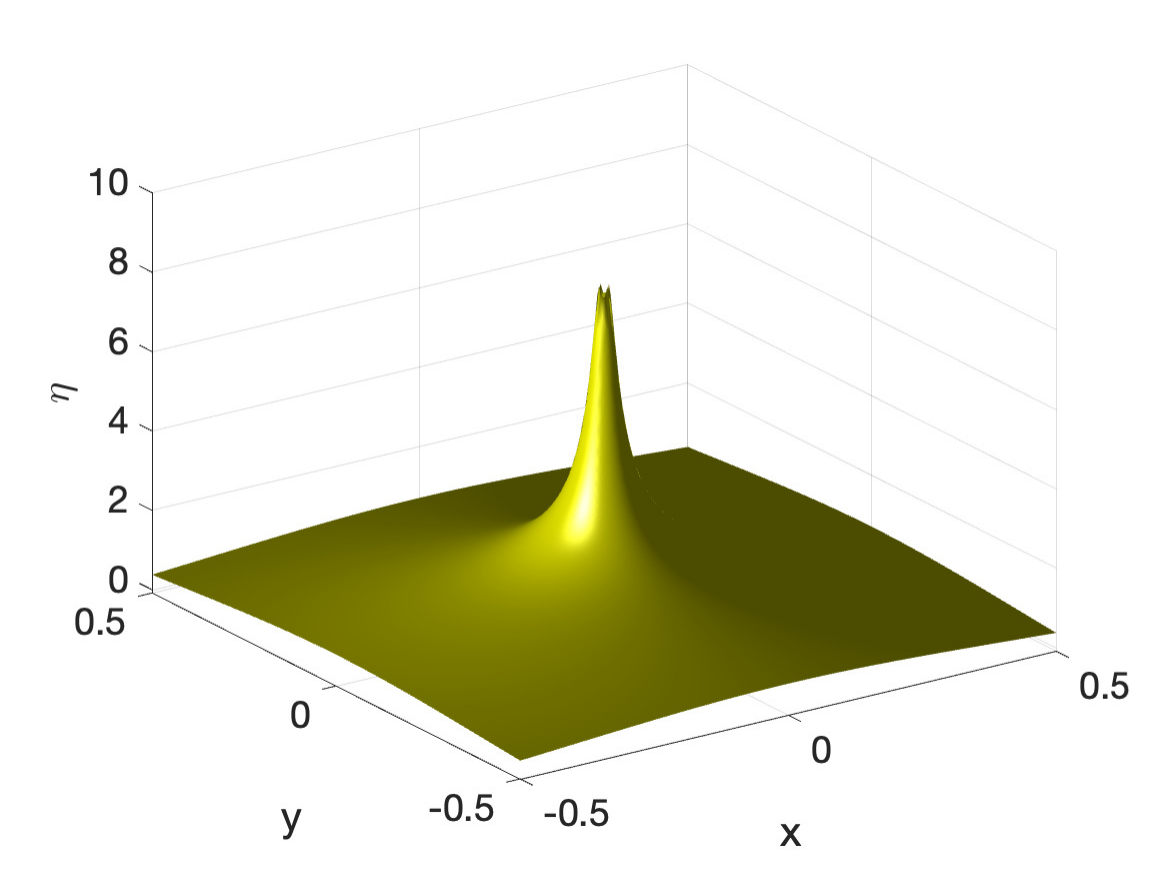}
 \includegraphics[width=0.49\textwidth]{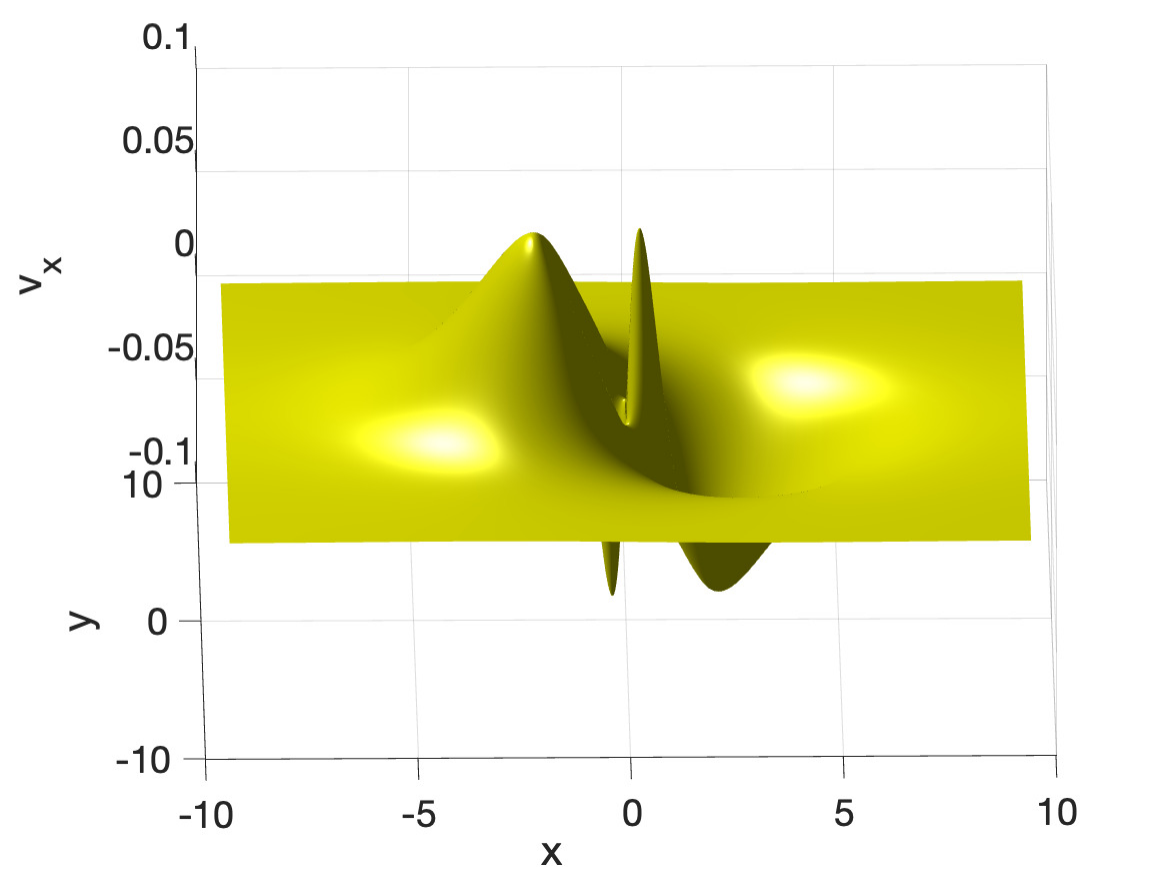}
 \caption{Solution  to the AS system (\ref{2D}) for the initial data 
 (\ref{NCini}) with $\kappa=-\alpha=-1$,  on the left $\eta$ for 
 $t=4.0857$,  on the right $v_{x}$ at the same time.}
 \label{figmgauss}
\end{figure}

Note further that the fitted values of the DFT 
coefficients for $v_{x}$ and $v_{y}$ appear to indicate that they are 
still regular at the critical time for $\eta$.
We show the solutions at the final recorded time in 
Fig.~\ref{figmgauss}. The behavior of $\eta$ near the origin is 
complicated, a strong peak with a caldera.

\begin{figure}[htb!]
 \includegraphics[width=0.49\textwidth]{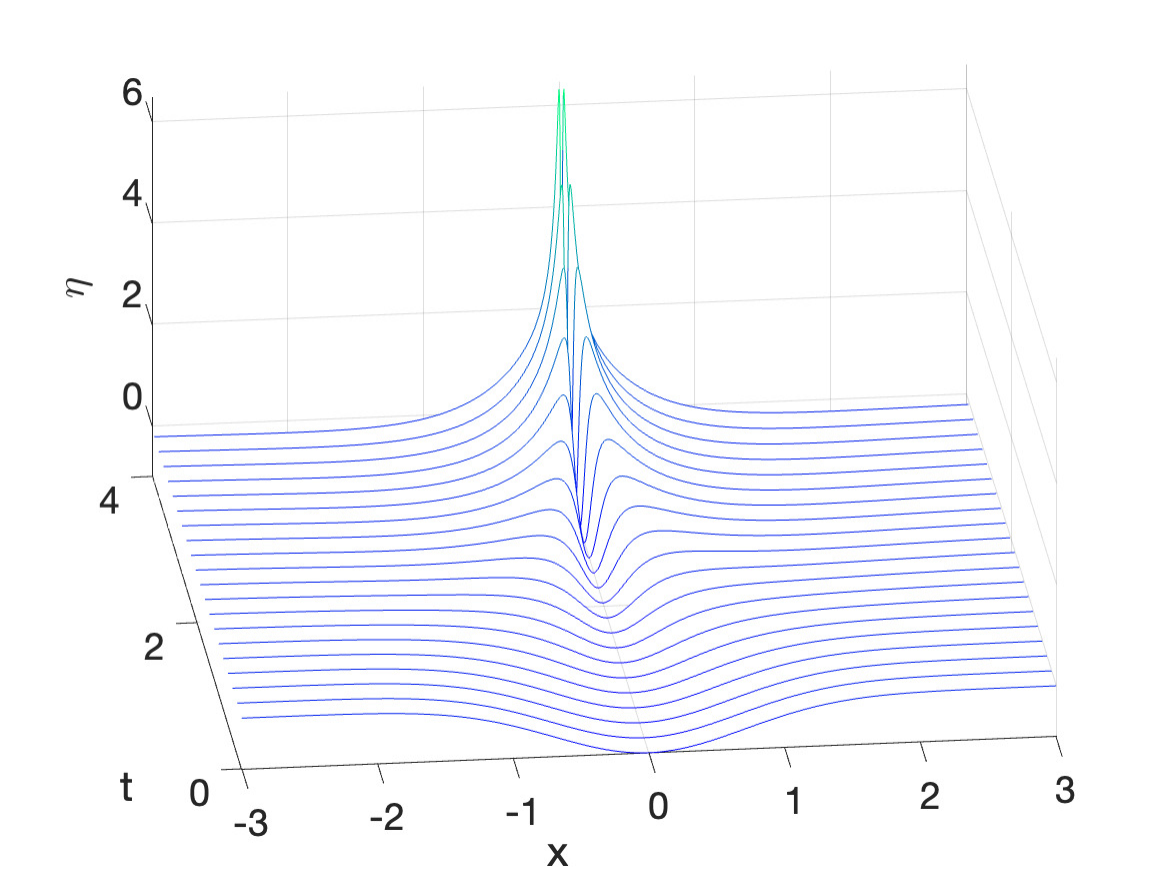}
 \includegraphics[width=0.49\textwidth]{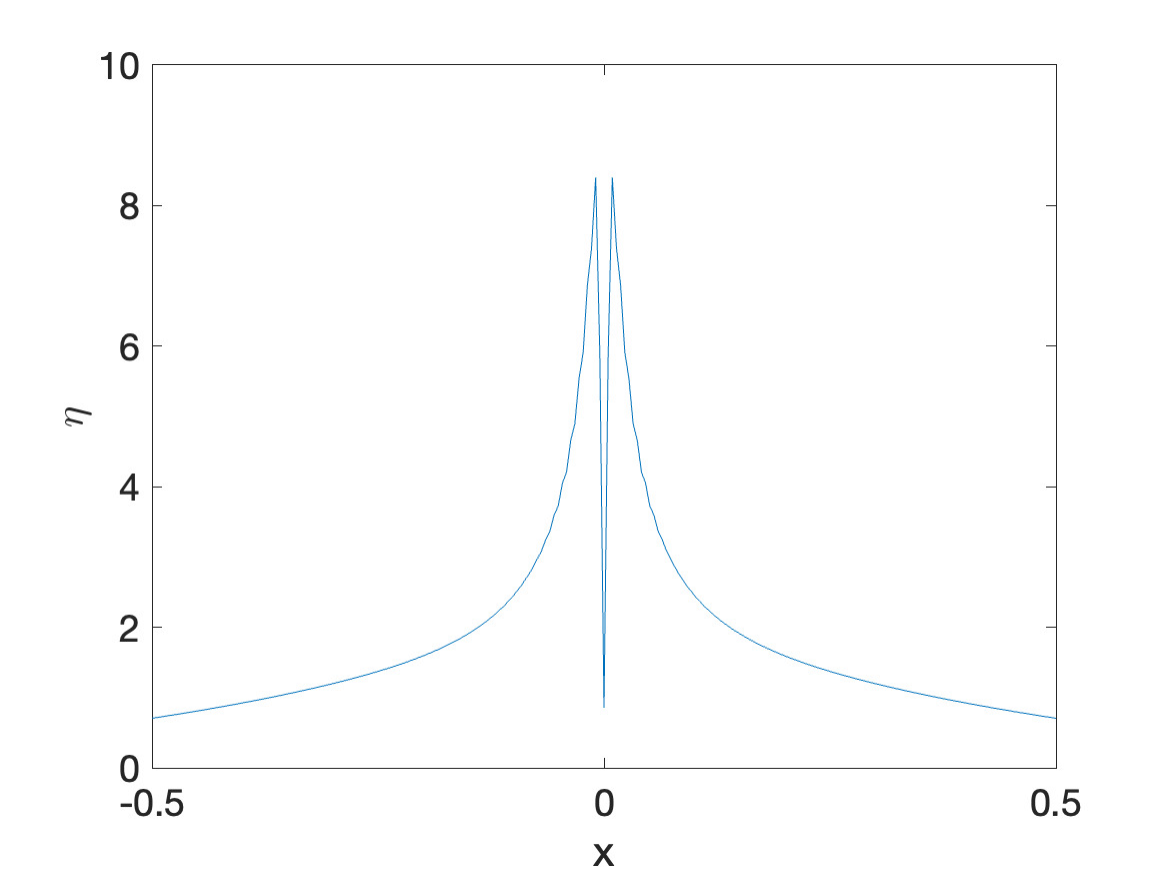}\\
 \includegraphics[width=0.49\textwidth]{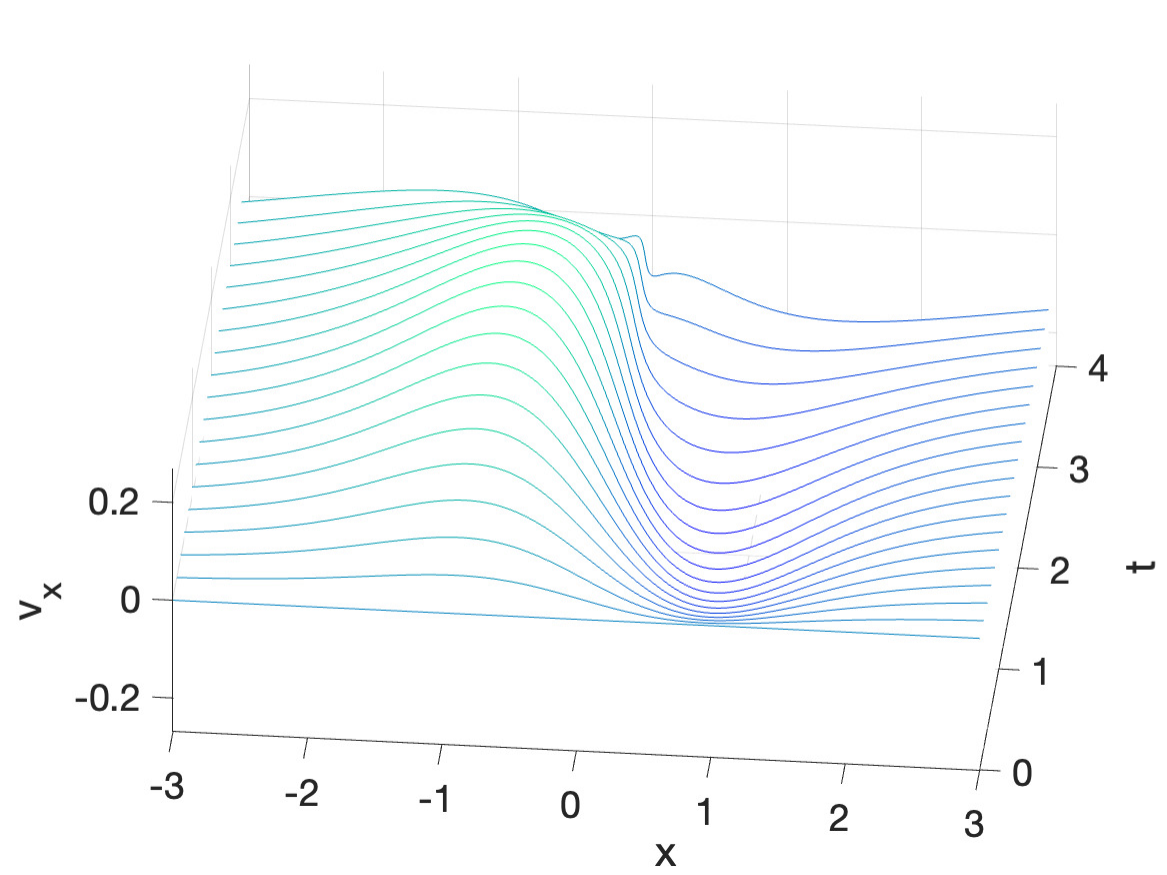}
 \includegraphics[width=0.49\textwidth]{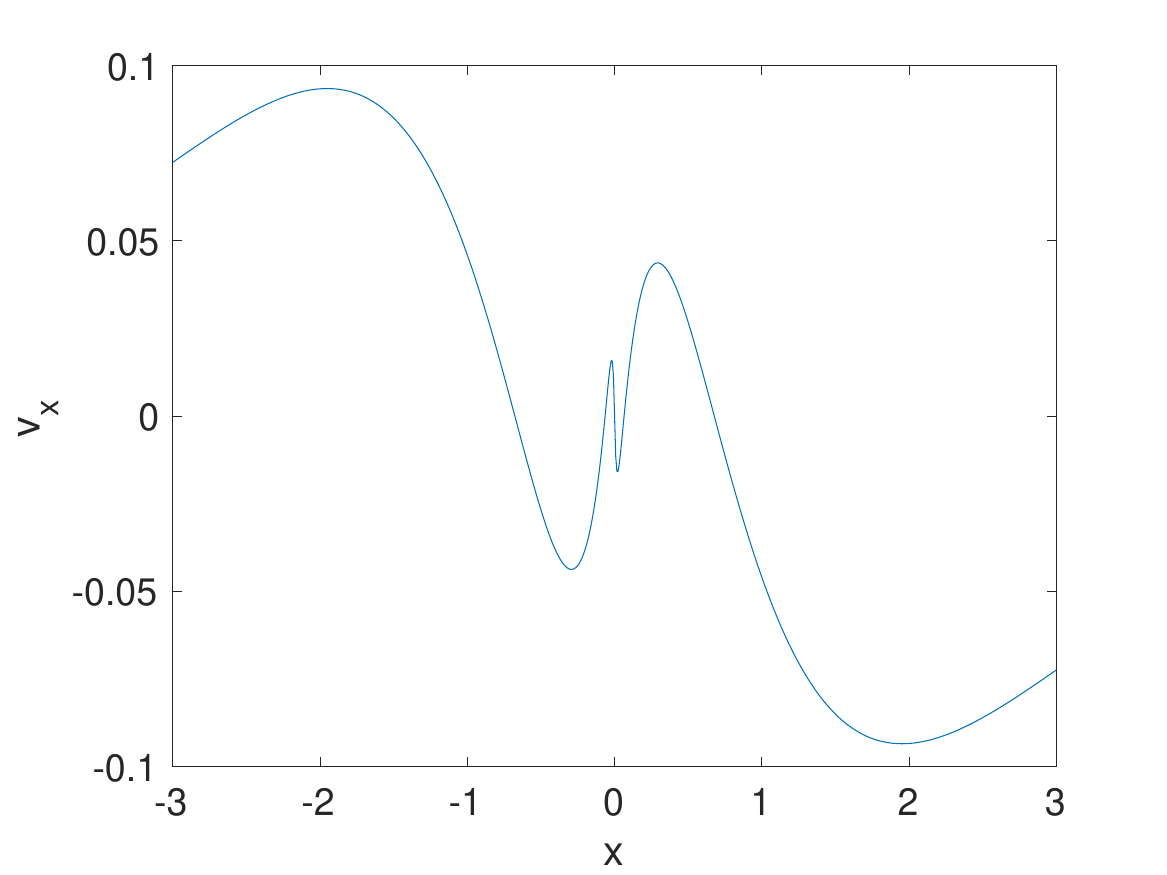}
 \caption{Solution to the AS system (\ref{2D}) for the initial data 
 (\ref{NCini}) with $\kappa=-\alpha=-1$ on the $x$-axis, in the upper row on 
 the left in dependence of time,   on the right at the critical time 
 in a close-up, 
in the lower row on the left $v_{x}$ and on the right 
a close-up at the critical time.}
 \label{figmgaussaxis}
\end{figure}

It appears best to visualize the behavior of the solution on the 
$x$-axis which is done in Fig.~\ref{figmgaussaxis}. It can be seen 
that the initial depression for $\eta$ is closed by two  
peaked maxima that are approaching each other, but that a gap remains 
in between these two 
peaks. The velocity $v_{x}$ is seen in the lower row of the same 
figure. It seems to develop a cusp near the strongly peaked minimum 
of $\eta$. 

\begin{figure}[htb!]
 \includegraphics[width=0.49\textwidth]{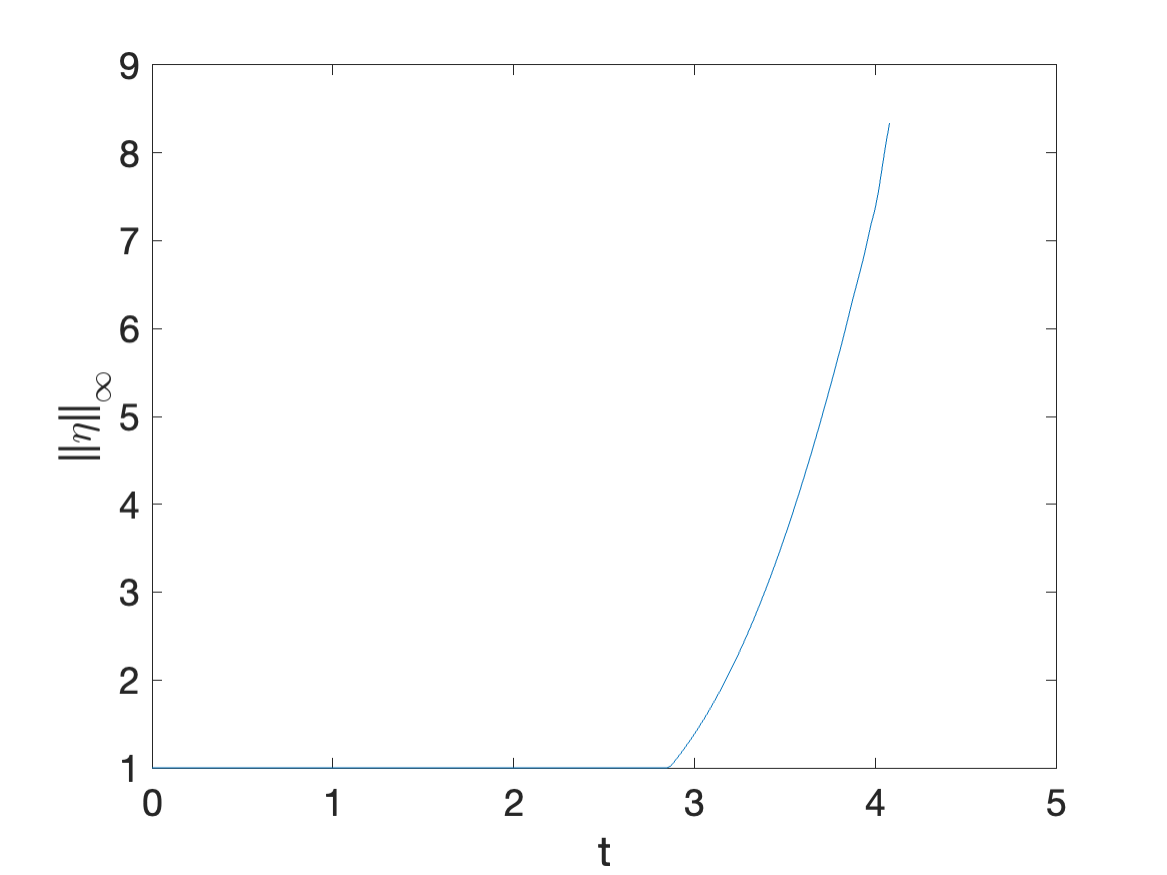}
 \includegraphics[width=0.49\textwidth]{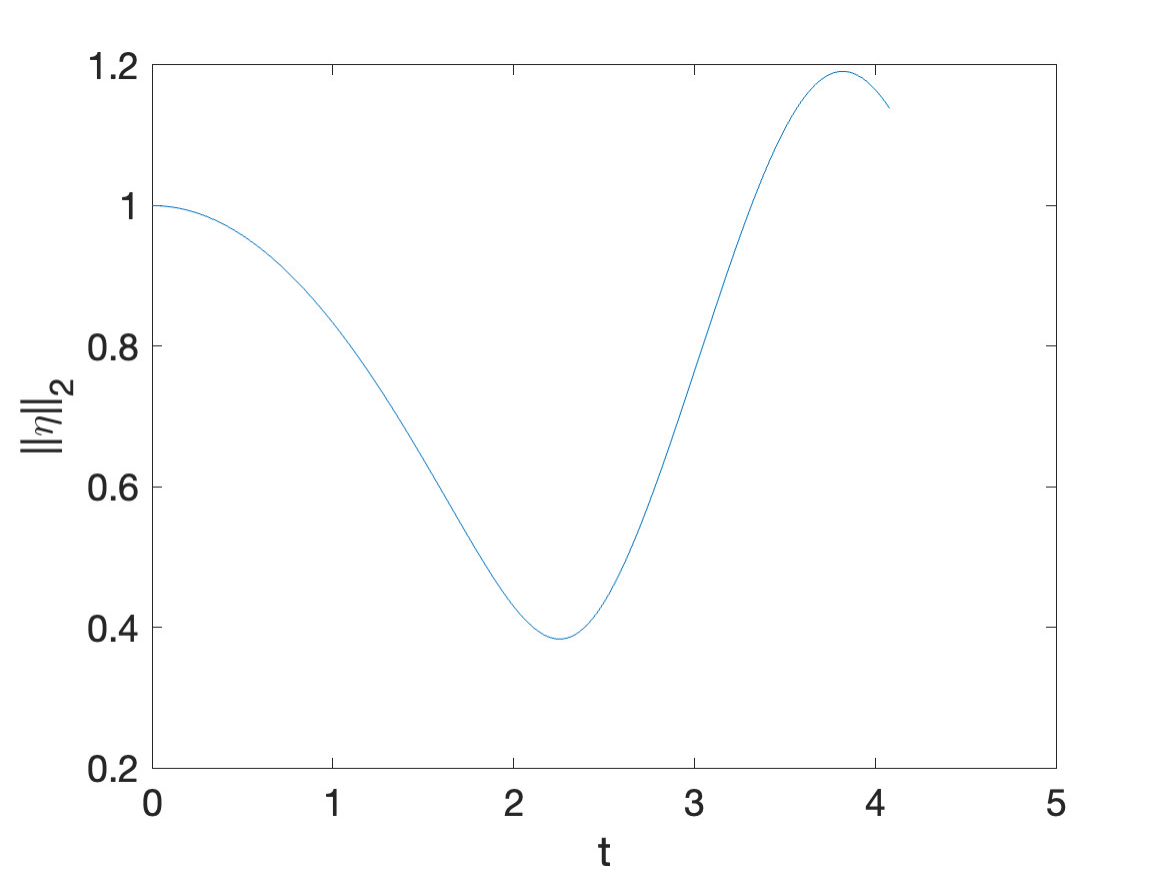}\\
 \includegraphics[width=0.49\textwidth]{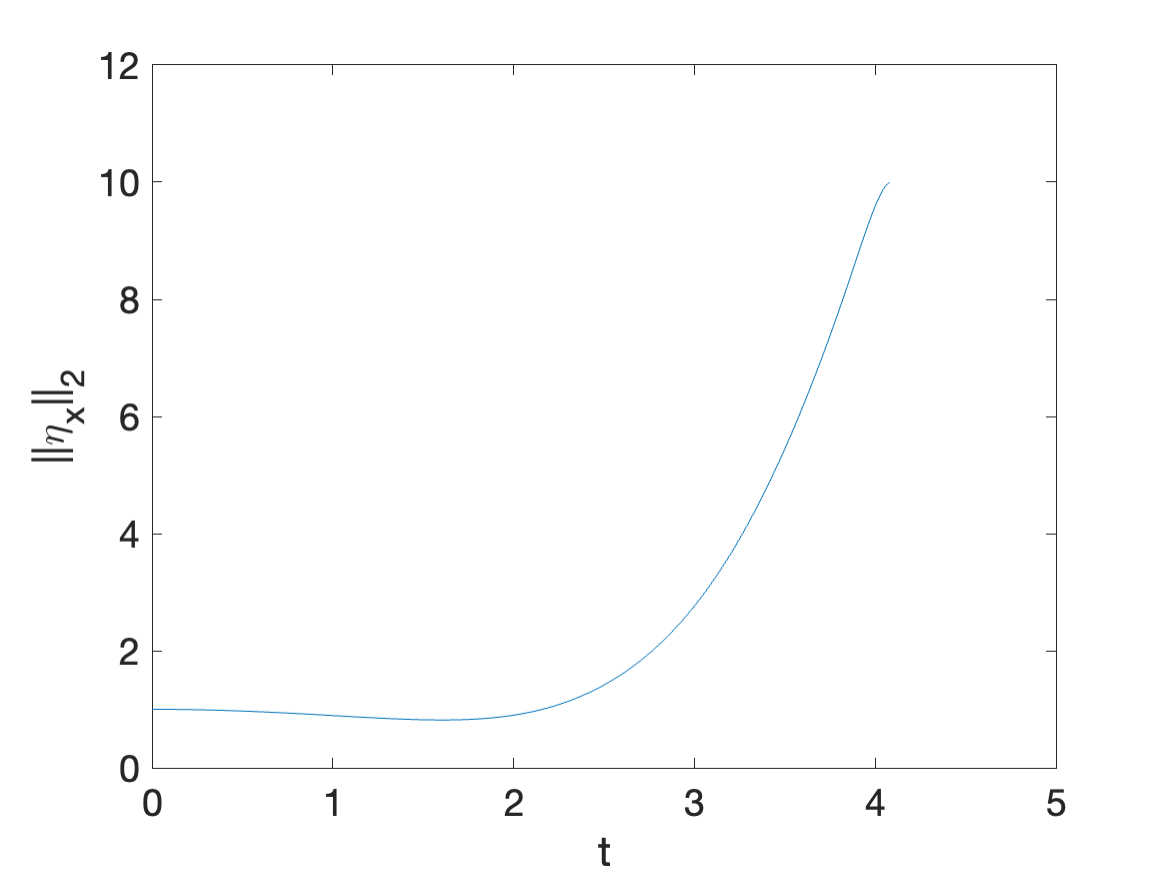}
 \includegraphics[width=0.49\textwidth]{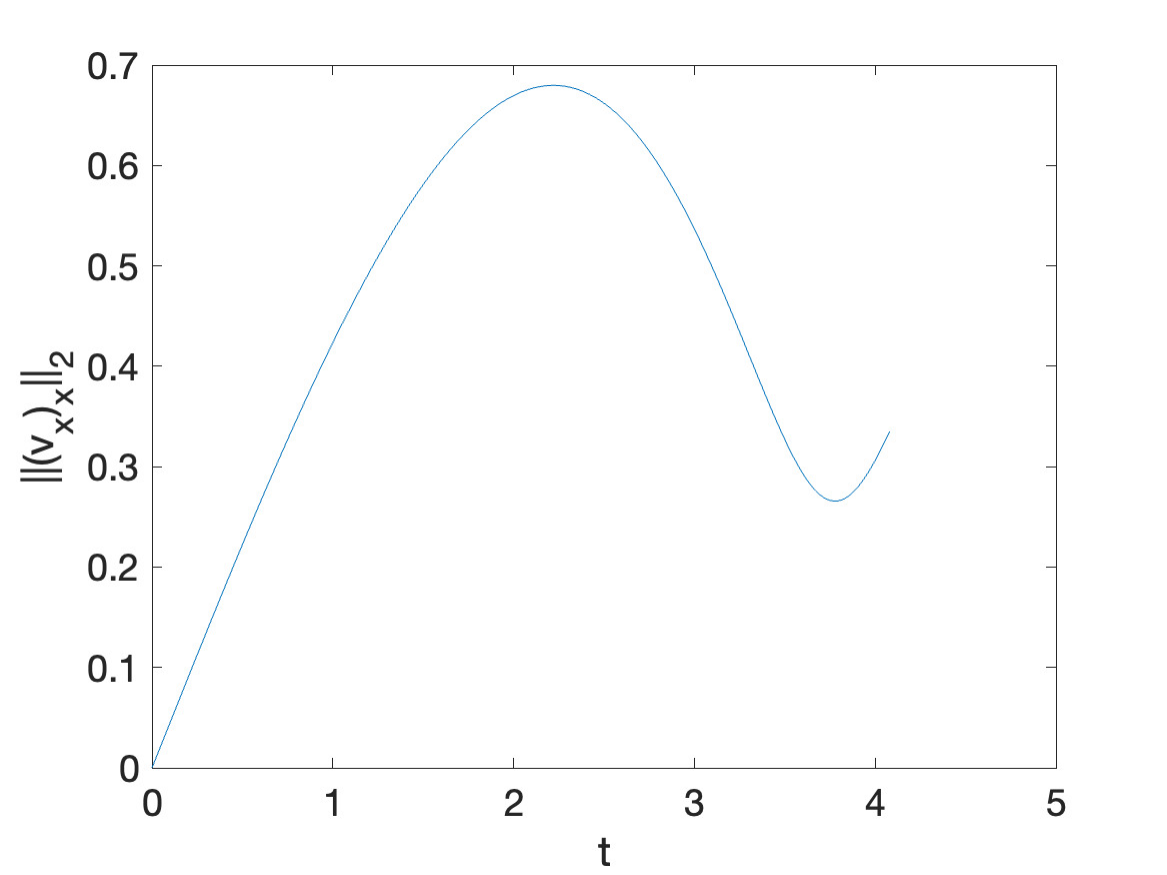}
 \caption{Various norms of the solution to the AS system (\ref{2D}) for the initial data 
 (\ref{NCini}) with $\kappa=-\alpha=-1$, in the upper row 
 the $L^{\infty}$ and the $L^{2}$ norm of $\eta$ (normalized to 1 at 
 $t=0$), in the lower row 
 the $L^{2}$ norms of the $x$-derivative of $\eta$  (normalized to 1 at 
 $t=0$) and of $v_{x}$.}
 \label{figmgaussnorms}
\end{figure}

As in the 1D case, we explore the possibly singular behavior via 
certain norms of the solution. In Fig.~\ref{figmgaussnorms}, the 
$L^{\infty}$ norm of $\eta$ appears to grow as can be seen on the 
left of the upper row. However as in the 1D case, a possible 
singularity only appears as a cusp at the minimum of the solution. 
This is confirmed by the $L^{2}$ norm of $\eta$ (higher $L^{p}$ norms are very 
similar) on the right. The $L^{2}$ norm of the $x$-derivative of 
$\eta$ can be seen in the lower row. It appears to grow, but we 
clearly do not have the needed resolution to reach the potential cusp 
singularity. The corresponding norm for $v_{x}$ does not indicate any 
singularity formation at the considered times. 

\begin{figure}[htb!]
 \includegraphics[width=0.49\textwidth]{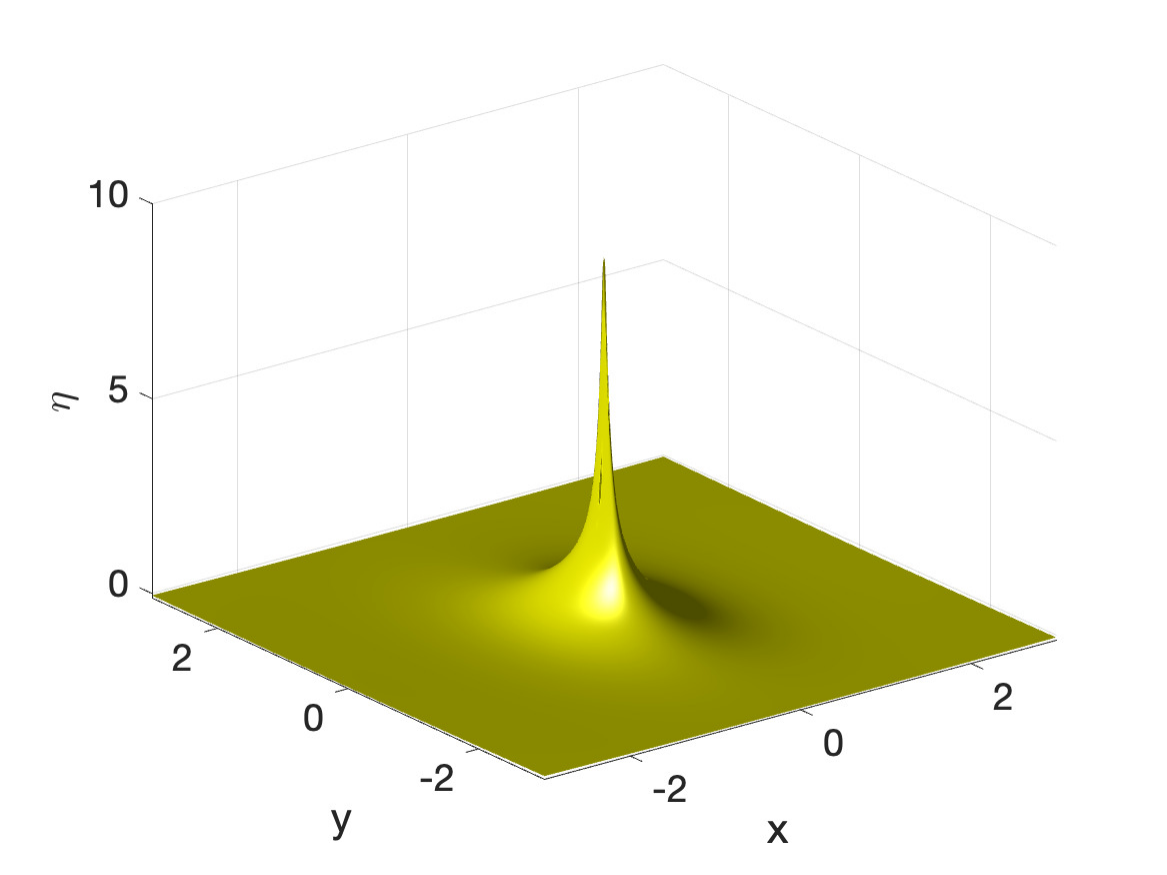}
 \includegraphics[width=0.49\textwidth]{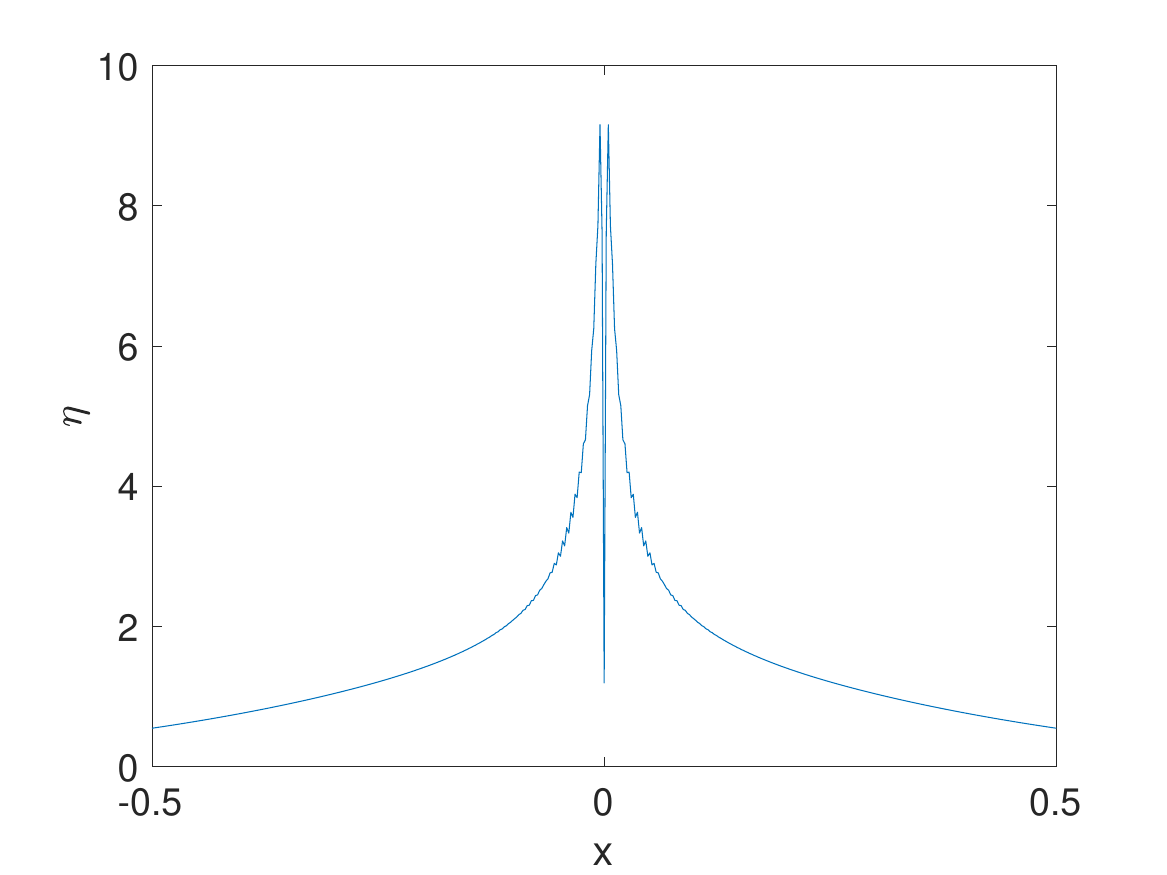}\\
 \includegraphics[width=0.49\textwidth]{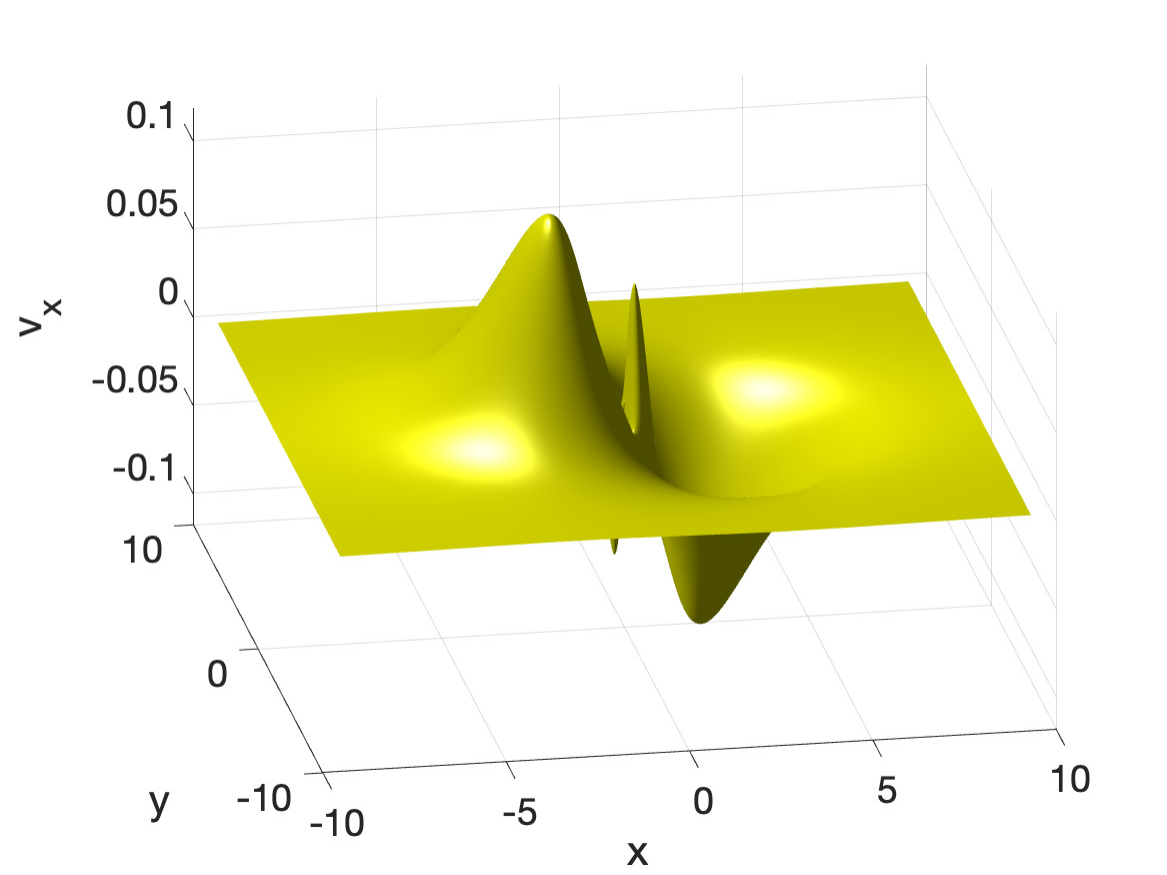}
 \includegraphics[width=0.49\textwidth]{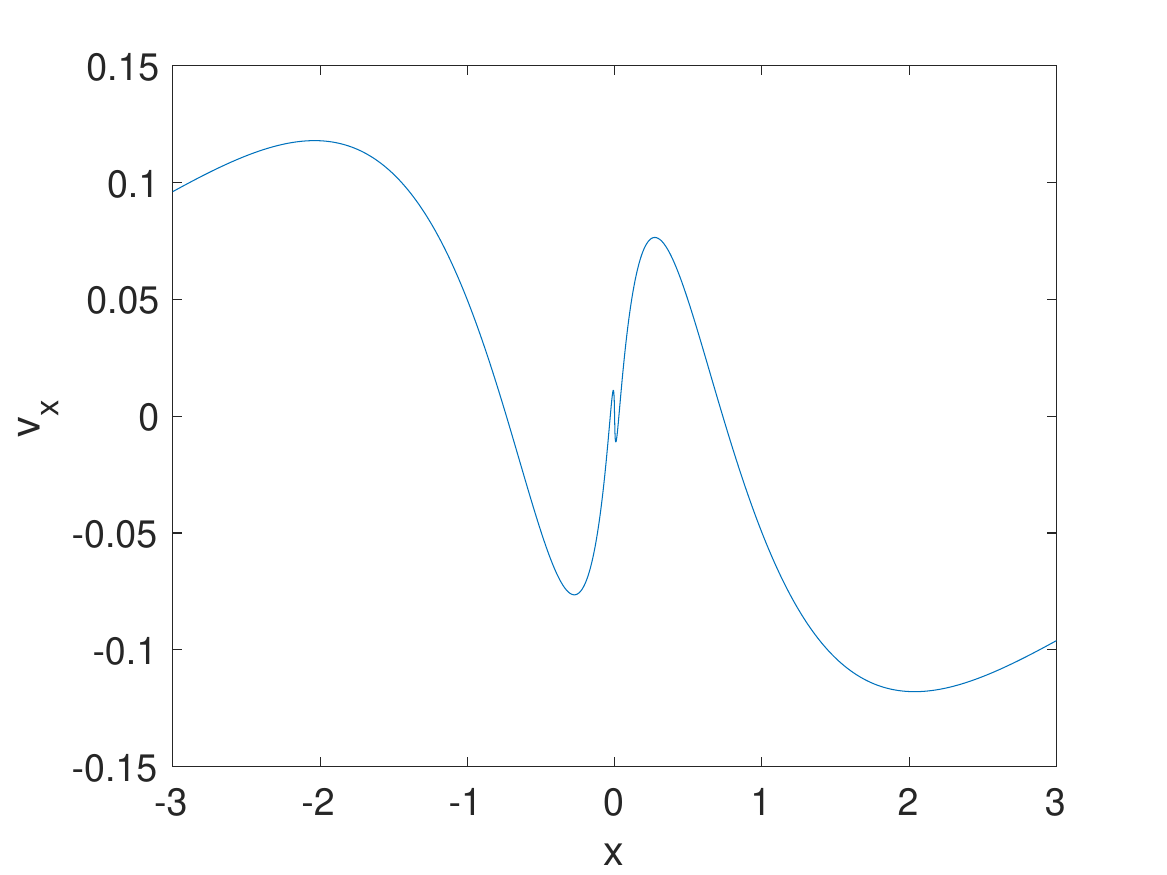}
 \caption{Solution to the AS system (\ref{2D}) for the initial data 
 (\ref{NCini}) with $\kappa=-1$ and $\alpha$=0.5. In the upper row on 
 the left the solution $\eta$ at the critical time, on the right a 
 close-up of the solution on the $x$-axis. In the lower row the 
 solution $v_{x}$ at the critical time on the left and on the right a 
 close-up on the $x$-axis.}
 \label{figmgaussy2}
\end{figure}

For completeness we consider also a situation without radial 
symmetry, initial data of the form (\ref{NCini}) with $\kappa=-1$ and 
$\alpha=0.5$. We use the same computational domain as above, but 
$N_{x}=2^{13}$, $N_{y}=2^{11}$ DFT modes. The code breaks for 
$t=4.0857$ since the fitted $\delta$ on the  $k_{x}$-axis vanishes. 
We show the 
solution at the final time in Fig.~\ref{figmgaussy2}.

\begin{figure}[htb!]
 \includegraphics[width=0.49\textwidth]{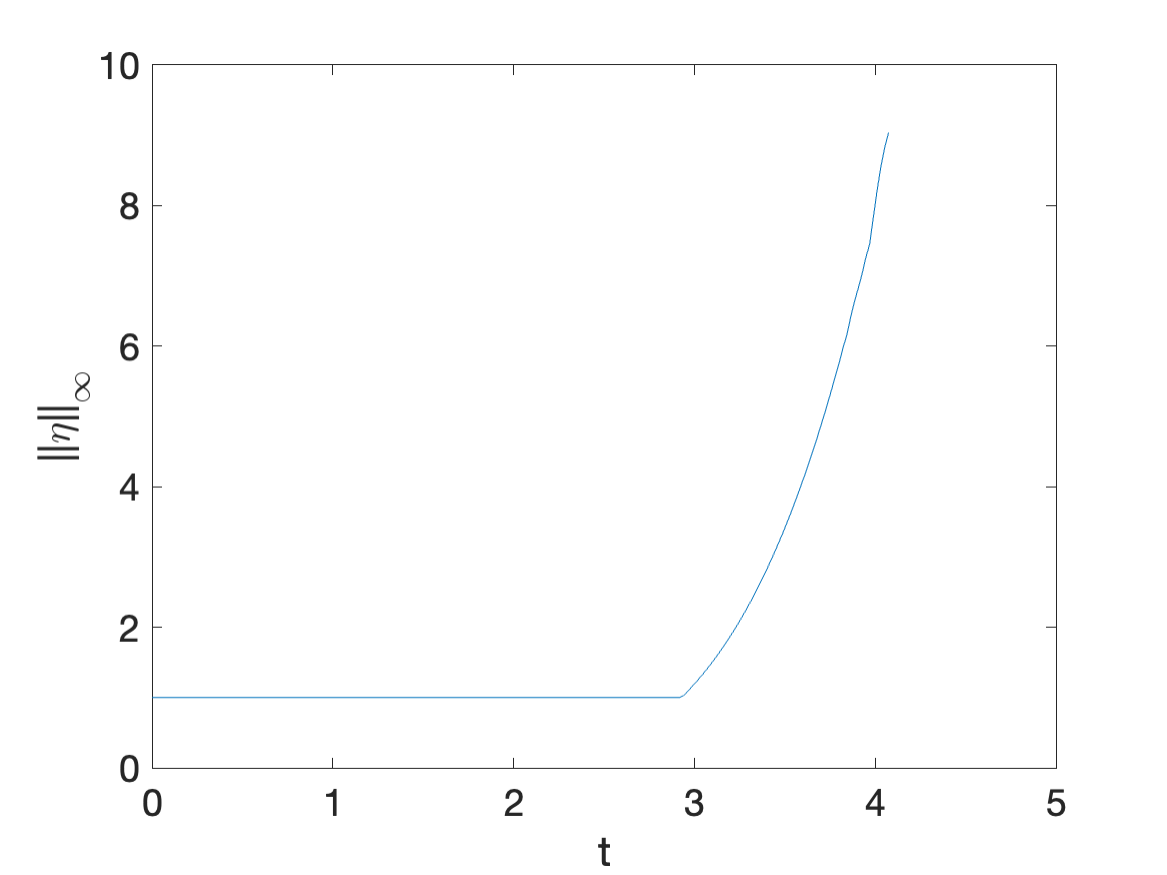}
 \includegraphics[width=0.49\textwidth]{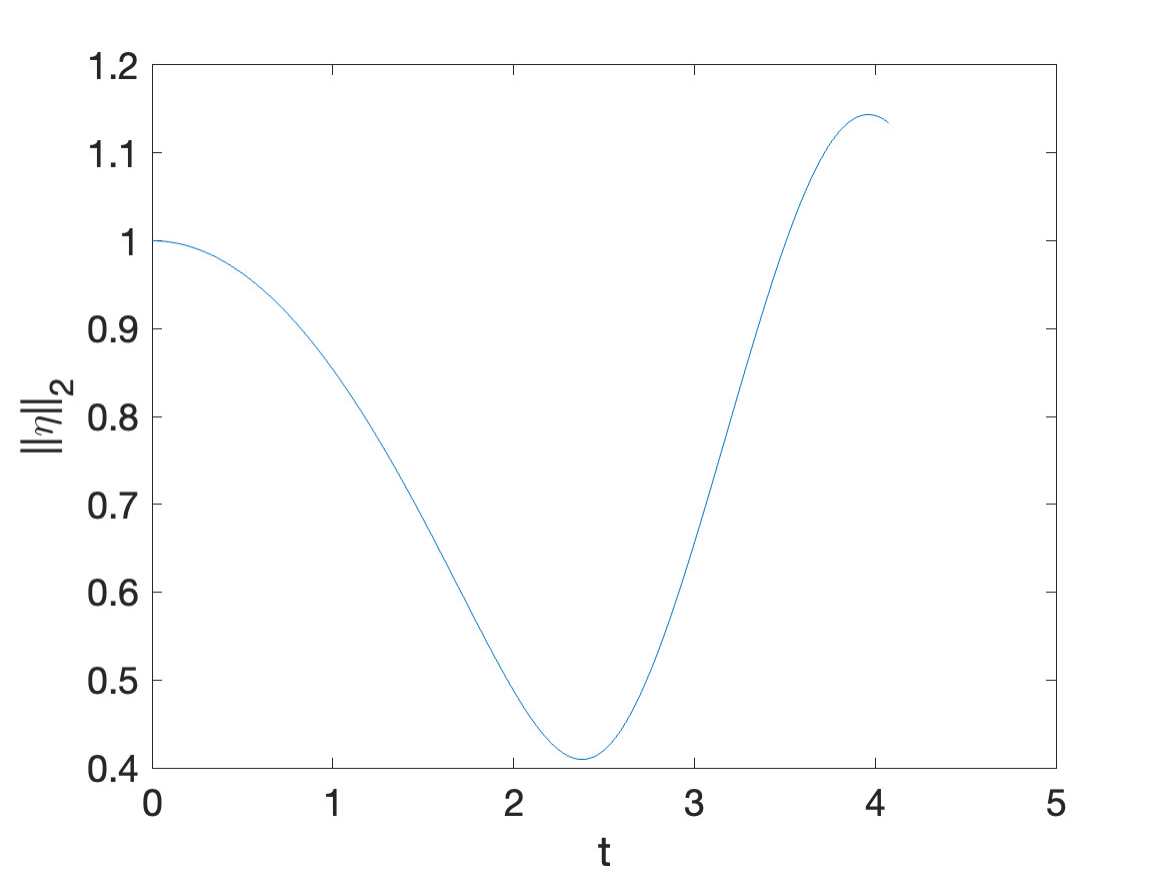}\\
 \includegraphics[width=0.49\textwidth]{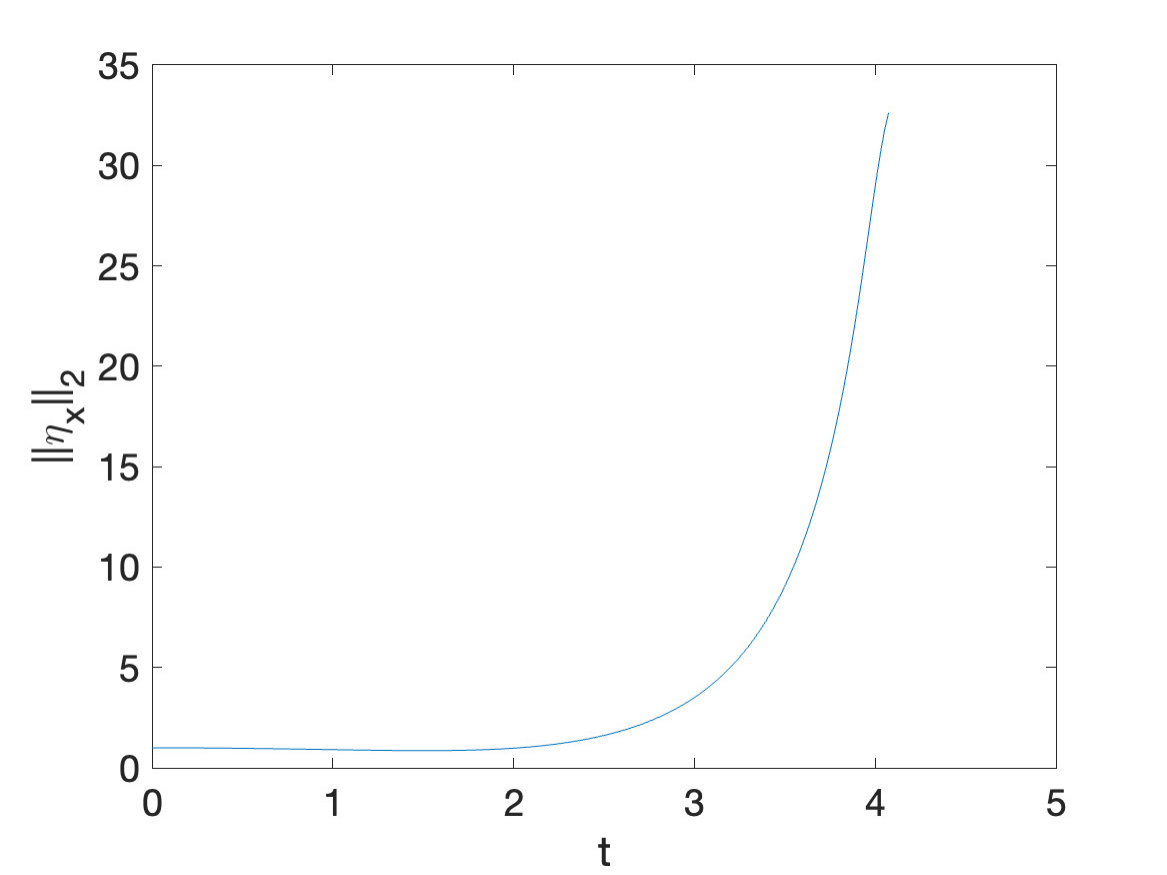}
 \includegraphics[width=0.49\textwidth]{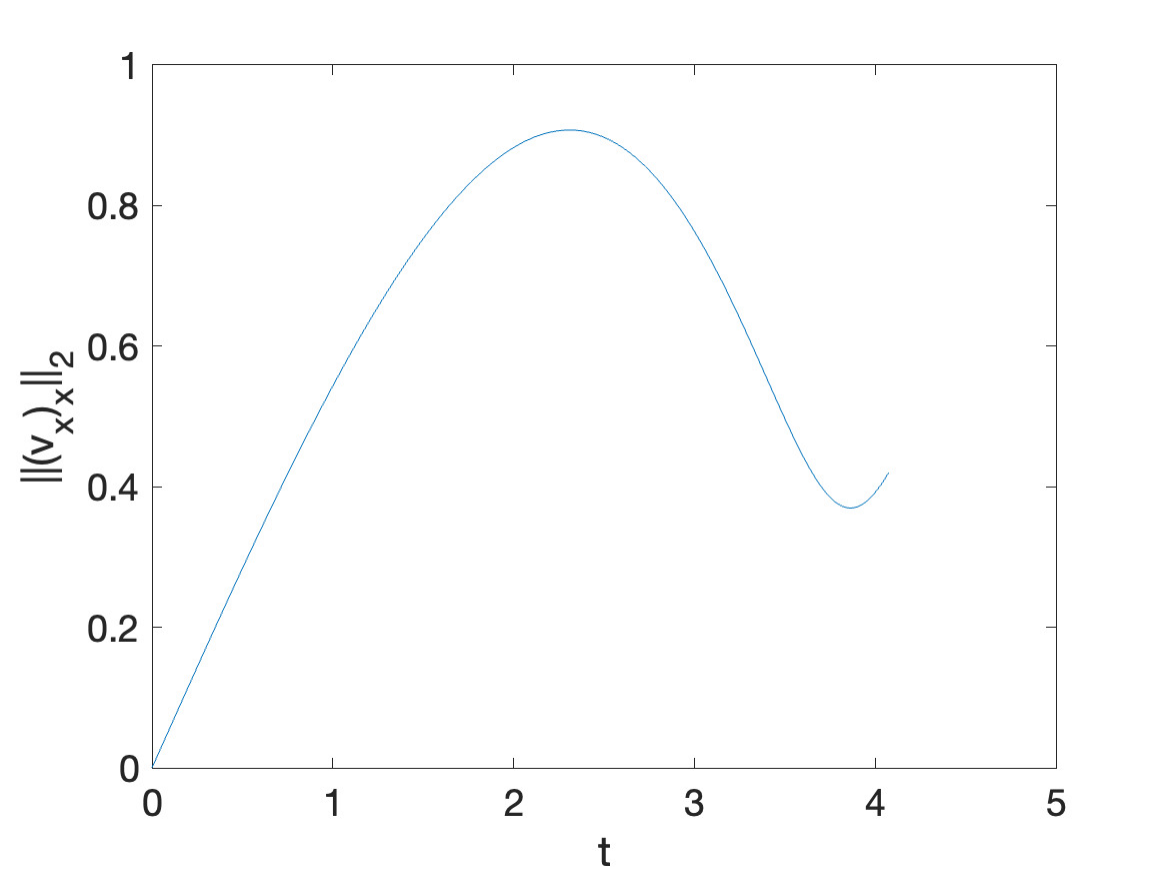}
 \caption{Various norms of the solution to the AS system (\ref{2D}) for the initial data 
 (\ref{NCini}) with $\kappa=-1$ and $\alpha=0.5$, in the upper row 
 the $L^{\infty}$ and the $L^{2}$ norm of $\eta$ (normalized to 1 at 
 $t=0$), in the lower row 
 the $L^{2}$ norms of the $x$-derivative of $\eta$  (normalized to 1 at 
 $t=0$) and of $v_{x}$.}
 \label{figmgaussy2norms}
\end{figure}

The corresponding  fitted value for $\mu\sim-.07$ is even negative, 
but as in the 1D case this is not reliable. Therefore we show in 
Fig.~\ref{figmgaussy2norms} some norms of the solution. 
The $L^{\infty}$ norm of $\eta$ grows as 
in the radially symmetric case, but the $L^{2}$ norm stays finite. 
The $L^{2}$ norm of the $x$-derivative of $\eta$ grows, but more 
slowly than in 1D which once more indicates that we lack resolution to get 
closer to the potential singularity. The $L^{2}$ norm of the 
$x$-derivative of $v_{x}$ stays finite for the studied times. 

\section{Localised initial data}
In this section we study the time evolution of localized initial 
data of the form (\ref{NCini}) with positive $\kappa$. As in the 
study of the transverse stability of the line solitary waves, no 
stable structures localized in two dimensions are observed. 

\begin{figure}[htb!]
  \includegraphics[width=0.49\textwidth]{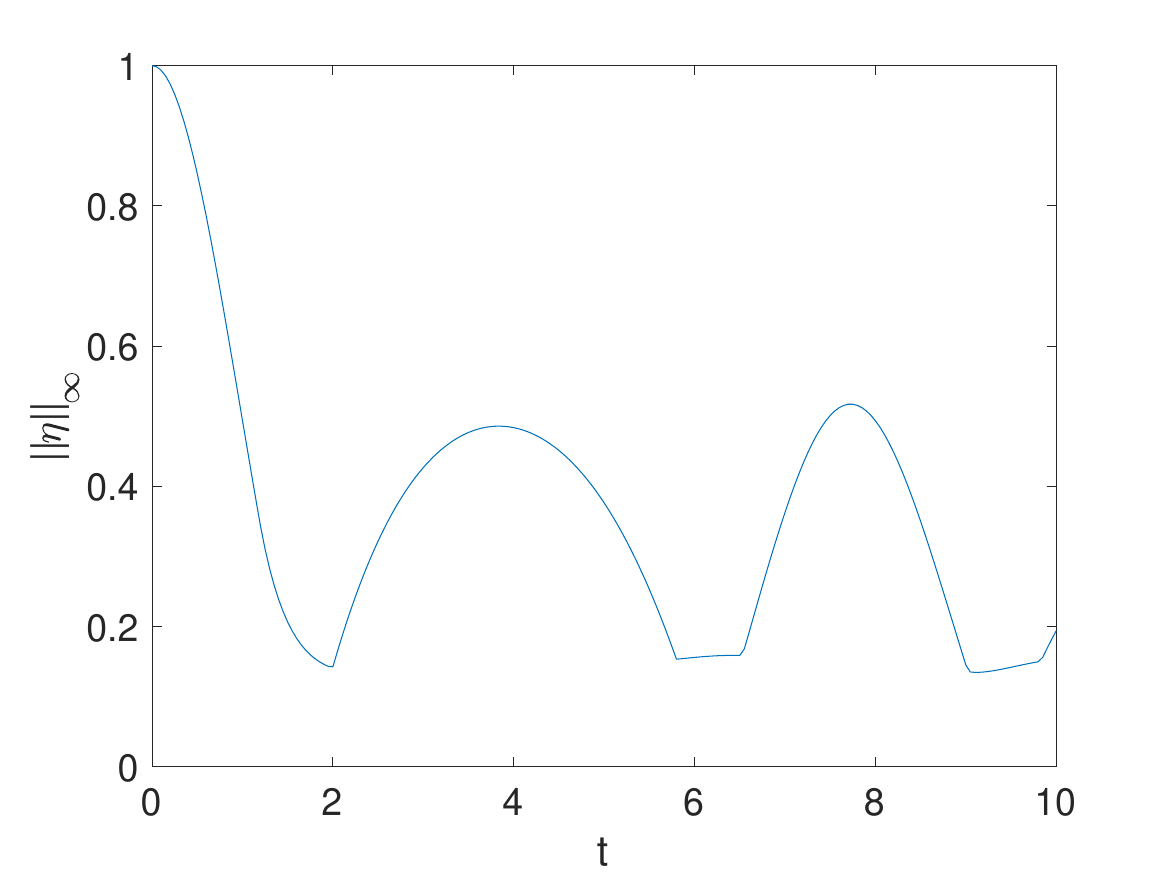}
  \includegraphics[width=0.49\textwidth]{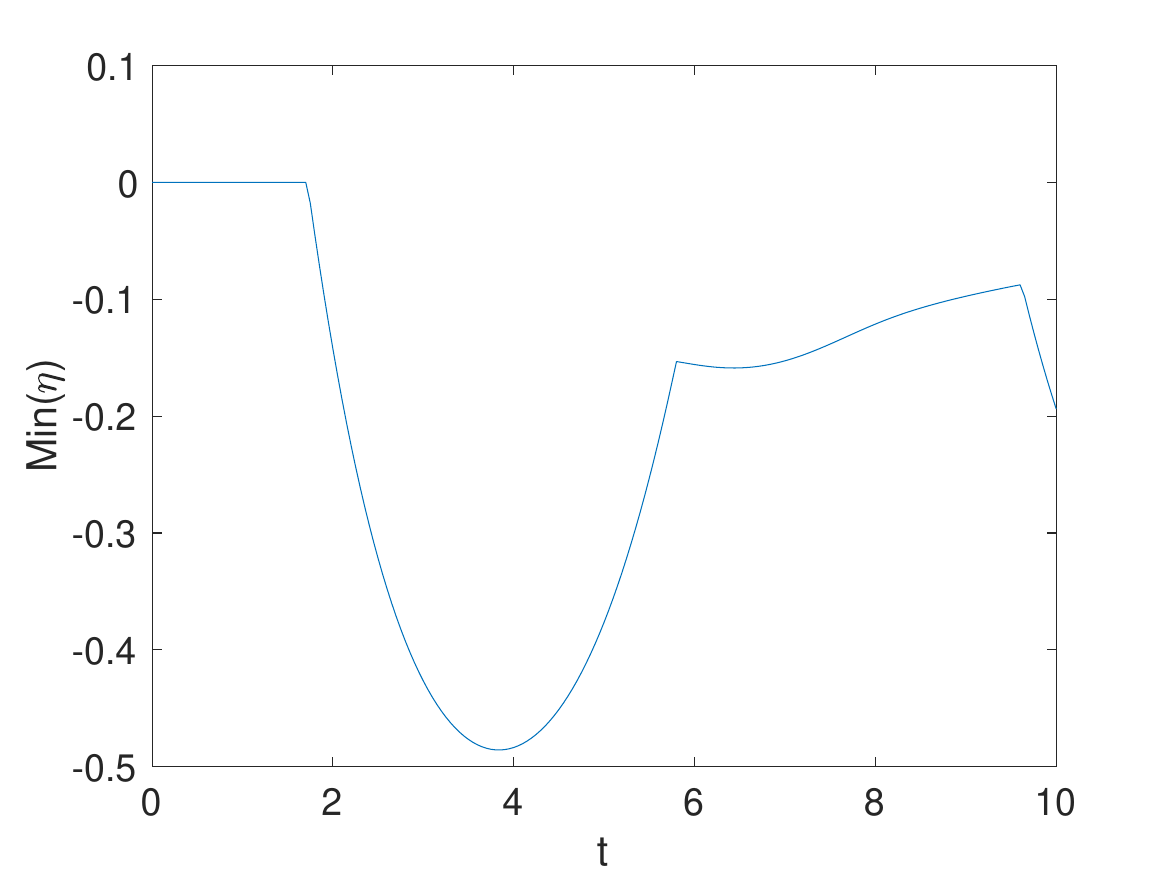} 
  \caption{Solution to the AS system (\ref{2D}) for the initial data 
  (\ref{NCini}) with $\kappa=\alpha=1$, on the left the $L^{\infty}$ norm of $\eta$, on the 
  right its global minimum, both in dependence of time.}
  \label{figgaussmax}
 \end{figure}

First we look at initial data with a small positive elevation, say 
$\kappa=1$ in (\ref{NCini}). The same numerical parameters as in the 
previous section are applied. In Fig.~\ref{figgaussmax} we show the 
global minimum of the solution $\eta$ and its $L^{\infty}$ norm in 
dependence of time. It can be seen that an annular structure forms 
near the original maximum, but that the observed minimum reaches only 
values of the order of $-0.5$ before a new peak is formed. Thus one 
observes a similar dynamics as in the previous section. There is no 
indication of a stable structure emerging from these initial 
conditions. 
 
\begin{figure}[htb!]
  \includegraphics[width=0.49\textwidth]{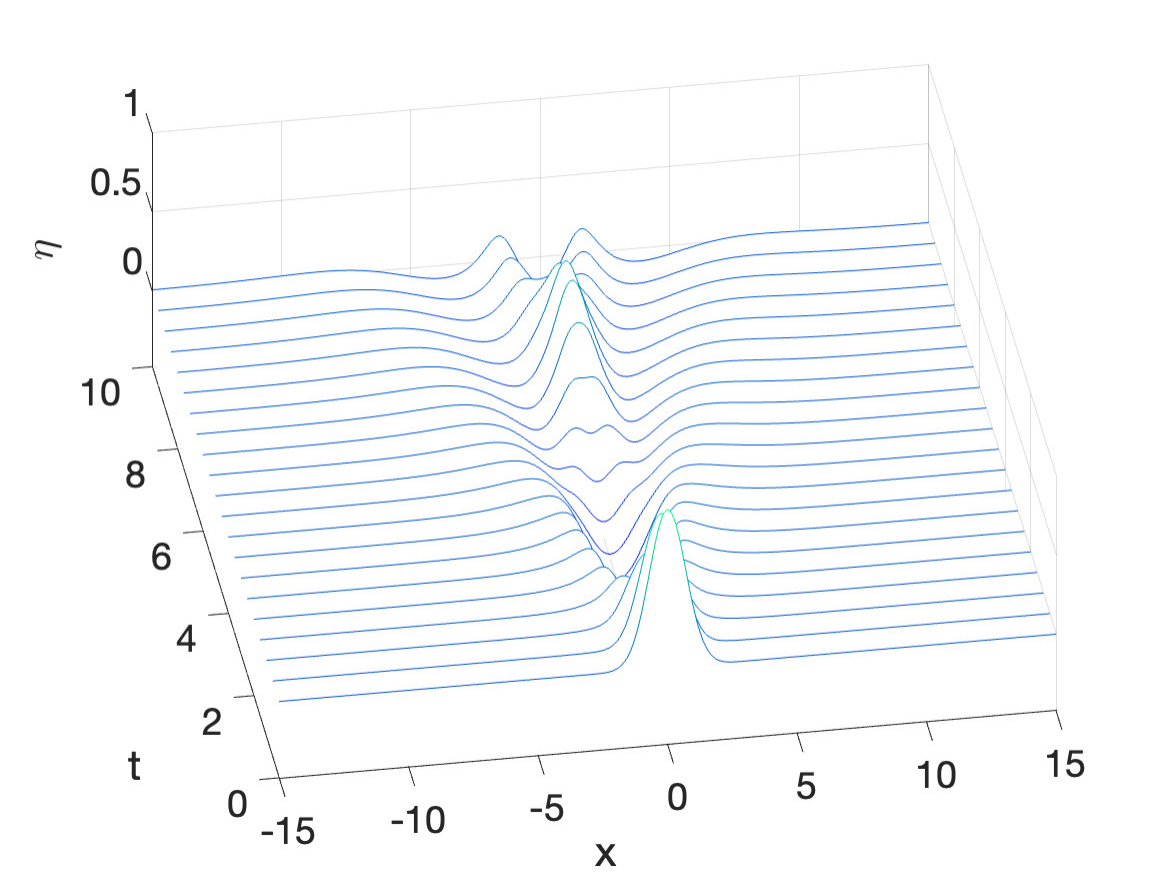}
  \includegraphics[width=0.49\textwidth]{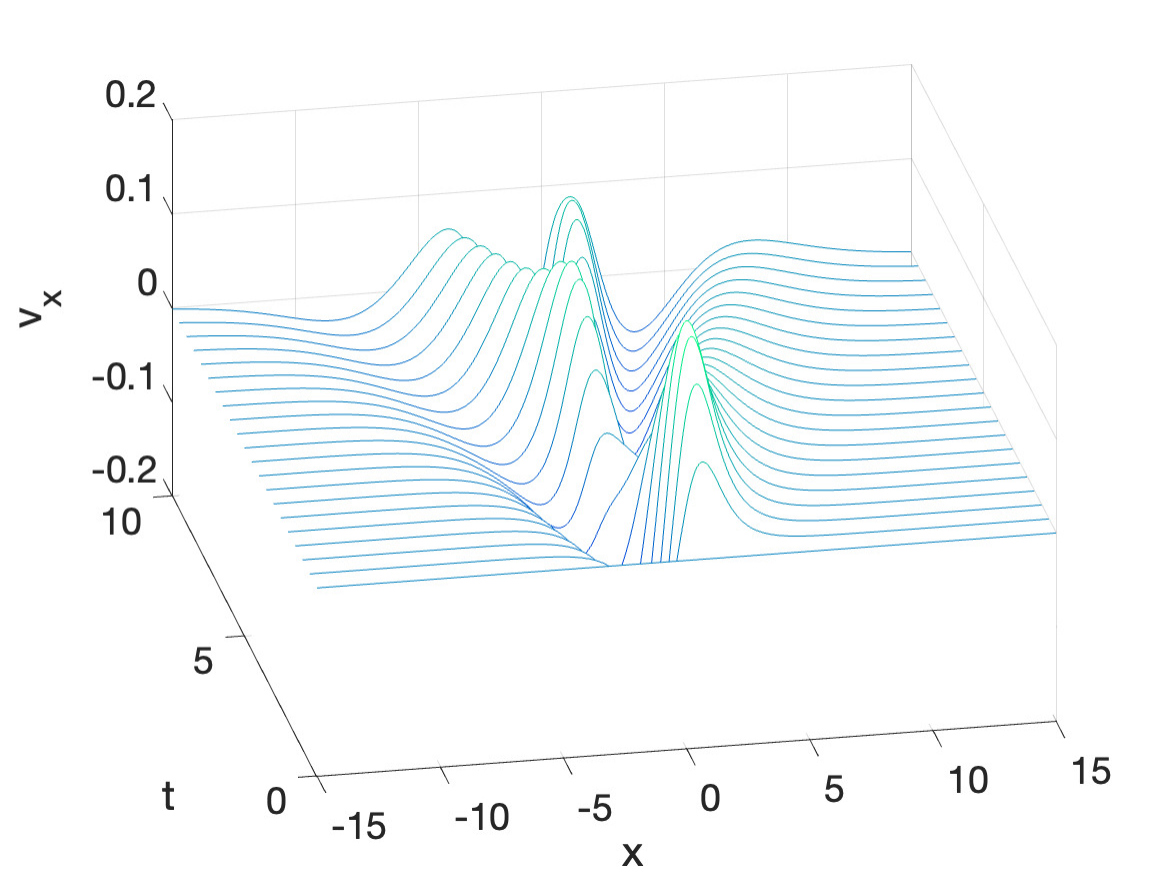} 
  \caption{Solution to the AS system (\ref{2D}) for the initial data 
 (\ref{NCini}) with $\kappa=\alpha=1$ on the $x$-axis in dependence of time,  on 
 the left $\eta$,   on the right  $v_{x}$.}
  \label{figgaussaxis}
 \end{figure}

This behavior of the solution can as before be best observed on the 
$x-$axis since the solution will be once more radially symmetric, see 
Fig.~\ref{figgaussaxis}. The initial hump for $\eta$ on the axis spreads 
overall, but several oscillations between annular structures and 
peaks can be observed. The corresponding plot for $v_{x}$ on the 
right of the same figure shows a similar behavior starting with the 
trivial initial condition. 
 
\begin{figure}[htb!]
  \includegraphics[width=0.49\textwidth]{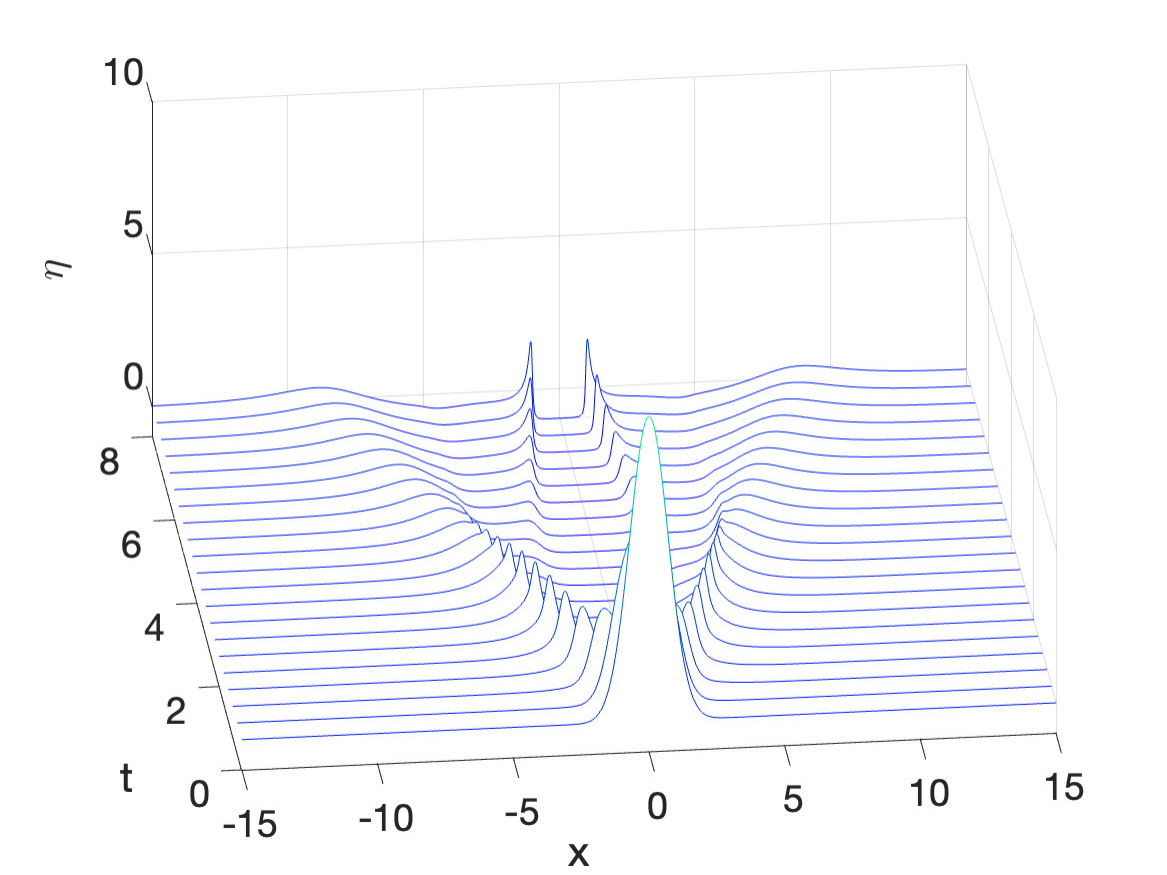}
  \includegraphics[width=0.49\textwidth]{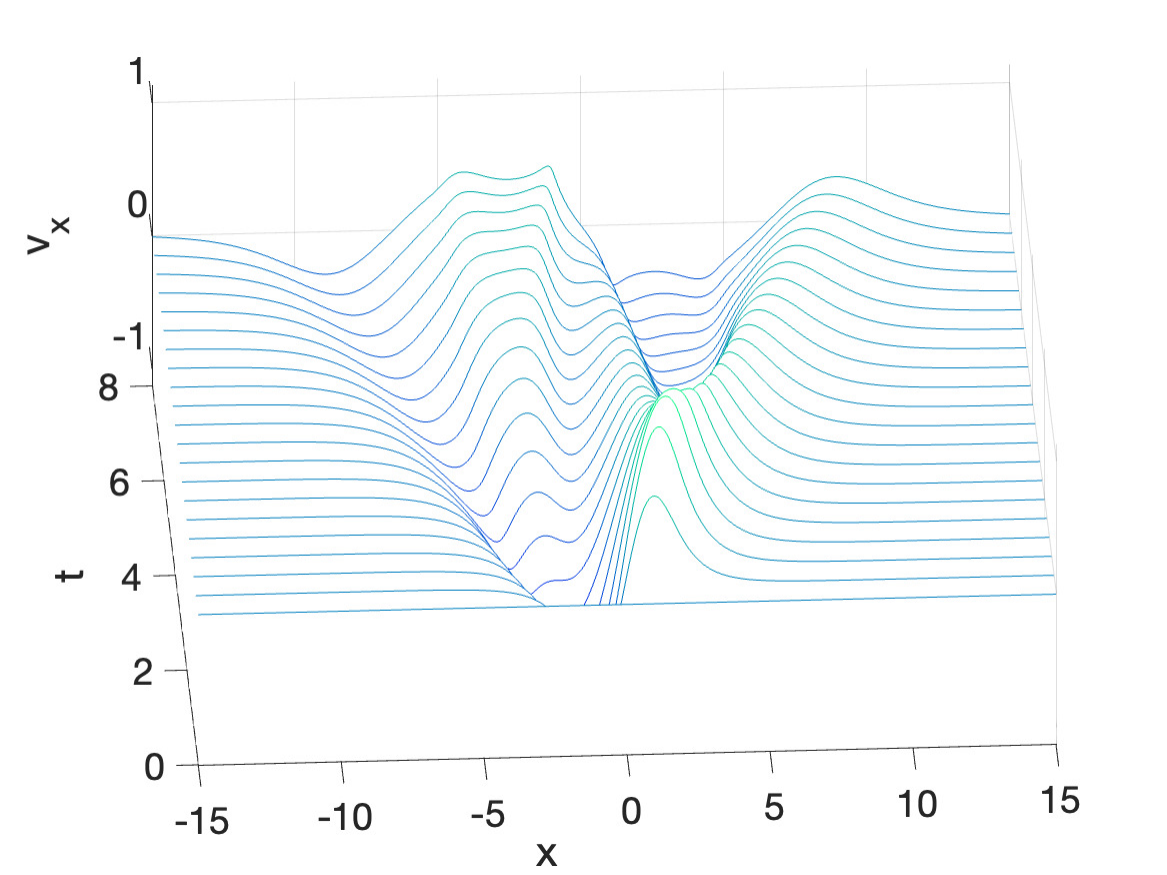} \\
    \includegraphics[width=0.49\textwidth]{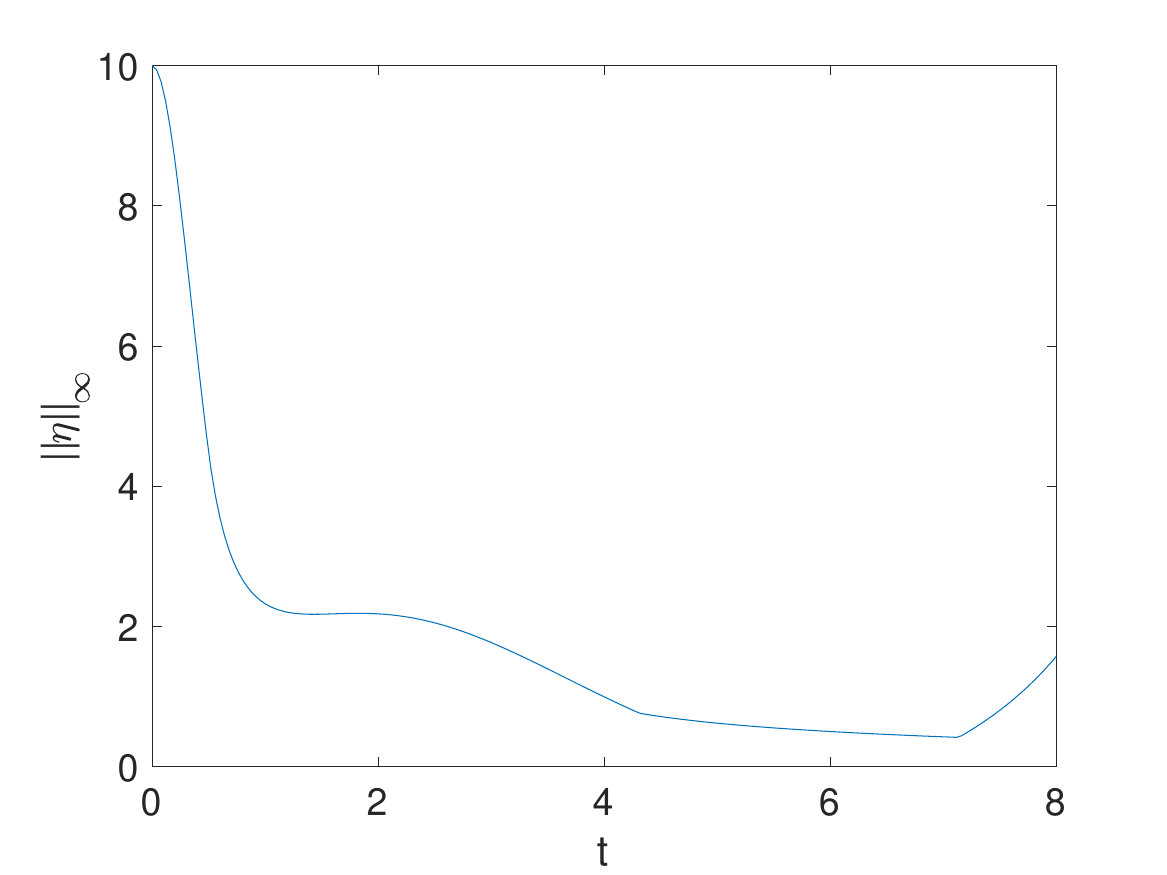}
  \includegraphics[width=0.49\textwidth]{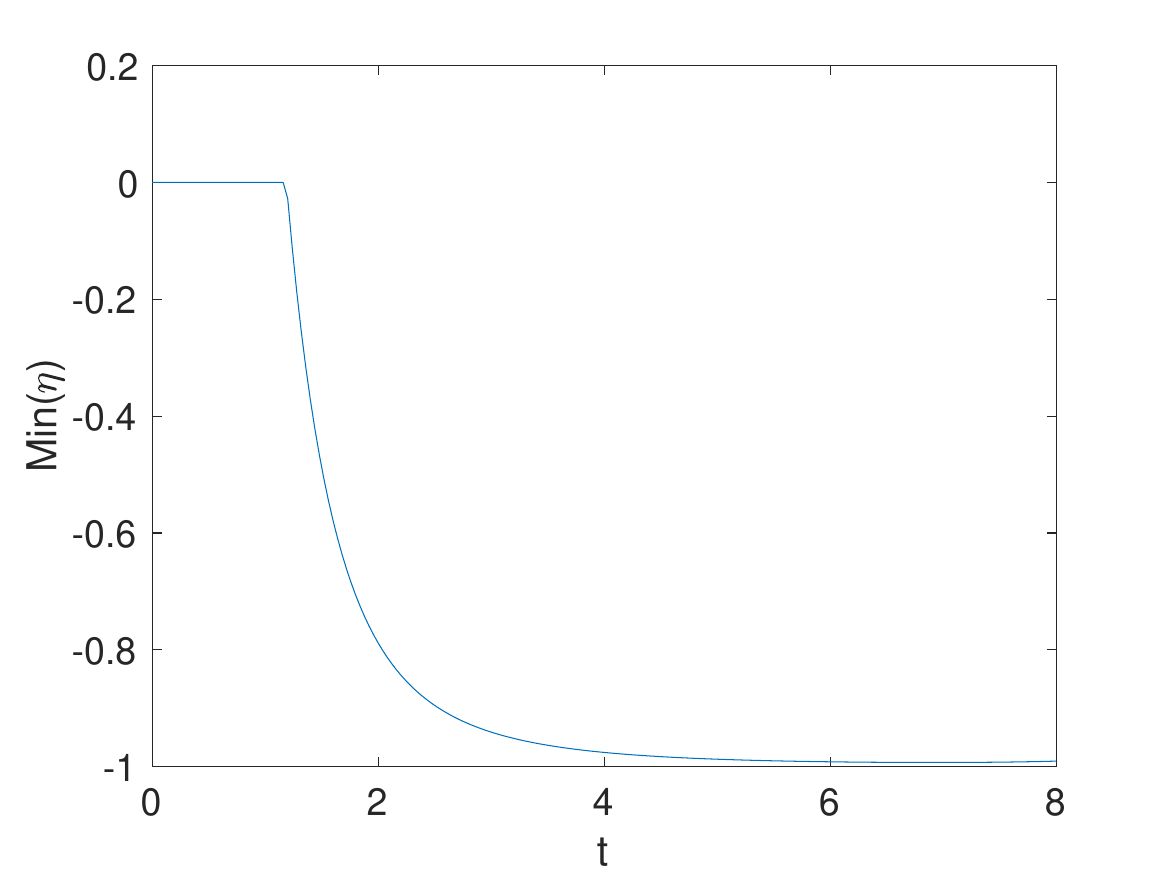} 
  \caption{Solution to the AS system (\ref{2D}) for the initial data 
 (\ref{NCini}) with $\kappa=10$, $\alpha=1$ on the $x$-axis in 
 dependence of time in the upper row,  on 
 the left $\eta$,   on the right  $v_{x}$, in the lower row the 
 $L^{\infty}$ norm of the solution on the left and the global minimum 
 on the right.}
  \label{fig10gaussaxis}
 \end{figure}

The situation becomes more agitated if initial data of larger  norm 
are considered, for instance $\kappa=10$ in (\ref{NCini}). The time 
evolution of the potential on the $x-$axis can be seen in 
Fig.~\ref{fig10gaussaxis}. The initial peak leads to an annular 
structure which develops a hole that approaches the ground to the 
order of $-0.993$, see the lower row of the same figure on the right. 
It is at this point that the dynamics of the last 
section for near cavitation initial data takes over which is why we 
stop the time evolution shortly afterwards. As can be seen on the right of the lower row 
of the figure for the $L^{\infty}$ norm, the near cavitation leads to 
the formation of a peak as discussed in the previous section. 

\begin{figure}[htb!]
  \includegraphics[width=0.49\textwidth]{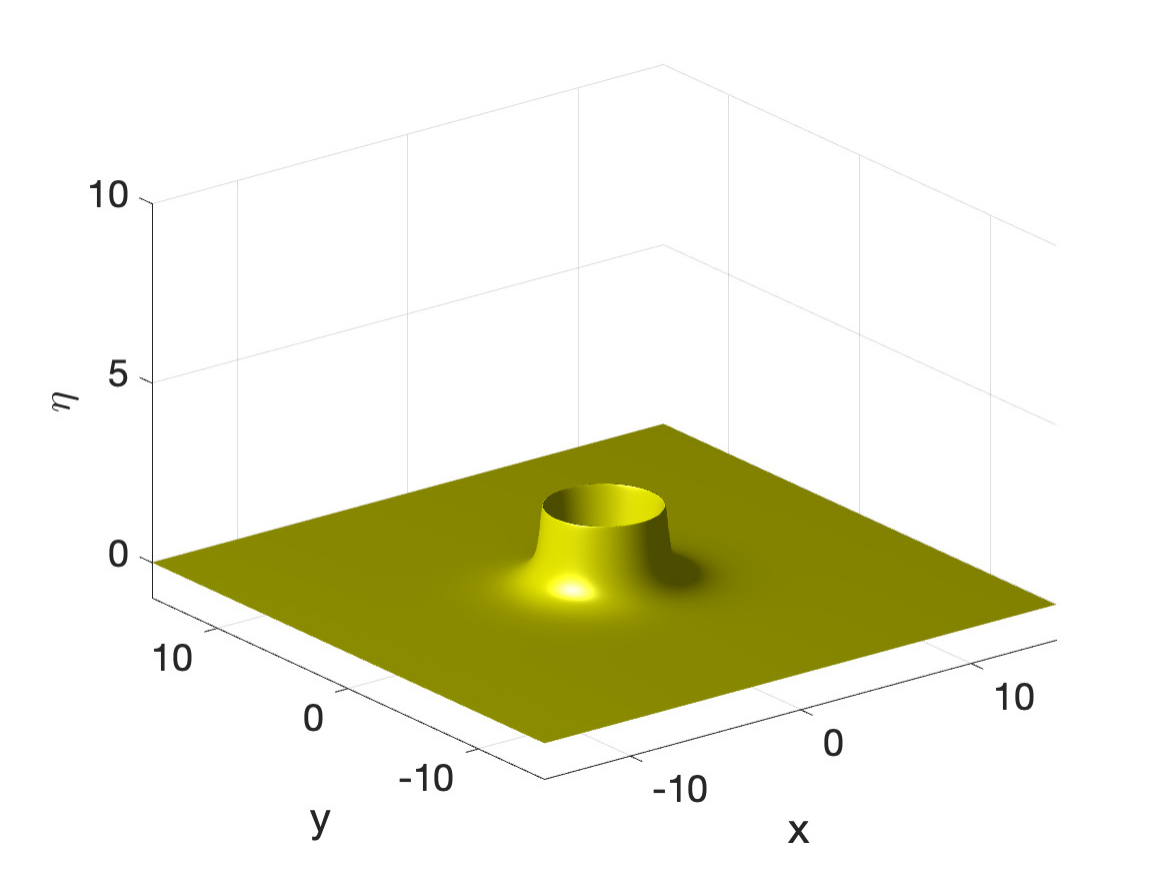}
  \includegraphics[width=0.49\textwidth]{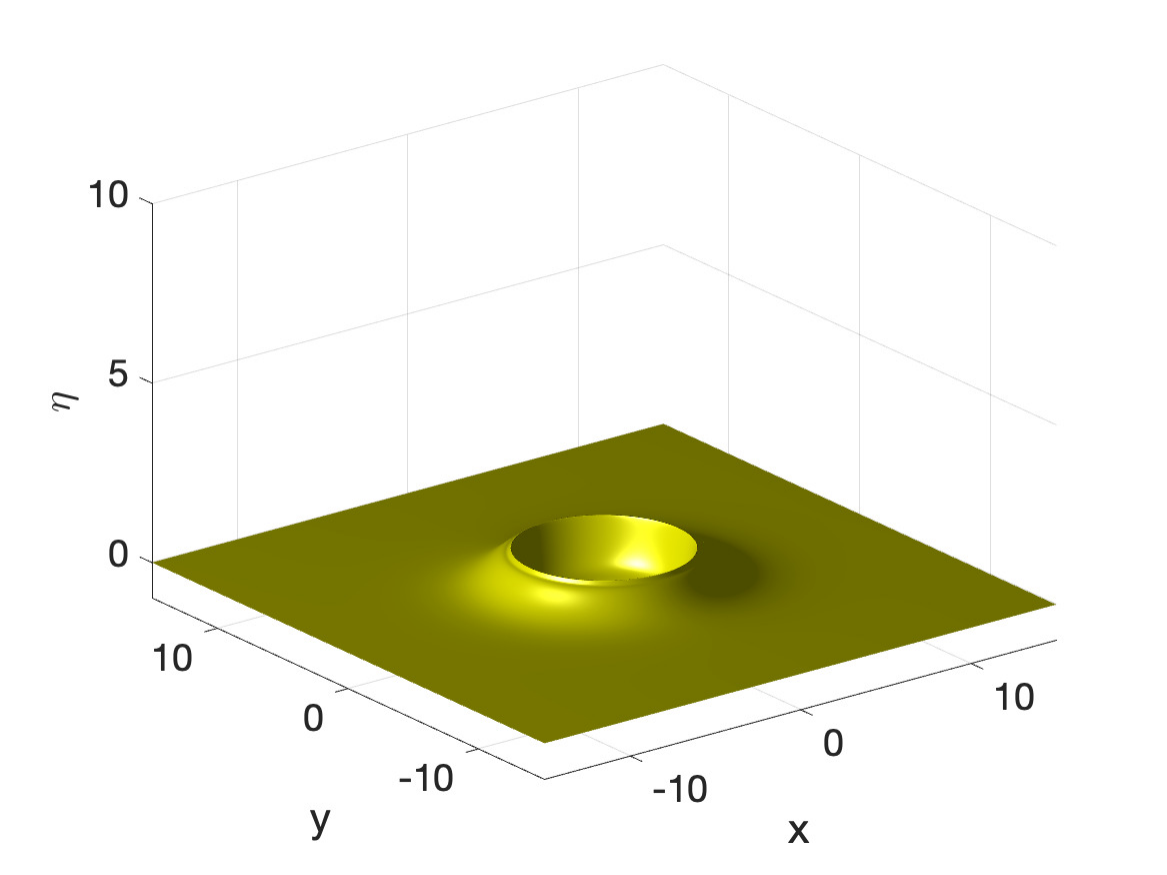} \\
    \includegraphics[width=0.49\textwidth]{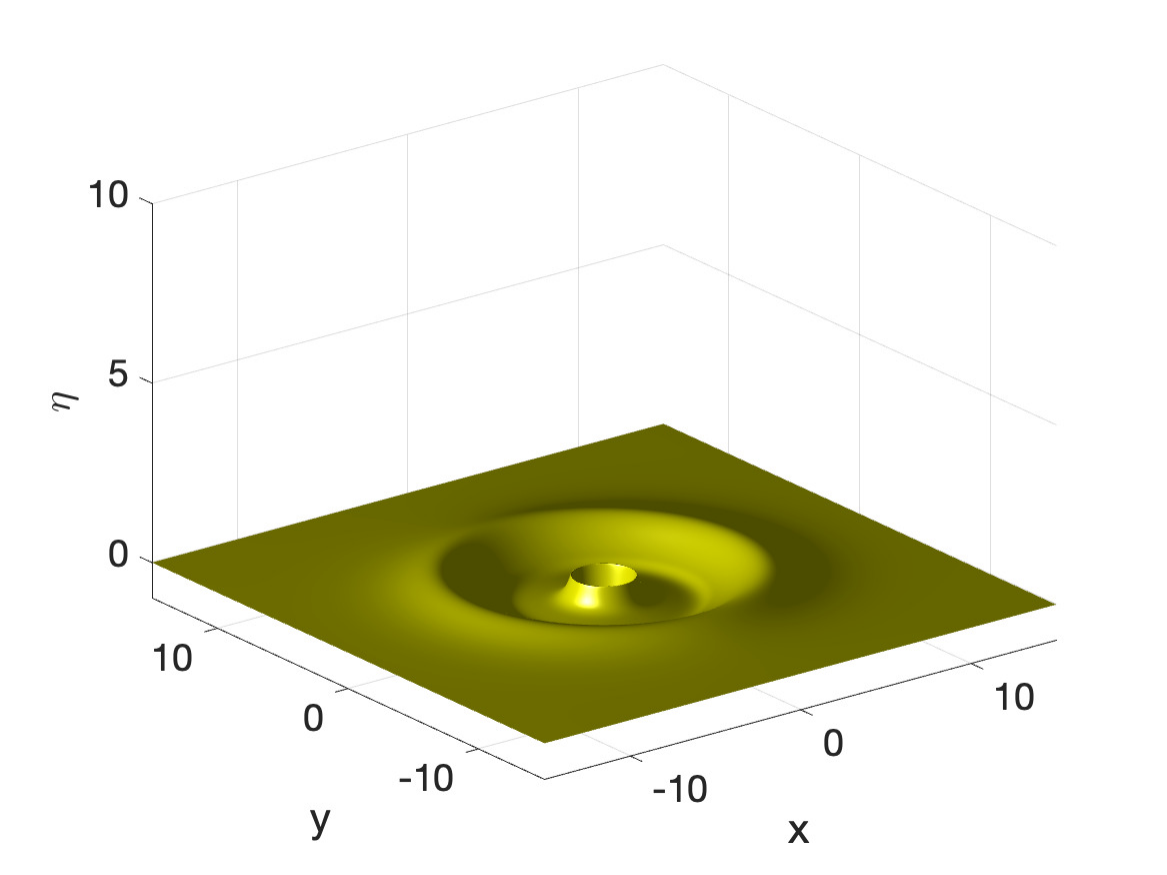}
  \includegraphics[width=0.49\textwidth]{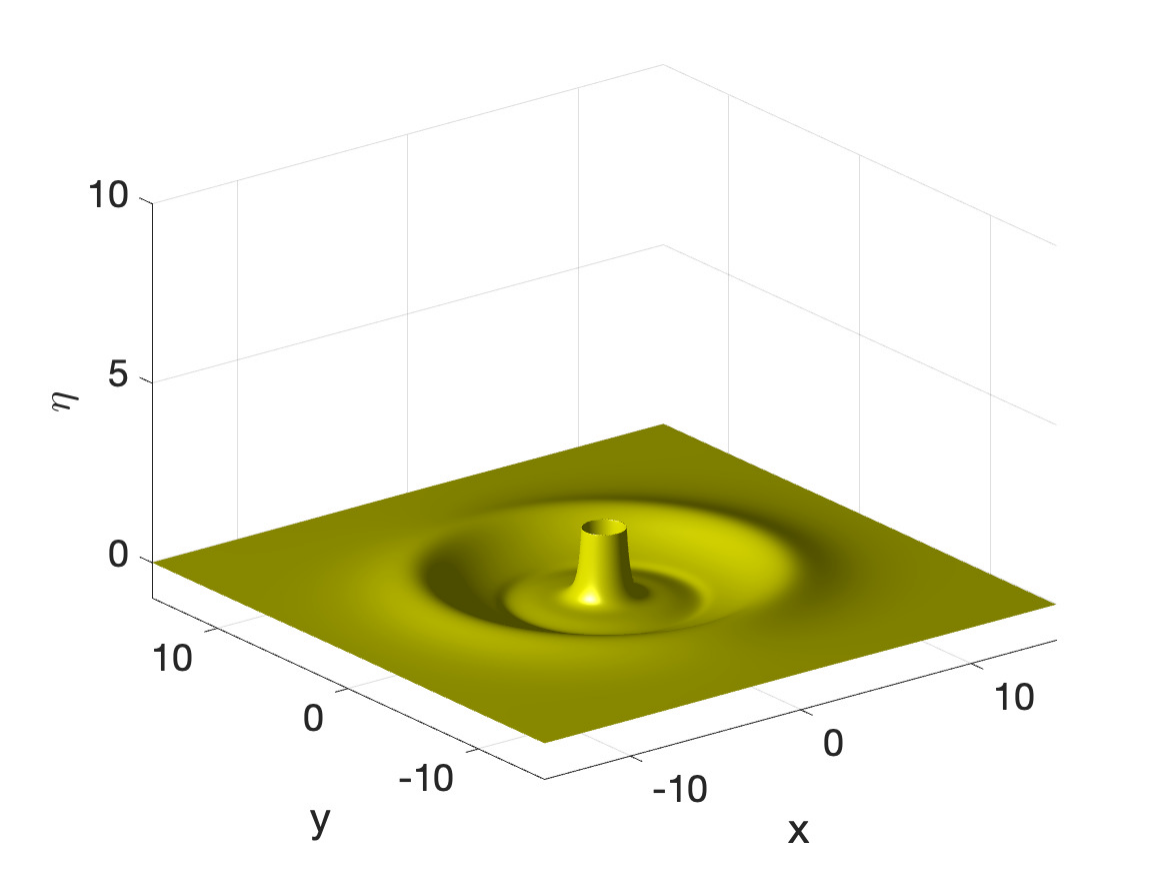} 
  \caption{Solution $\eta$ to the AS system (\ref{2D}) for the initial data 
 (\ref{NCini}) with $\kappa=10$, $\alpha=1$ for several values  of 
 time, in the upper row for $t=2$ and $t=4$, in the lower row for 
 $t=6.92$ (time close to the minimal value) and $t=8$.}
  \label{fig10gausseta}
 \end{figure}

To get a better understanding of the solution in the 2D setting, 
i.e., ignoring the radial symmetry, we show the solution $\eta$ 
at several times in Fig.~\ref{fig10gausseta}.

\begin{figure}[htb!]
  \includegraphics[width=0.49\textwidth]{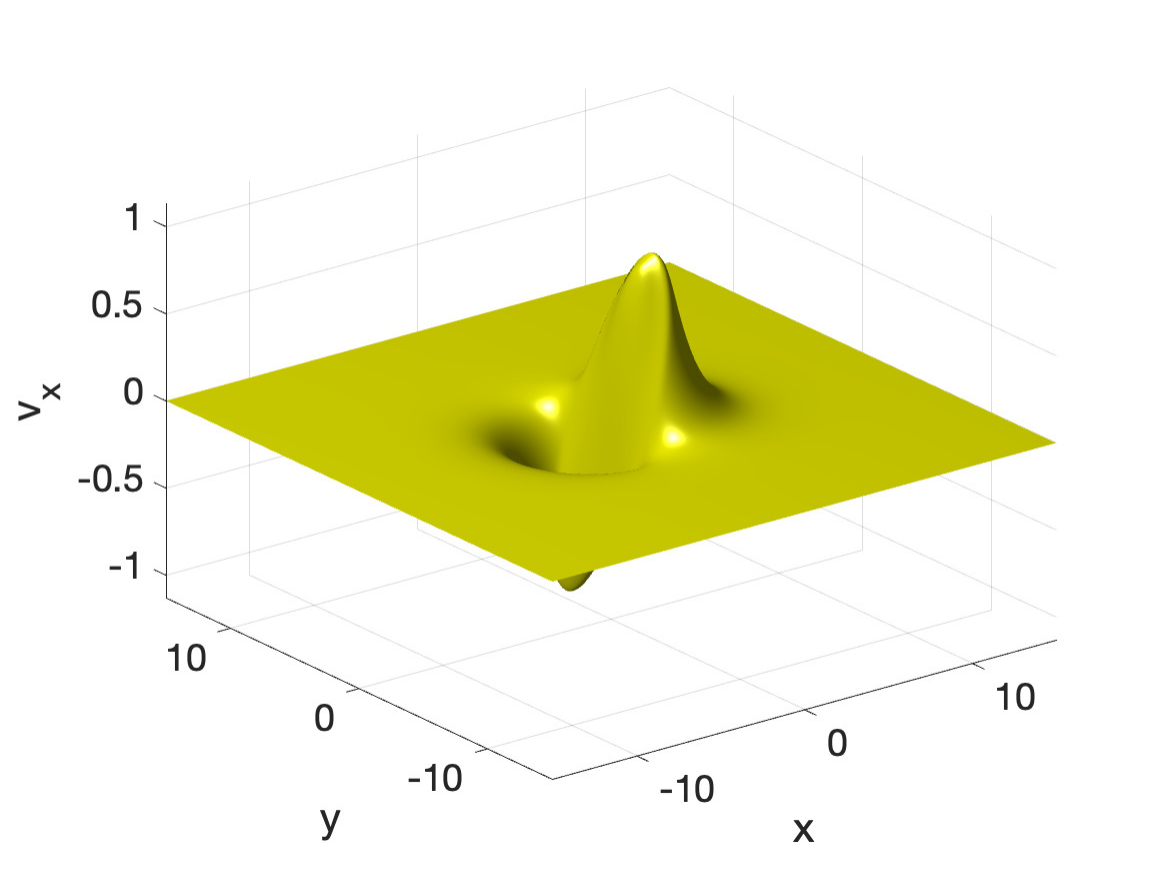}
  \includegraphics[width=0.49\textwidth]{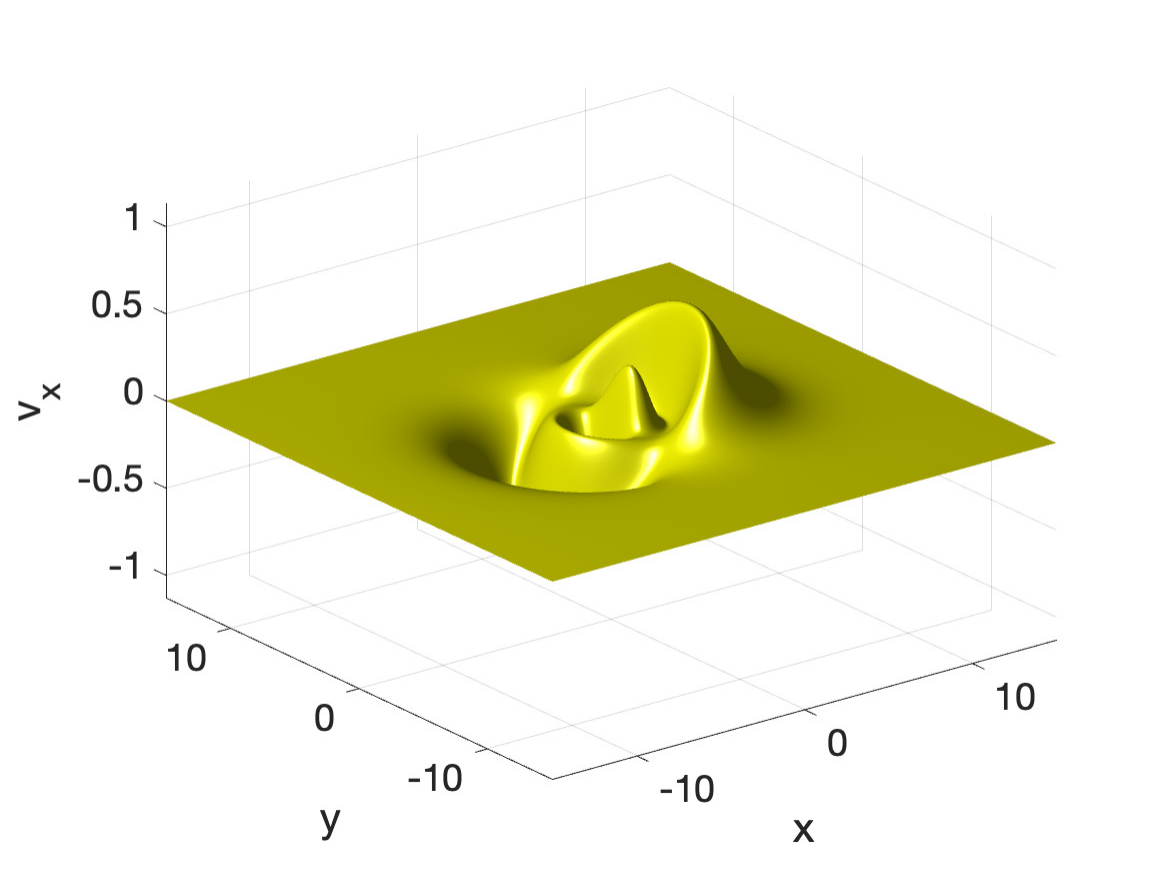} \\
    \includegraphics[width=0.49\textwidth]{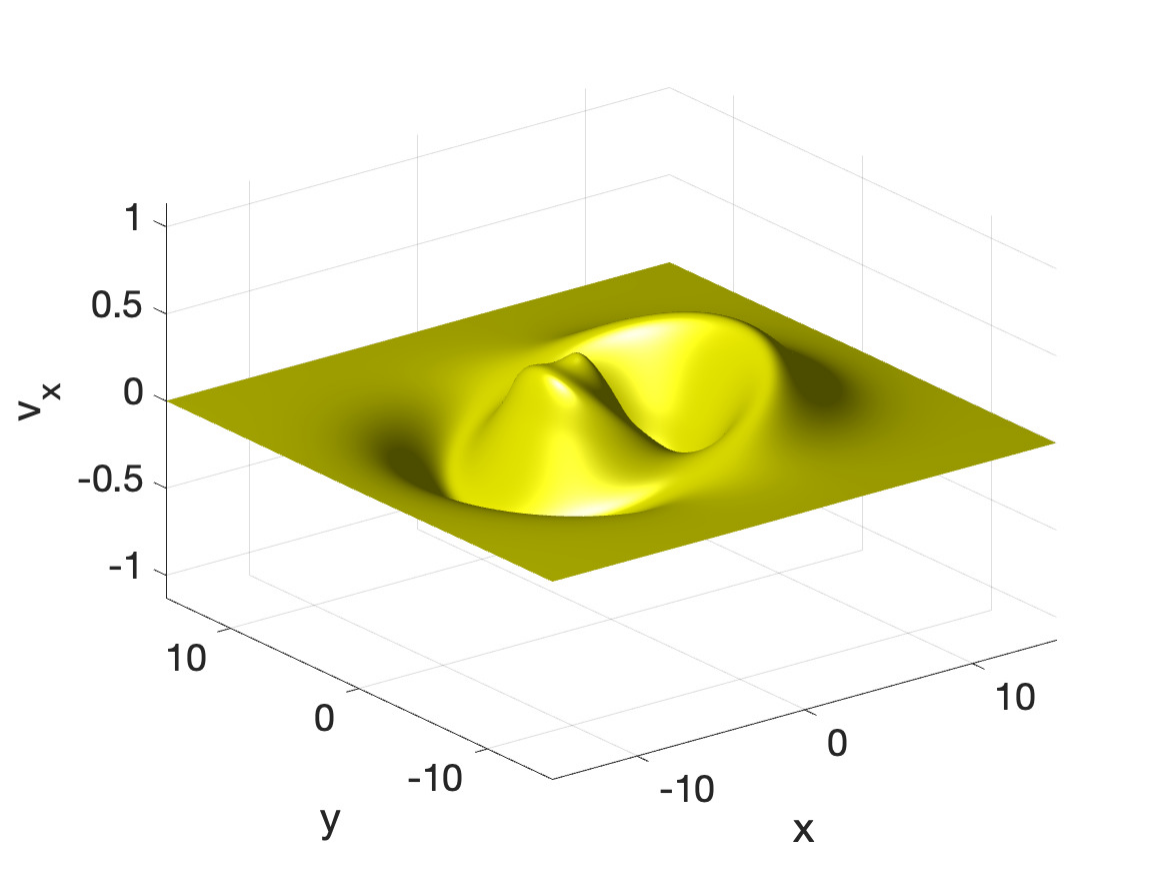}
  \includegraphics[width=0.49\textwidth]{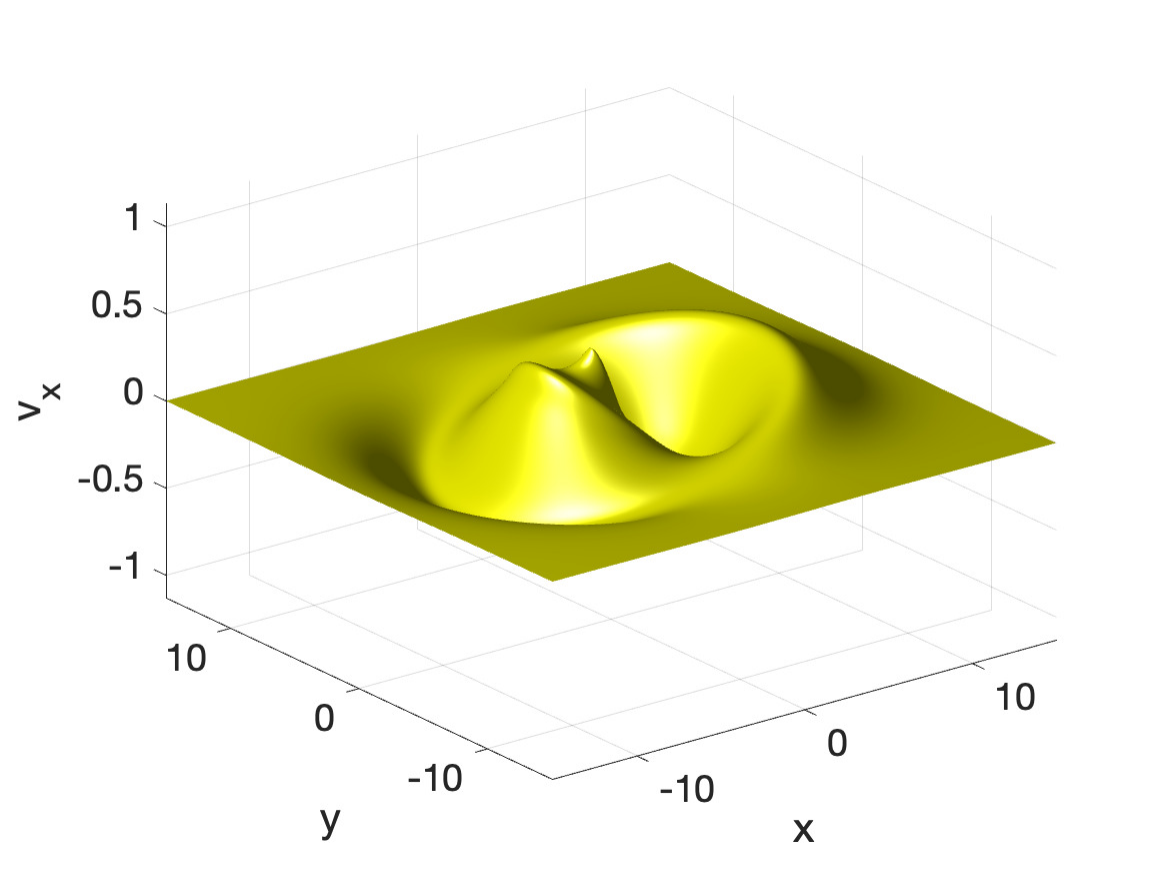} 
  \caption{Solution $v_{x}$ to the AS system (\ref{2D}) for the initial data 
 (\ref{NCini}) with $\kappa=10$, $\alpha=1$ for several values  of 
 time, in the upper row for $t=2$ and $t=4$, in the lower row for 
 $t=6.92$  and $t=8$.}
  \label{fig10gaussvx}
 \end{figure}

The corresponding plots for $v_{x}$ can be seen in 
Fig.~\ref{fig10gaussvx}. 

\section{Dispersive shock waves}
It is well known that the solutions to nonlinear dispersive PDEs can 
have zones of rapid modulated oscillations near shocks of solutions 
for the same initial data to the corresponding dispersionless 
equations. These zones are called \emph{dispersive shock waves} 
(DSWs). A convenient way to study DSWs is to consider initial data 
slowly varying on length scales of order $1/\varepsilon$ where 
$\varepsilon\ll1$ for times of order $1/\varepsilon$. This can be 
addressed by a rescaling of both time and spatial variable with 
$\varepsilon$. In an abuse of notation (we keep the same symbols as 
before in order to avoid cluttered notation) we get for (\ref{2D})
\begin{equation}
    \label{2De}
    \left\lbrace
    \begin{array}{l}
    \eta_t+\nabla \cdot {\bf v}+ \lbrack\nabla\cdot(\eta {\bf v})=0 \\
    {\bf v}_t+\nabla \eta+\lbrack \frac{1}{2}\nabla |{\bf v}|^2-\varepsilon\Delta {\bf v }_t\rbrack=0.
\end{array}\right.
    \end{equation}
The formal limit $\varepsilon\to0$ leads to a Saint-Venant system 
expected to have shock solutions,
\begin{equation}
    \label{2DSV}
    \left\lbrace
    \begin{array}{l}
    \eta_t+\nabla \cdot {\bf v}+ \lbrack\nabla\cdot(\eta {\bf v})=0 \\
    {\bf v}_t+\nabla \eta+\lbrack \frac{1}{2}\nabla |{\bf v}|^2\rbrack=0.
\end{array}\right.
    \end{equation}

\begin{figure}[htb!]
  \includegraphics[width=0.49\textwidth]{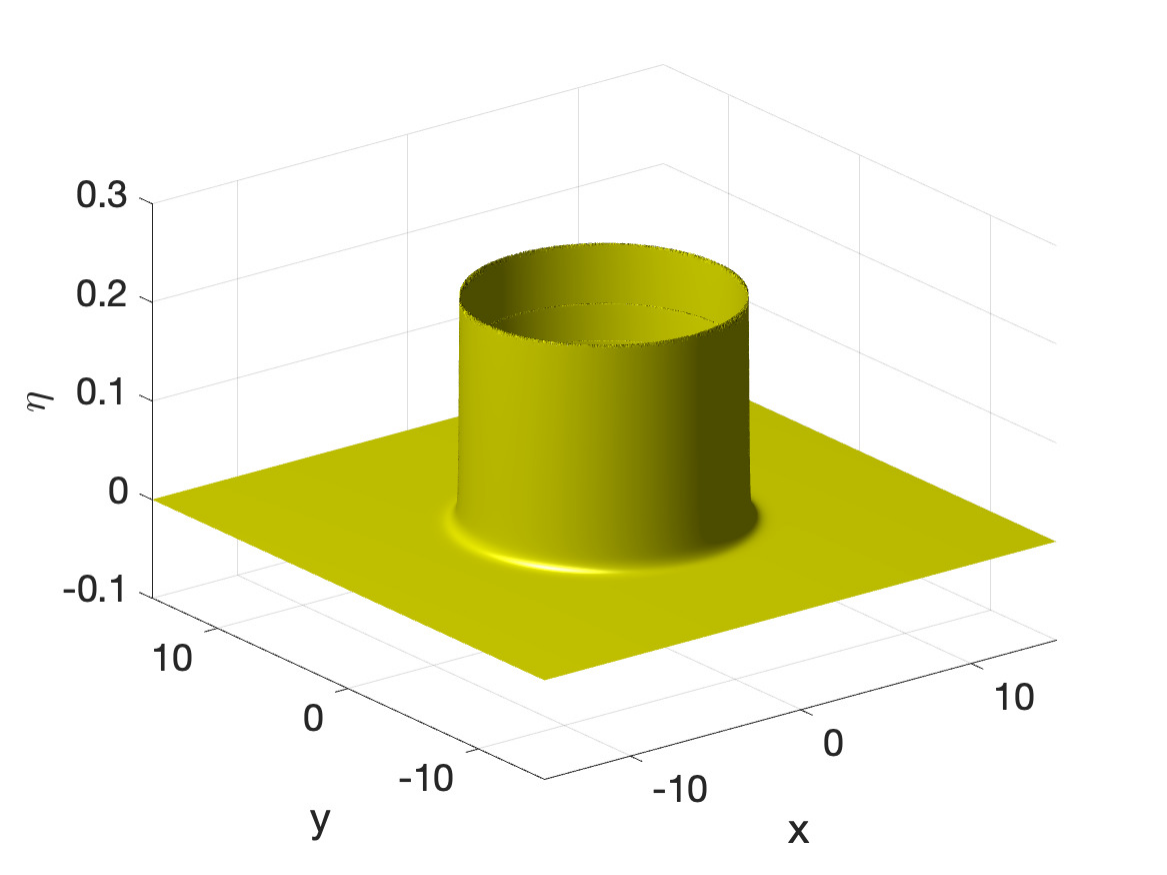}
  \includegraphics[width=0.49\textwidth]{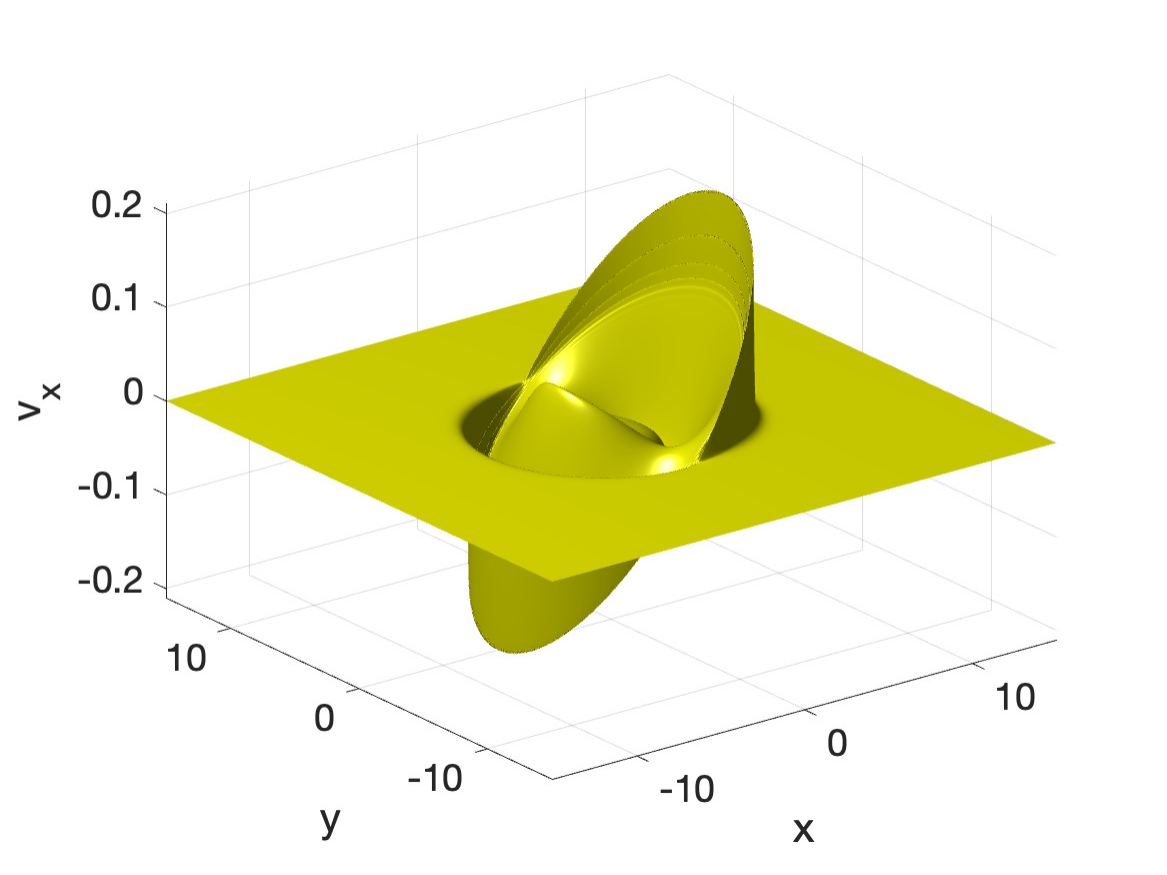} 
  \caption{Solution to the AS system (\ref{2De}) for the initial data 
 (\ref{NCini}) with $\kappa=\alpha=1$ and $\varepsilon=10^{-2}$ for 
 $t=5$, on the left  $\eta$,   on the right  $v_{x}$.}
  \label{figgausse2}
 \end{figure}

To illustrate the behavior of the solutions for small values of 
$\varepsilon$, we consider once more Gaussian initial data for 
$\varepsilon=10^{-2}$. The solution for $t=5$ is shown in 
Fig.~\ref{figgausse2}. It can be seen that the time evolution is 
similar to what was shown in the previous section, an initial peak 
develops into some annular structure for $\eta$. But the small value 
$\varepsilon$ leads to very strong gradients in an almost tubular 
structure, and similarly for $v_{x}$. 
 
\begin{figure}[htb!]
  \includegraphics[width=0.49\textwidth]{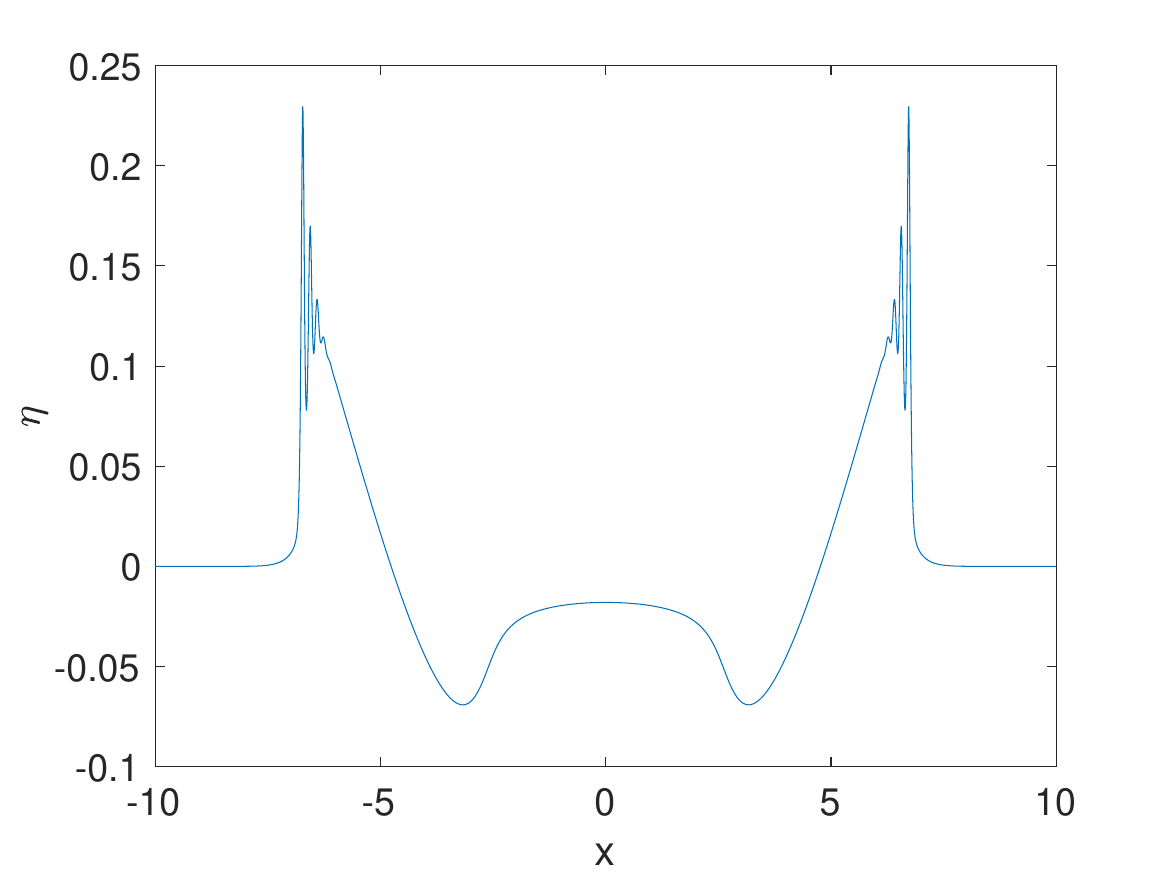}
  \includegraphics[width=0.49\textwidth]{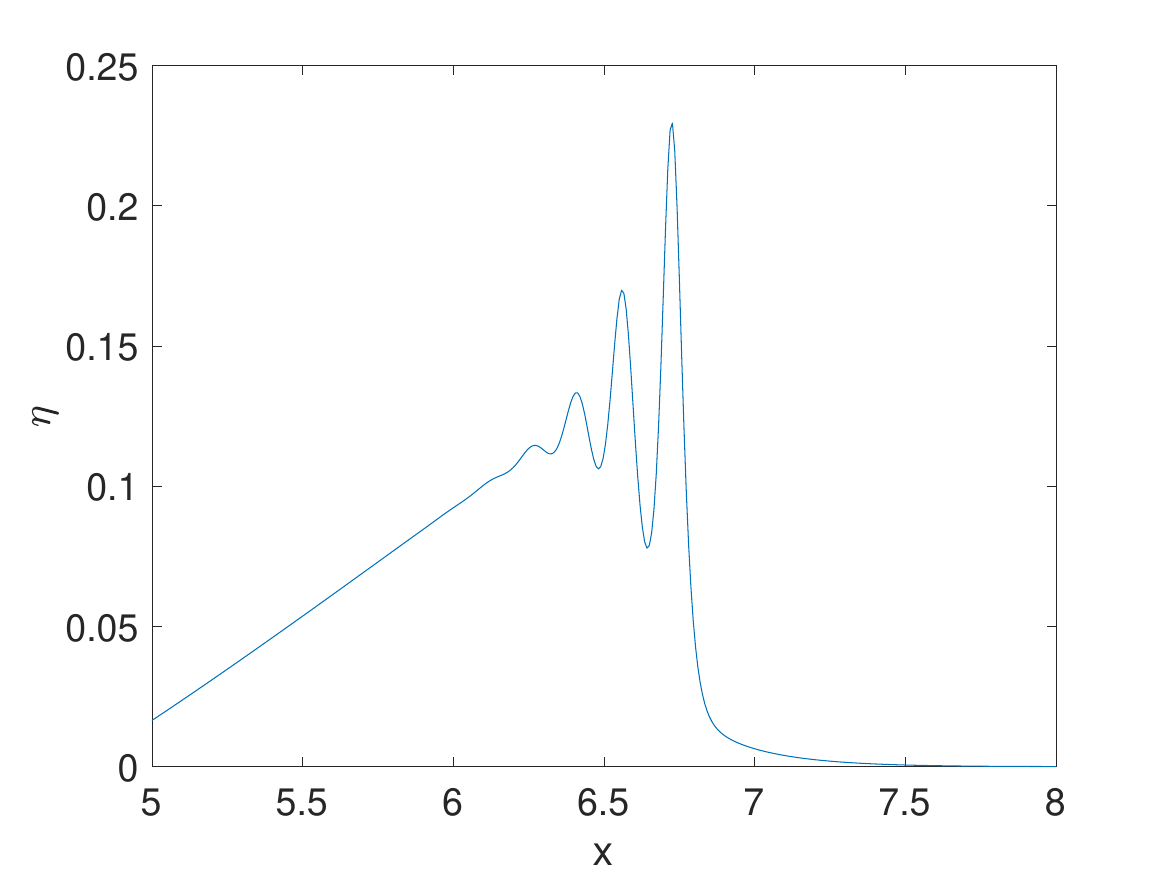} \\
    \includegraphics[width=0.49\textwidth]{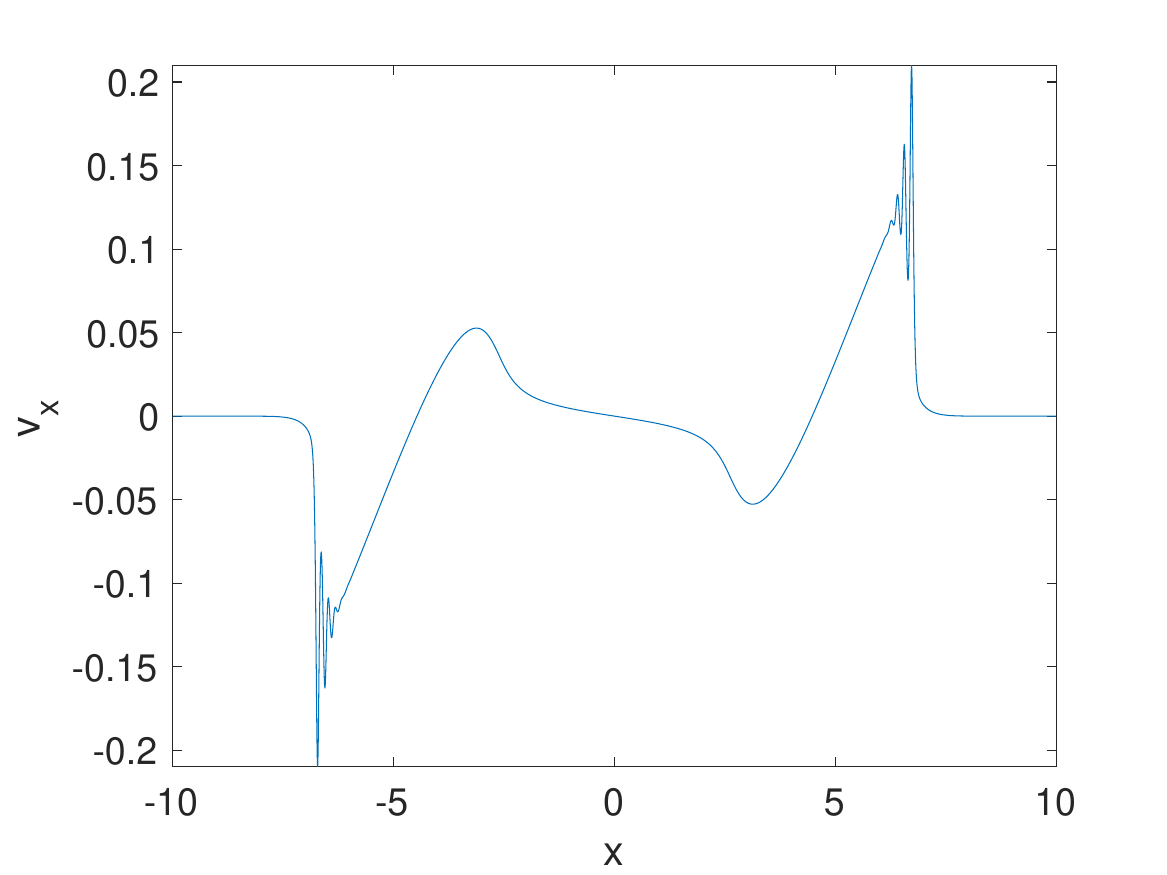}
  \includegraphics[width=0.49\textwidth]{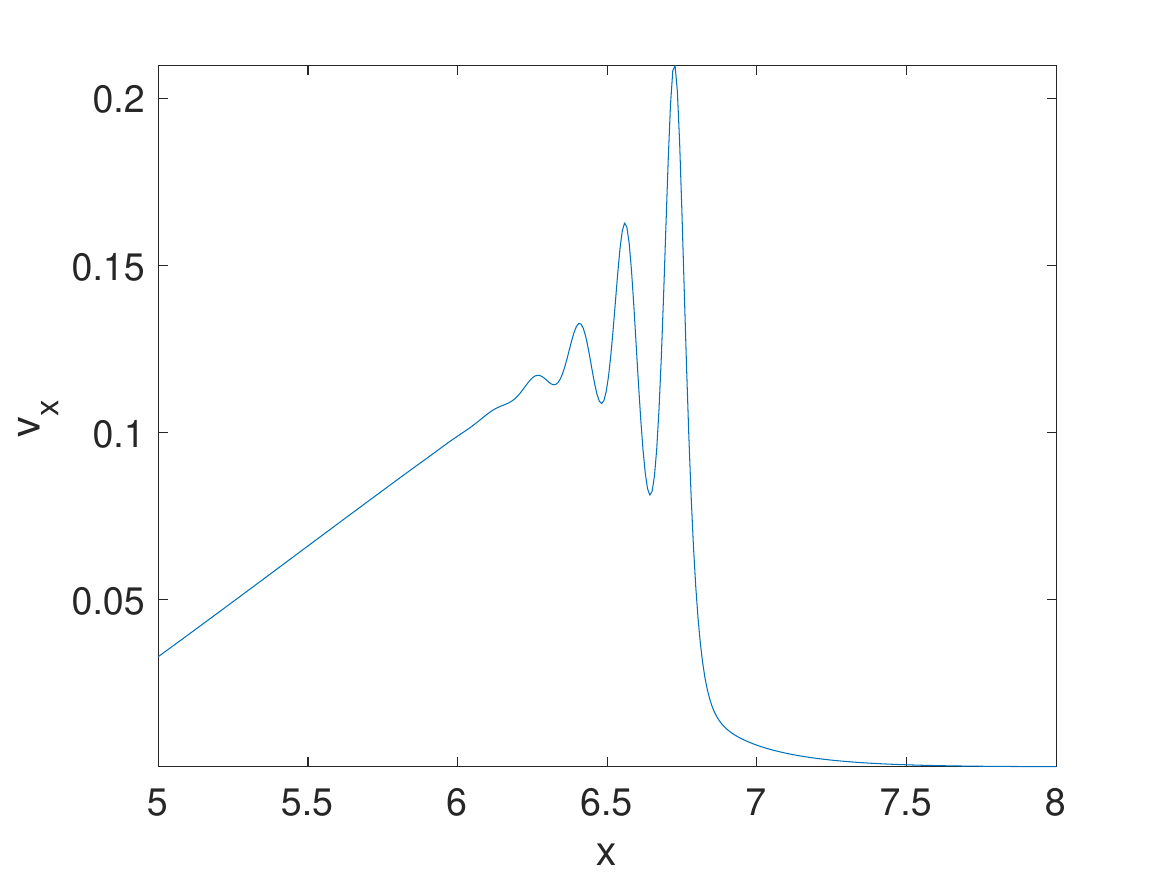} 
  \caption{Solution to the AS system (\ref{2De}) for the initial data 
 (\ref{NCini}) with $\kappa=\alpha=1$ and $\varepsilon=10^{-2}$ for 
 $t=5$ on the $x$-axis, in the upper row $\eta$, in the lower $v_{x}$, 
 both with a close-up of the oscillatory zone on the right.}
  \label{figgausse2zoom}
 \end{figure}

 However there are small oscillation on what appears to be a steep 
 front in Fig.~\ref{figgausse2}. To illustrate this, we show the 
 potentials in Fig.~\ref{figgausse2} on the axis together with a 
 close-up of the oscillatory zone in Fig.~\ref{figgausse2zoom}. It 
 can be  seen that there are modulated oscillations near the strong 
 gradients. 

\section{Conclusion}
In this paper we have presented a numerical study of the 2D 
Amick-Schonbek system (\ref{2D}). Numerical evidence for the 
transverse stability of line solitary waves was presented as well as 
for the absence of stable lump solitary waves, i.e., for solitary 
waves localised in both spatial dimensions. It was shown that 
solutions to the system can exhibit zones of rapid modulated 
oscillations known as dispersive shock waves. It was argued that 
initial data not satisfying the non-cavitation condition lead to the 
formation of cusps in finite time. A precise characterization of 
these cusp was not possible due to a local of numerical resolution. 
This will be the subject of further research. 

Similar numerical studies will be performed for the Boussinesq system 
integrable in 1D and numerically studied in \cite{KSKBK}. The results 
are to be compared to other equations appearing in the context of 
water waves as the Serre-Green-Naghdi equations \cite{GN1,GN2,Ser}. 
The goal is to expand the study of \cite{GK} to the questions 
addressed in the present paper. This will be the subject of future 
work.  

\begin{merci}  
 The work of the first author  was partially  supported 
by the ANR project ANR-17-EURE-0002 EIPHI. 

Both authors were partially supported by the ANR project ISAAC-ANR-23-CE40-0015-01.
\end{merci}

\bibliographystyle{amsplain}

\end{document}